\documentclass[10pt]{article}
\usepackage{color}
\usepackage{amsmath}
\usepackage{amssymb}
\usepackage{mathrsfs}
\usepackage[utf8]{inputenc}
\usepackage[english]{babel}
\usepackage{amsthm}
\usepackage{tikz}
\usepackage{mathdots}

\usepackage{graphicx}
\graphicspath{{./images/}}

\usepackage[toc,title,page]{appendix}

\usepackage{epsfig}

\usepackage[scaled=1]{helvet}
\usepackage{titlesec}
\usepackage{multido}
\usepackage{multirow}
\usepackage{blindtext}
\usepackage{color}
\usepackage{lipsum}
\usepackage{mathtools}
\usepackage[vcentermath]{youngtab}
\usepackage{young}
\usepackage[all,ps,dvips,graph]{xy}
\usepackage[misc,geometry]{ifsym}

\usepackage[left=3cm,top=2cm,right=3cm,bottom=2cm]{geometry}
\usepackage{graphicx}
\numberwithin{equation}{subsection}

\let\OLDthebibliography\thebibliography
\renewcommand\thebibliography[1]{
  \OLDthebibliography{#1}
  \setlength{\parskip}{0pt}
  \setlength{\itemsep}{0pt plus 0.3ex}
}

\newtheorem{theorem}{Theorem}
\newtheorem{proposition}[theorem]{Proposition}
\newtheorem{corollary}[theorem]{Corollary}
\newtheorem{lemma}[theorem]{Lemma}
\newtheorem{remark}[theorem]{Remark}
\newtheorem{example}[theorem]{Example}
\newtheorem{definition}[theorem]{Definition}

\numberwithin{theorem}{section}

\newcommand{\A}{\mathcal{A}}
\newcommand{\aaa}{\mathfrak{a}}

\newcommand{\Basis}{\mathcal{C}_n}
\newcommand{\BasisOne}{\mathcal{C}^{1,\K}_n}
\newcommand{\B}{\mathcal{B}_{l,m}}
\newcommand{\BB}{\mathcal{B}^{\F}_{l,m}(\kappan,\een)}
\newcommand{\BBnu}{\mathcal{B}^{\F}_{l,m}(\bnu,\een)}
\newcommand{\Bb}{\pmb{\mathfrak{b}}}
\newcommand{\BBB}{\mathcal{B}_{n+1}}
\newcommand{\bbb}{\mathfrak{b}}
\newcommand{\BBk}{\mathcal{B}_k}
\newcommand{\balpha}{\boldsymbol{\alpha}}
\newcommand{\bbeta}{\boldsymbol{\beta}}
\newcommand{\bgamma}{\boldsymbol{\gamma}}
\newcommand{\bGamma}{\boldsymbol{\Gamma}}
\newcommand{\bbkap}{\boldsymbol{\kappa}}
\newcommand\bbn{\mathbf{b}}
\newcommand\bbs{\mathsf{s}}
\newcommand\bbt{\boldsymbol{\mathsf{t}}}
\newcommand\bbu{\mathsf{u}}
\newcommand\bbv{\mathsf{v}}
\newcommand{\Bchico}{\mathcal{B}_{n-1}}
\newcommand\be{\mathbb{E}}
\newcommand\ben{\boldsymbol{e}}
\newcommand{\belongs}{ \in }
\newcommand{\Belongs}{ \ni}
\newcommand{\Bg}{\pmb{\mathfrak{g}}}
\newcommand\bi{\boldsymbol{i}}
\newcommand{\bif}{{\underline{\boldsymbol{i}}}}
\newcommand{\bjf}{{\underline{\boldsymbol{j}}}}
\newcommand\bj{\boldsymbol{j}}
\newcommand\blambda{{\boldsymbol\lambda}}
\newcommand\bmu{{\boldsymbol\mu}}
\newcommand\brho{\boldsymbol{\rho}}
\newcommand\bzeta{{\boldsymbol\zeta}}
\newcommand\bn{\boldsymbol{n}}
\newcommand\bnu{{\boldsymbol\nu}}
\newcommand\boxbluek{\color{blue}\boldsymbol{[k]}}
\newcommand\boxbluej{\color{blue}\boldsymbol{[j]}}
\newcommand\boxredk{\color{red}\boldsymbol{[k]}}
\newcommand\boxredj{\color{blue}\boldsymbol{[j]}}
\newcommand{\bpartial}{\boldsymbol{\partial}}
\newcommand{\bCpartial}{\mathbf{C}\boldsymbol{\partial}}
\newcommand{\Bs}{\pmb{\mathfrak{s}}}
\newcommand\bS{\Sigma}
\newcommand\bs{\mathbf{s}}
\newcommand{\bT}{\pmb{\mathfrak{t}}}
\newcommand\bt{\mathbf{t}}
\newcommand{\bTI}{ \bT^{-1} (\bT_{\theta}^{\blambda}(1))}
\newcommand{\bTII}{ \bT^{-1} (\bT_{\theta}^{\blambda}(2))}
\newcommand{\bTj}{ \bT^{-1} (\bT_{\theta}^{\blambda}(j))}
\newcommand{\bTn}{ \bT^{-1} (\bT_{\theta}^{\blambda}(n))}
\newcommand\btau{{\boldsymbol\tau}}
\newcommand{\Bu}{\pmb{\mathfrak{u}}}
\newcommand{\Bv}{\pmb{\mathfrak{v}}}
\newcommand\bu{\mathbf{u}}
\newcommand\bv{\mathbf{v}}
\newcommand\bcalJ{\boldsymbol{\mathcal{J}}}
\newcommand\bcalL{\boldsymbol{\mathcal{L}}}
\newcommand\bcalm{\boldsymbol{\mathcal{m}}}

\newcommand\bcalNT{\boldsymbol{\mathcal{NT}}}
\newcommand\bcalNGAT{\boldsymbol{\mathcal{NT}}}
\newcommand\bcalD{\boldsymbol{\mathcal{D}}}
\newcommand\bcalR{\boldsymbol{\mathcal{R}}}
\newcommand\bcalV{\boldsymbol{\mathcal{V}}}

\newcommand\bcalU{\boldsymbol{\mathcal{U}}}
\newcommand\bcalY{\boldsymbol{\mathcal{Y}}}
\newcommand\bcalW{\boldsymbol{\mathcal{W}}}
\newcommand\bnabla{\boldsymbol{\nabla}}

\newcommand{\C}{\mathcal{C}_{\rm{YH}}}
\newcommand\calA{\mathcal{A}}
\newcommand\calB{\mathcal{B}}
\newcommand\calC{\mathcal{C}}
\newcommand\calD{\mathcal{D}}
\newcommand\calE{\mathcal{E}}
\newcommand\calF{\mathcal{F}}
\newcommand\calG{\mathcal{G}}
\newcommand\calH{\mathcal{H}}
\newcommand\calL{\mathcal{L}}
\newcommand\calM{\mathcal{M}}
\newcommand\calN{\mathcal{N}}

\newcommand\calO{\mathcal{O}}
\newcommand\calP{\mathcal{P}}
\newcommand\calQ{\mathcal{Q}}
\newcommand\calR{\mathcal{R}}
\newcommand\calS{\mathcal{S}}
\newcommand\calT{\mathcal{T}}

\newcommand\calU{\mathcal{U}}
\newcommand\calV{\mathcal{V}}
\newcommand\calW{\mathcal{W}}
\newcommand\calX{\mathcal{X}}
\newcommand\calY{\mathcal{Y}}
\newcommand\calZ{\mathcal{Z}}

\newcommand{\catorce}{ 14}
\newcommand{\catorceB}{\color{red} 14}
\newcommand{\CC}{ \mathbb C }
\newcommand{\ccc}{\mathfrak{c}}
\newcommand{\ch}{{\rm char}}
\newcommand{\cincuentacinco}{55}
\newcommand{\cincuentacincoR}{\color{red}55}
\newcommand{\Comp}{{\mathcal Comp}_n}
\newcommand{\cuarentacuatro}{ 44}
\newcommand{\cupdot}{\mathbin{\mathaccent\cdot\cup}}
\newcommand{\Cs}{\overleftarrow{C}}
\newcommand{\Cd}{\overrightarrow{C}}

\newcommand{\Der}{{\rm Der}}
\newcommand{\diez}{10}
\newcommand{\dieciseis}{16}
\newcommand{\diezR}{\color{red}10}
\newcommand{\doce}{12}
\newcommand{\doceB}{\color{blue} 12}
\newcommand{\doceR}{\color{red} 12}

\newcommand{\E}{ {\mathcal E}_n(q)}
\newcommand{\e}{\mathfrak{e}}
\newcommand{\EE}{ {\mathcal E}_n}
\newcommand\een{\mathbf{e}}
\newcommand\es{\mathbbm{s}}
\newcommand{\End}{{\rm End}}
\newcommand\et{\mathbbm{t}}
\newcommand\eu{\mathbbm{u}}
\newcommand\ev{\mathbbm{v}}
\newcommand{\Exp}{ {\rm \bf exp} }

\newcommand{\F}{ { \mathbb F}}
\newcommand{\FF}{ {\mathcal F}_n}

\newcommand{\g}{  \mathfrak{g}}
\newcommand{\gl}{\mathfrak{gl}}

\newcommand{\h}{{h}}
\newcommand{\HH}{ \mathcal{H}_n}
\newcommand{\HHO}{ \mathcal{H}^{\OO}_n}
\newcommand{\HHK}{ \mathcal{H}^{\K}_n}
\newcommand{\HHtwo}{ \mathcal{H}_2}
\newcommand{\HHKtwo}{ \mathcal{H}^{\K}_2}
\newcommand{\HHOtwo}{ \mathcal{H}^{\OO}_2}
\newcommand{\HHKOne}{ \mathcal{H}^{1,\K}_n}

\newcommand{\II}{I_{\mathbf{e}}}
\newcommand{\IIa}{I_{e}}
\newcommand{\id}{{\rm id}}
\newcommand{\ind}{{\rm ind}}
\newcommand{\inv}{{\rm inv}}

\newcommand{\JM}{ \mathcal L }

\newcommand{\K}{\mathcal{K}}
\newcommand{\kk}{\mathcal{K}}
\newcommand{\kkk}{k-1}
\newcommand{\kappan}{\boldsymbol{\kappa}}

\newcommand{\Li}{\mathcal{L}}
\newcommand{\LL}{\mathbb{L}}

\newcommand{\m}{\mathfrak{m}}
\newcommand{\MC}{{ {\rm Comp}}_{l,n}}
\newcommand{\MCm}{{ {\rm Comp}}_{l,m}}
\newcommand{\MP}{{\rm Par }_{l,n}}
\newcommand{\mfra}{{\mathfrak{a}}}
\newcommand{\mfrb}{{\mathfrak{b}}}
\newcommand{\mfrc}{{\mathfrak{c}}}
\newcommand{\mfrt}{{\mathfrak{t}}}
\newcommand{\mfrs}{{\mathfrak{s}}}
\newcommand{\mfru}{{\mathfrak{u}}}
\newcommand{\mfrv}{{\mathfrak{v}}}
\newcommand{\mfrw}{{\mathfrak{w}}}
\newcommand{\mfrx}{{\mathfrak{x}}}
\newcommand{\mfry}{{\mathfrak{y}}}
\newcommand{\mfrz}{{\mathfrak{z}}}

\newcommand{\N}{ { \mathbb N}}
\newcommand{\No}{ { \mathbb N}_0}
\newcommand{\nstd}{{\rm NStd}}
\newcommand{\nSLTM}{{$n-$\rm SLTM}}
\newcommand{\nSLTMn}{$\mathcal{TM}_n$}
\newcommand{\NB}{\mathbb{N}\mathcal{B}_k}

\newcommand{\once}{11}
\newcommand{\onceB}{\color{blue}11}
\newcommand{\OnePar}{{ \rm Par}^1_{l,m}}
\newcommand{\OneParn}{{ \rm Par}^1_{l,n}}
\newcommand{\OO}{\mathcal{O}}
\newcommand{\op}{\otimes}

\newcommand{\Par}{{\rm Par}_{l,m}}

\newcommand{\q}{\hat{q}}
\newcommand{\quince}{15}
\newcommand{\quinceB}{\color{red} 15}

\newcommand{\R}{ \mathcal{R}_m}
\newcommand{\Rad}{{\rm Rad}}
\newcommand{\res}{ \textrm{res} }
\newcommand\rrn{\mathbf{r}}

\newcommand{\s}{\mathfrak{s}}
\newcommand{\seq}{{\rm seq}_n}
\newcommand{\shape}{\textsf{shape}}
\newcommand{\Si}{\mathfrak{S}}
\newcommand{\sign}{-}
\newcommand{\sixteenB}{\color{red} 16}
\newcommand{\Snake}{{\rm Snake}}
\newcommand{\spa}{{\rm span}}
\newcommand{\std}{{\rm Std}}
\newcommand{\SLTM}{{\rm SLTM}}
\newcommand{\SLTMn}{$\mathcal{TM}$}

\newcommand{\T}{  \mathfrak{t}}
\newcommand{\tab}{{\rm Tab}}
\newcommand{\tr}{{\rm \textbf{tr}}}
\newcommand{\tR}{ \overline{\R}}
\newcommand{\Tr}{{\rm Tr}}
\newcommand{\trece}{13}
\newcommand{\treintatres}{33}
\newcommand{\treintatresR}{\color{red}33}
\newcommand{\trunc}{ {\B} (\blambda ) }
\newcommand{\truncPrime}{ {\mathbb B}_n^{\prime} (\blambda ) }
\newcommand{\truncSing}{ {\mathbb B}_{\bar n} (\overline{\blambda} ) }
\newcommand{\TT}{{\mathfrak T}}
\newcommand{\TTc}{  \mathcal{T}}

\newcommand{\U}{\mathfrak{u}}
\newcommand{\UU}{\mathbb{U}}
\newcommand{\UUc}{\mathcal{U}}
\newcommand{\ulu}{\underline{u}}
\newcommand{\ulv}{\underline{v}}
\newcommand{\ulw}{\underline{w}}
\newcommand{\ulx}{\underline{x}}
\newcommand{\uly}{\underline{y}}
\newcommand{\ulz}{\underline{z}}

\newcommand{\V}{\mathfrak{v}}
\newcommand{\veintidos}{22}
\newcommand{\VV}{\mathbb{V}}

\newcommand{\Y}{\mathcal{Y}}
\newcommand{\Yy}{\mathcal{Y}_{r,n}}
\newcommand{\yvc}{\Yvcentermath1}
\newcommand{\YY}{\mathcal{Y}_{r,n}(q)}

\newcommand{\Z}{\mathbb{Z}}

\begin{document}
  \Yvcentermath1
\title{Invariants for isomorphism classes in the category $\bcalNT$ }
\author{\sc  Diego Lobos Maturana }

\maketitle


\begin{abstract}
 The category $\bcalNT$ is a category of certain commutative graded algebras over a field.  It was introduced in \cite{Lobos2} as a generalization of algebras generated by Jucys-Murphy elements in the many ${\textrm{\textbf{{End}}}}$ algebras of the diagrammatic Soergel category of Elias and Williamson. In the first part of this article we define certain \emph{Invariants} for the isomorphism classes in $\bcalNT,$ following in the same spirit of \cite{Lobos3}, where a series of \emph{Isomorphism criteria} were found. At the end, we use our invariants to provide a new lower bound for the number of isomorphism classes, improving a similar result obtained in \cite{Lobos3}.
\end{abstract}
keywords: Commutative graded algebras, Jucys-Murphy elements, Soergel Calculus.

\pagenumbering{arabic}

\section{Introduction}

In general, when we study the objects of certain category, one of the most natural goal is to obtain, if it is possible, a complete classification (up to isomorphism) of them. Depending on the nature of the objects and the complexity of the morphisms in the category, this goal sometimes can be achieved but most of the time it becomes an idealistic dream. However, most of the greater open problems in Pure Mathematics, are related (explicitly or implicitly) to the search of such a classification in many important categories. In the context of \emph{Algebra}, there are some well known classifications, for example, in the category of all finite generated Abelian groups, we have the well known \emph{Fundamental Theorem of Finitely Generated Abelian Groups}, where the desired complete classification was obtained (see \cite{Hungerford} for example). In the category of \emph{Coxeter groups}, there is a well known classification of their objects, based in \emph{Coxeter graphs} (see \cite{Bjorner-Brenti}, \cite{Hump1} for example). In a more general context, we have the so-called \emph{Enormous Theorem}, where a complete classification of all the finite simple groups was achieved (see \cite{Gor-Lyons-Solomon}).

On the other side, roughly speaking, the main objective of \emph{Knot Theory}, is to classify all the different \emph{Knots} up to \emph{continuous deformations}. In this case, a complete classification still seem very distant; however, having this objective in mind, researchers have been forced to produce new mathematical techniques, which, on the one hand, have allowed progress toward the desired goal and on the other, have led to the interaction of different branches of study. It is well known that the search for algebraic invariants in knot theory has generated wonderful interactions between topologists, algebraists and even physicists (see \cite{Aicardi-Juyu-2021}, \cite{Birman1993}, \cite{Juyumaya-Lambropoulou-2009}, \cite{Kauffman1991}, \cite{Murasugi2009}, for example).

Another example that we would like to mention comes from \emph{Algebraic Geometry}, where, roughly speaking, the main objective is to classify all the \emph{Algebraic Varieties}. Over the years, mathematicians have developed several tools to study these objects. Again, starting from a natural search, one can find the motivation for the development of new knowledge, that have impact in the original purpose as in other problems of mathematics. For example, Algebraic Geometry have motivated the development of a great amount of knowledge in Commutative Algebra, Homological algebra, Representation Theory and Number Theory (see \cite{Bar-Guay-Sche}, \cite{Eisenbud1995}, \cite{Ginzburg}, \cite{Kostr-Shafa}, \cite{Liu} for example).

The category $\bcalNT$ was defined in \cite{Lobos2}, is a subcategory of the category of \emph{commutative graded algebras}, where morphisms are understood as \emph{preserving degree} homomorphism of algebras. Naturally we fix as an innocent objective the search of a complete classification of the objects of $\bcalNT.$ In \cite{Lobos2} we initiated this long-term search, where first valuable remarks and results were given. Later in \cite{Lobos3} we continued this study, finding a series of isomorphism criteria. The present article can be seen as a third episode of this saga, where we define an interesting set of invariants. The nature of the objects and morphisms of $\bcalNT,$ invite us to think that the search of a complete classification of its objects could contribute to generate valuable knowledge not only for the category $\bcalNT$ itself, but for Commutative Algebra in general and for related areas such that Algebraic Geometry and Number Theory as well. Although in this article we only provide a purely algebraic approach to the study of $\bcalNT,$ we expect that our results could be applied in those other areas in future works.

Originally, the definition of $\bcalNT$ was motivated by the study of \emph{Jucys-Murphy} elements in the diagrammatic Soergel category $\bcalD,$ of Elias and Williamson (defined in \cite{EW}, as a diagrammatic incarnation of the Hecke category. See also \cite{EW-intro}, \cite{LibGentle}, \cite{EliasKhov10}, \cite{So}, \cite{So1}). More specifically,  the objects of $\bcalNT$ appear as a generalization of the algebras generated by the Jucys-Murphy elements, found by Ryom-Hansen in \cite{Steen2}, relative to the Libedinsky's \emph{cellular bases} of the many $\textrm{\textbf{End}}$ algebras of $\bcalD$ (the so-called \emph{light leaves bases}, originally defined in  \cite{LibLightLeaves}).  Some of those $JM-$subalgebras were studied previously by Espinoza and Plaza in \cite{Esp-Pl}, connecting two realizations of them, one in the context of the category $\bcalD$ and the other in the context of the \emph{blob algebras} of Martin and Saleur (see \cite{Mart-Sal}). Our first approach to this algebras was made in \cite{Lobos1}, but that time, only in the context of the \emph{generalized blob algebras} of Martin and Woodcook (see \cite{Mart-Wood}), realized as certain quotients of \emph{KLR algebras} (see  \cite{Lobos-Ryom-Hansen} and \cite{PlazaRyom}).

It is well known from the work of Mathas, that the subalgebra generated by the \emph{Jucys-Murphy} elements (from now $JM-$subalgebra) in a cellular algebra $\mathcal{A},$ is a commutative subalgebra that under certain conditions, it could be a maximal commutative subalgebra of $\calA$  (see \cite{Mat-So}). Therefore, the $JM-$subalgebra contains important information on the structure of the cellular algebra $\mathcal{A}$ itself. \emph{Jucys-Murphy} elements also play an important role in the study of the \emph{Representation Theory} of the corresponding cellular algebra $\calA$ (see \cite{hu-mathas}, \cite{Lobos-Ryom-Hansen}, \cite{Mat}, \cite{Mat-So}, \cite{Murphy1}, \cite{Murphy2}, \cite{Murphy3}, \cite{Murphy4}, \cite{Ormeno-RH-2025}, \cite{PlazaRyom}, for example). In the case of the category of Elias and Williamson, we expect that all the information that we can find from the $JM-$subalgebras, could be relevant for a better understanding of the structure and representation theory of the corresponding $\textrm{\textbf{End}}$ algebras and therefore a better understanding of the category $\bcalD$ itself.

The category $\bcalNT$ provides a systematic and general approach to the study of the $JM-$algebras arising from the category $\bcalD.$ A typical object of $\bcalNT$ is a graded commutative algebra $\calA=\calA(T),$ given by a presentation, whose relations depends on a strictly lower matrix $T$ with entries in the ground field (see Definition \ref{Def-A(T)} below). The existence of an isomorphism between two of those algebras $\calA(T),\calA(S),$ naturally define an equivalence relation $T\sim S$ on the defining matrices. Therefore to study isomorphisms classes in $\bcalNT$ is equivalent to study the equivalent classes of matrices under that relation.  In \cite{Lobos2}, \cite{Lobos3} and the present article, we concentrate our effort on looking for relevant information coming from the matrices, that allow to develop \emph{isomorphism criteria} for the algebras in $\bcalNT.$  In some way, we are emulating the classification of Coxeter groups, where presentations and valuable information is taken from codifying elements, such as \emph{Coxeter graphs}, or equivalently, \emph{Coxeter matrices} (see \cite{Bjorner-Brenti} or \cite{Hump1}).

The following are some advances in our search:

\begin{description}
  \item[1] From \cite{Lobos2}, we know that if $T$ is a $n\times n-$matrix, then the dimension of $\calA(T)$ is $2^{n},$ therefore, we have a first criterium: If $T,S$ are two matrices of different size, then we know that $\calA(T)$ and $\calA(S)$ are not isomorphic. That is, the size of the matrix $T$ acts as a first invariant for our category.
  \item[2] In \cite{Lobos3}, we have pointed out, that if $T,S$ are two matrices of the same size, and we assume that $\bgamma:\calA(T)\rightarrow \calA(S)$ is an isomorphism, then there is an invertible matrix $\Gamma$ (of the same size of $T$ and $S$) whose entries satisfy a particular system of equations involving the entries of $T$ and $S$ (See equation \ref{Eq1-theo-Key-condition}). Although in practice, to find a solution of that system, is not something easy to reach, specially when the size of the matrices is too large, if we are able to find a solution, we will know that the algebras $\calA(T)$ and $\calA(S)$ are isomorphic. Conversely if we are able to prove that the system in equation \ref{Eq1-theo-Key-condition}, has no solution, then we will know that the algebras $\calA(T)$ and $\calA(S)$ are not isomorphic. This is an absolute criterium, but very difficult to apply.
  \item[3] In \cite{Lobos3}, we define a series of \emph{matrix operations} (denominated ETO's), emulating the classical matrix operations of linear algebra. We have proved that if there is a finite sequence of ETO's that transform $T$ into $S,$ then the algebras $\calA(T)$ and $\calA(S)$ are isomorphic (see Theorem \ref{theo-sim-vs-approx}). In this case, ETO's, provide a way to navigate inside of an isomorphism class. Theoretically, one can describe a big part of a class of a an algebra $\calA(T)$ by applying different sequences of ETO's. Unfortunately, if we have two matrices $T,S$ of the same size, in general, is very difficult to know if there is, or not, a sequence of ETO's that transform $T$ on $S.$
\item[4] In the present article (Subsection \ref{ssec-invariants-1REF}, Definition \ref{def-first-red-form} and Proposition \ref{Prop-Agorithm-1REF}), we define the \emph{First reduced eschelon form} ($1-$REF) associated to a matrix $T,$ a new matrix obtained from $T$ by a sequence of ETO's. As a consequence of Theorem \ref{theo-sim-vs-approx}, if two matrices $T,S$ are associated to the same $1-$REF, then we can conclude that $T\sim S.$

\item[5] In the present article (Subsection \ref{ssec-invariants-1REF}, Definition \ref{def-Blocks-1REF}), we associate, for each matrix $T$ a vector called the \emph{wall} of $T,$ and denoted by  $\mathbb{W}(T).$ The wall of a matrix acts as an invariant for the isomorphism class of the corresponding algebra $\calA(T).$ That is, if $T,S$ are two matrices of the same size and we know that $\mathbb{W}(T)\neq \mathbb{W}(S),$ then we can conclude that the algebras $\calA(T)$ and $\calA(S)$ are not isomorphic (see Theorem \ref{theo-invariants-1REF} and Corollary \ref{coro2-theo-invariants-1-REF}).

\item[6] In the present article (Subsection \ref{ssec-invariants-2REF}, Definition \ref{def-second-red-form} and Proposition \ref{Prop-Agorithm-2REF}), we define the \emph{Second reduced eschelon form} ($2-$REF) associated to a matrix $T,$ a new matrix obtained from $T$ by a sequence of ETO's. As a consequence of Theorem \ref{theo-sim-vs-approx}, if two matrices $T,S$ are associated to the same $2-$REF, then we can conclude that $T\sim S.$

\item[7] In the present article (Subsection \ref{ssec-invariants-2REF}, Definition \ref{def-Blocks-2REF}),  we associate for each matrix $T$ a sequence of matrices, called the \emph{measure sequence} of $T,$ and denoted by $\mathbb{M}(T).$ The measure sequence of a matrix acts as an invariant for the isomorphism class of the corresponding algebra $\calA(T)$ (See Theorem \ref{theo-invariants-2REF} and Corollary \ref{coro2-theo-invariants-2REF}). This invariant is stronger than the wall, in the sense, that we can have $\mathbb{W}(T)=\mathbb{W}(S)$ but $\mathbb{M}(T)\neq \mathbb{M}(S)$ (implying that $T\nsim S$). On the other hand $\mathbb{M}(T)= \mathbb{M}(S)$ implies $\mathbb{W}(T)=\mathbb{W}(S),$ by definition.
 \end{description}

The content of this paper is organized as follows; we start Section 2 by fixing some notations and terminology, later we recall definitions and previous results coming from \cite{Lobos2} and \cite{Lobos3}. In Section 3, we define the main object of study in this article, two invariants for the isomorphisms classes in the category $\bcalNT,$ the \emph{Wall} and the \emph{Measure sequence}. Finally in Section 4, we use our invariants to provide a new lower bound to the number of isomorphism classes in the category $\bcalNT.$ The main results of this article are Theorems \ref{theo-invariants-1REF} and \ref{theo-invariants-2REF}, where we prove that the \emph{Wall} and the \emph{Measure sequence} are  effectively invariants for the isomorphism classes of $\bcalNT.$ On the other hand, Theorem \ref{theo-new-lower-bound} is an important application of our invariants, concerning with the number of isomorphism classes.


\medskip


\medskip

\section{Generalities}\label{ssec-generalities}
\subsection{Notations and terminology}\label{ssec-Notations}
\begin{itemize}
  \item $\F$ is a fixed field of characteristic different from $2.$ We denote $\F^{\times}=\F-\{0\}.$
  \item A \nSLTM,  is a strictly lower $n\times n-$matrix,  with entries in $\F.$ A \SLTM, is a strictly lower matrix (of any size), with entries in $\F.$ We denote by \nSLTMn, the set of all the \nSLTM,  and by \SLTMn,  the set of all the \SLTM.
  \item Let $a,b\in\mathbb{Z}$, with $a\leq b.$ we define the intervals $[a,b]:=\{m\in\mathbb{Z}:a\leq m \leq b\},$ $]a,b]=\{m\in\mathbb{Z}:a<m\leq b\}.$ Analogously we define $[a,b[,]a,b[,$ etc. We also denote $[a]=[a,a]=\{a\}.$
  \item If $T$ is a \nSLTM,  $I,J\subset [1,n]$ two intervals, then we denote by $T^{I}_{J}$ the submatrix of $T$ comprised by the entries $t_{ij}$ where $i\in I$ and $j\in J.$ We also denote $T^{I}:=T^{I}_{[1,n]}$ and $T_{J}:=T^{[1,n]}_{J}.$ In  particular $T^{[r]}$ denotes the $r-$th row of $T$ and $T_{[c]}$ denotes the $c-$th column of $T$
  \item If $T$ is a \nSLTM, then for each $r\in[1,n]$ we define:
  \begin{equation*}
    c_r=\left\{\begin{array}{cc}
                 \max\{j:t_{rj}\neq 0\}, & \textrm{if}\quad T^{[r]}\neq 0. \\
                 0, & \textrm{otherwise}
               \end{array}\right.
  \end{equation*}
  When $c_r>0$ we say that the entry $t_{r,c_r}$ is the \emph{leader} of the $r-$th row of $T.$ The pair $(r,c_r)$ is the \emph{leader position} of the $r-$th row of $T.$
  \item Given a \nSLTM, $T=[t_{ij}]$ and $\alpha\in\F.$ For each $1\leq i<j<k\leq n,$ we denote $\Delta_{i,j,k}^{(\alpha)}(T)=\alpha t_{ki}+t_{kj}t_{ji}.$ If there is no possible confusion we only write $\Delta_{i,j,k}^{(\alpha)}.$
\end{itemize}

\subsection{Definitions and previous results}
In this subsection we recall the category $\bcalNT$  and some of the results obtained in \cite{Lobos2} and \cite{Lobos3}:

\begin{definition}\label{Def-A(T)}
  Given a \nSLTM, $T=[t_{ij}],$ the \emph{nil graded algebra associated to} $T$ is the commutative unitary $\F-$algebra defined by generators $X_1,\dots,X_n$ and relations:
  \begin{equation}\label{EQ-def-A(T)}
    \left\{\begin{array}{c}
      X_1^{2}=0 \\
       \quad\\
      X_i^{2}=\sum_{j<i}{t_{ij}}X_jX_i,\quad (2\leq i \leq n)
    \end{array}\right.
  \end{equation}
  We denote this algebra by $\calA(T).$  We define a grading on $\calA(T),$ by $\deg(X_i)=2,$ for any $1\leq i\leq n.$
\end{definition}

\begin{theorem}\label{theo-monomial-basis}
  If $T=[t_{ij}]$ is a \nSLTM, then the dimension of the algebra $\calA(T)$ is $2^n$ and the set $\{X_1^{a_1}\cdots X_n^{a_n}:a_i\in\{0,1\}\}$ is a basis for $\calA(T).$
\end{theorem}
\begin{proof}
  See \cite{Lobos2}.
\end{proof}

\begin{definition}\label{def-category-bcalNT}
  We define the category $\bcalNT$  as the one whose objects are the many algebras  $\calA(T),$ where $T$ is a \SLTM,  with entries in the ground field $\F,$ and whose morphisms are all the \emph{preserving degree homomorphisms} between them.
\end{definition}

For the rest of the article, when we say that $\bgamma:\calA(T)\rightarrow\calA(S)$ is a morphism (resp. isomorphism), we shall understand that we refer to a morphism (resp. isomorphism) in the category $\bcalNT.$ When we say that the algebras $\calA(T)$ and $\calA(S)$ are isomorphic, we mean that they are isomorphic as objects of the category $\bcalNT.$

Note that if $\bgamma:\calA(T)\rightarrow \calA(S)$ is a morphism, then $\bgamma$ defines a matrix $\Gamma=[\gamma_{ij}]$ given by the equation
\begin{equation*}
  \bgamma(X_j)=\sum_{i}\gamma_{ij}Y_i
\end{equation*}
where $\{X_j\},\{Y_i\}$ are the set of generators (as in definition \ref{Def-A(T)}) of $\calA(T)$ and $\calA(S)$ respectively. If $\bgamma$ is an isomorphism of the category $\bcalNT,$ then $\Gamma$ is an invertible matrix. Moreover we have:

\begin{theorem}\label{theo-Key-condition}
Let $T,S$ be two \nSLTM.  Given a $n\times n$ invertible matrix, $\Gamma=[\gamma_{ij}],$ then the assignment
  \begin{equation*}
    X_j\mapsto \sum_{i=1}^{n}{\gamma_{ij}}{Y_i}
  \end{equation*}
  defines an isomorphism $\bgamma:\calA(T)\rightarrow\calA(S),$ if and only if the following system of equations holds:
  \begin{equation}\label{Eq1-theo-Key-condition}
    2\gamma_{ir}\gamma_{kr}+\gamma_{kr}^{2}s_{ki}=
    \sum_{j<r}{t_{rj}}(\gamma_{kj}\gamma_{kr}s_{ki}
    +\gamma_{kj}\gamma_{ir}+\gamma_{ij}\gamma_{kr})
    ,\quad (1\leq r\leq n;1\leq i<k\leq n).
  \end{equation}
\end{theorem}
\begin{proof}
  See \cite{Lobos3}.
\end{proof}

Given $T,S$ two \nSLTM,  we denote by $T\sim S$ whenever the algebras $\calA(T),\calA(S)$ are isomorphic. Relation $\sim$ is an equivalence relation on the set of all \nSLTM. We denote by $\overline{T}$ the class of the matrix $T$ under relation $\sim.$

The class $\overline{0_n}$ of the zero matrix $0_n$ was completely described in \cite{Lobos3} via the following theorem:

\begin{theorem}\label{theo-zero-class-n-gen}
  Given a \nSLTM, $U\neq 0_n,$ then $U\in \overline{0_n}$ if and only if one of the following conditions holds
  \begin{enumerate}
    \item\label{theo-zero-class-n-gen-It1} The leader of any nonzero row is in the first column.
    \item\label{theo-zero-class-n-gen-It2} For every position $(r,c)$ such that $c>1$ and $u_{rc}\neq 0,$ we have that $\Delta^{(2)}_{i,c,r}(U)=0,$ for each  $i\in [1,c[.$
  \end{enumerate}
\end{theorem}

\begin{proof}
  See \cite{Lobos3}.
\end{proof}

Another criterium to determine when two matrices are equivalent under relation $\sim$ is given by that we call \emph{Elementary triangular operations} (ETO's), we now recall the definition of the ETO's and the Theorem obtained in \cite{Lobos3}, that connect ETO's with relation $\sim.$


\begin{definition}\label{def-ETOs}
Let $T=[t_{rk}]$ be a \nSLTM:
\begin{enumerate}
  \item  Fix an $1\leq r_1\leq n$ and an scalar $\alpha\neq 0.$ We define the matrix $S=\calP_{r_1}(T,\alpha),$ as the \nSLTM, whose entries are given by equation \ref{eq-def-P-operation}
  \begin{equation}\label{eq-def-P-operation}
    s_{rk}=\left\{\begin{array}{cc}
                    \alpha^{-1} t_{rk} & \textrm{if}\quad r=r_1 \\
                    \alpha t_{rk} & \textrm{if}\quad k=r_1 \\
                    t_{rk} & \textrm{otherwise}
                  \end{array}\right.
  \end{equation}
  \item Fix $1\leq r_1<r_2\leq n.$ If the matrix $T$ satisfies the following conditions:
  \begin{enumerate}
    \item The entry $t_{r_2j}=0$ for each $r_1\leq j$
    \item The entry $t_{rr_1}=0$ for each $r_1\leq r\leq r_2$
  \end{enumerate}
  then we define the matrix $S=\calF_{(r_1,r_2)}(T)$ as the \nSLTM, whose entries are given by the following equation:
  \begin{equation}\label{eq-def-F-operation}
    s_{rk}=\left\{\begin{array}{cc}
                    t_{r_2k} & \textrm{if}\quad r=r_1 \\
                    t_{r_1k} & \textrm{if}\quad r=r_2  \\
                    t_{rr_2} & \textrm{if}\quad k=r_1 \\
                    t_{rr_1} & \textrm{if}\quad k=r_2 \\
                    t_{rk} & \textrm{otherwise}
                  \end{array}\right.
  \end{equation}
  \item Fix a position $(r_0,k_0)$ of $T,$ such that $r_0>k_0$ and $t_{r_0,k}=0$ for all $k>k_0.$
  \begin{enumerate}
    \item If $k_0=1,$ then for any $\beta\in \F^{\times},$ we define the matrix $S=\calQ_{(r_0,1)}(T,\beta)$ as the \nSLTM, whose entries $s_{rk}$ are given by equation \ref{eq-def-Q-operation}.
    \item If $k_0>1$ and there is a $\beta\in \F^{\times},$ such that $\Delta^{(1)}_{i,k_0,r_0}(T)=\beta t_{k_0,i},$ for each $1\leq i<k_0,$ then we define the matrix $S=\calQ_{(r_0,k_0)}(T,\beta)$ as the \nSLTM, whose entries $s_{rk}$ are given by equation \ref{eq-def-Q-operation}.
  \end{enumerate}
  \begin{equation}\label{eq-def-Q-operation}
    s_{rk}=\left\{\begin{array}{cc}
                    t_{r_0,k_0}-2\beta & \textrm{if}\quad (r,k)=(r_0,k_0) \\
                    \quad & \quad\\
                    t_{r,k_{0}}+\beta t_{r,r_{0}} & \textrm{if}\quad r>r_0\quad \textrm{and}\quad k=k_0  \\
                    \quad & \quad\\
                    t_{rk} & \textrm{otherwise}
                  \end{array}\right.
  \end{equation}
\end{enumerate}



We call \emph{Elementary triangular operations} (ETO for short) any of the following matrix assignments (whenever they are well defined):

\begin{equation}\label{eq-ETO-assignments}
  \begin{array}{cc}
    \calP_r(\quad,\alpha):T\mapsto \calP_r(T,\alpha)&(\calP-\textrm{operation} ) \\
    \quad&\quad  \\
     \calF_{(r_1,r_2)}(\quad):T\mapsto \calF_{(r_1,r_2)}(T)& (\calF-\textrm{operation}) \\
    \quad &\quad \\
    \calQ_{(r_0,k_0)}(\quad,\beta):T\mapsto \calQ_{(r_0,k_0)}(T,\beta)& (\calQ-\textrm{operation})
  \end{array}
\end{equation}
\end{definition}

\begin{theorem}\label{theo-sim-vs-approx}
  Let $T,S$ be two \nSLTM.  If there exists a (finite) sequence of ETO's that transform $T$ into $S,$ then $T\sim S.$
\end{theorem}
\begin{proof}
  See \cite{Lobos3}.
\end{proof}

\section{Invariants}\label{sec-Invariants}

In this section we describe a series of invariants for the isomorphism classes of objects in the category $\bcalNT.$ All the results of this section are true for any field $\F$ of characteristic different from $2,$ despite that, in all the examples of this section we shall assume that $\F$ is a field of characteristic equal to zero.

\subsection{Walls and Bricks}\label{ssec-invariants-1REF}

\begin{definition}\label{def-first-red-form}
  Let $S\neq 0_n$ be a \nSLTM, we say that $S$ is in the \emph{First Reduced Echelon Form} ($1-$REF), if it satisfies the following conditions:
  \begin{enumerate}
    \item The zero rows of $S$ are above of the nonzero rows.
    \item For each leader position $(r,c)$ we have that $c>1$ and there is a $i_0\in [1,c[$ such that $\Delta^{(2)}_{i_0,c,r}\neq 0.$
    \item If $(r,c_r)$ and $(r+1,c_{r+1})$ are the leader positions of the $r-$th and $(r+1)-$th rows of $S,$ then $c_r\leq c_{r+1}.$
  \end{enumerate}

  We also say that the zero matrix, $0_n$ is in the first reduced echelon form.
\end{definition}

\begin{remark}\label{remark-def-1-red-form}
Let $S$ be a \SLTM, then
\begin{itemize}
  \item If $S$ is in the $1-$REF, then $S$ does not contains leaders located in the first column.
  \item If we had a leader position $(r,c)$ with $c>1$ satisfying that $\Delta_{i,c,r}^{(2)}=0,$ for each $i\in[1,c[,$ then $s_{ri}=-\left(\frac{s_{rc}}{2}\right)s_{ci}$ and therefore $\Delta^{(1)}_{i,j,k}=s_{ri}+s_{rc}s_{ci}=\left(\frac{s_{rc}}{2}\right)s_{ci}.$ The last equation implies that we are able to define $U=\calQ_{(r,c)}\left(S,\frac{s_{rc}}{2}\right)$ and this new matrix have a zero in its position $(r,c)$ (see definition \ref{def-ETOs}). Therefore the first condition in definition \ref{def-first-red-form} implies that if $S$ is in the $1-$REF, then we are not able to form a zero in a leader position of $S.$
\end{itemize}
\end{remark}

\begin{example}\label{ex-1-red-form}
  The following matrices are not in the first reduced echelon form:
  \begin{equation}\label{eq1-ex-1-red-form}
   T_1=\left[\begin{matrix}
            0 & 0 & 0 & 0& 0 \\
            0 & 0 & 0 & 0& 0 \\
            0& 0& 0& 0& 0\\
            2 & 2 & 0 & 0& 0 \\
            -1 & -1 & 1 & 0& 0
          \end{matrix}\right],\quad
   T_2=\left[\begin{matrix}
            0 & 0 & 0 & 0& 0 \\
            0 & 0 & 0 & 0& 0 \\
            0 & 0 & 0 & 0& 0 \\
            1 & 1 & 1 & 0& 0\\
            1& 1& 0& 0& 0
          \end{matrix}\right],\quad T_3=\left[\begin{matrix}
            0 & 0 & 0 & 0& 0 \\
            2 & 0 & 0 & 0& 0 \\
            1 & 2 & 0 & 0& 0 \\
            0 & 4 & 3 & 0& 0\\
            0& 0& 0& 0& 0
          \end{matrix}\right],
  \end{equation}
  \begin{itemize}
   \item In $T_1$ the leader of the fourth row is in the position $(4,3),$ but in this case $\Delta^{(2)}_{i,3,4}=0,$ for all $i\in [1,3[,$ then the second condition fails.
    \item In $T_2$ the leader of the third row is in the third column, while the leader of the fourth row is in the second column, then the third condition fails.
    \item Finally in $T_3$ the second row has its leader in the first column and also the last row is a zero row, but it is below a nonzero row, then the first and the second conditions fail.
  \end{itemize}
  The following matrices are in the first reduced echelon form:
   \begin{equation}\label{eq2-ex-1-red-form}
   S_1=\left[\begin{matrix}
            0 & 0 & 0 & 0& 0 \\
            0 & 0 & 0 & 0& 0 \\
            0& 0& 0& 0& 0\\
            2 & 2 & 0 & 0& 0 \\
            -1 & -1 & 0 & 0& 0
          \end{matrix}\right],\quad
   S_2=\left[\begin{matrix}
            0 & 0 & 0 & 0& 0 \\
            0 & 0 & 0 & 0& 0 \\
            0 & 0 & 0 & 0& 0 \\
            1 & 1 & 0 & 0& 0\\
            1& 1& 1& 0& 0
          \end{matrix}\right],\quad
   S_3=\left[\begin{matrix}
            0 & 0 & 0 & 0& 0 \\
            0 & 0 & 0 & 0& 0 \\
            0 & 0 & 0 & 0& 0 \\
            3 & 2 & 0 & 0& 0\\
            4& 4& 3& 0& 0
          \end{matrix}\right],
  \end{equation}
\begin{itemize}
    \item Since $\Delta^{(2)}_{i,3,4}(T_1)=0,$ for all $i\in [1,3[,$ then we can apply a convenient $\calQ-$operation  obtaining $S_1:$
        \begin{equation*}
          T_1\rightarrow S_1=\calQ_{(4,3)}\left(T_2,\frac{1}{2}\right).
        \end{equation*}
    \item In $T_2$ we can apply a suitable $\calF-$operation obtaining $S_2:$
    \begin{equation*}
      T_2\rightarrow S_2=\calF_{(4,5)}(T_2).
    \end{equation*}
    \item Note that $S_3$ can be obtained from $T_3$ by applying a convenient sequence of ETO's:
        \begin{equation*}
          T_3\rightarrow U_1=\calQ_{(2,1)}(T_3,1)\rightarrow U_2=\calF_{(4,5)}(U_1)\rightarrow S_3=\calF_{(3,4)}(U_2).
        \end{equation*}
  \end{itemize}
\end{example}

In example \ref{ex-1-red-form} we show three matrices that were not in the first reduced echelon form, but after a process of ETO's, we could see that they were in the same class (under relation $\sim$) with a matrix in the first reduced echelon form. The following proposition generalized this fact:

\begin{proposition}\label{Prop-Agorithm-1REF}
  Let $T$ be a \nSLTM, then there is a sequence of ETO's that transform $T$ into a matrix $R$ that is in the first reduced echelon form.
\end{proposition}

Recall that for each nonzero row of $T,$ we denote by $(r,c_r)$ the corresponding leader position of the $r-th$ row of $T.$ If $T^{[r]}=0,$ then we define $c_r=0.$ (See subsection \ref{ssec-Notations}).

\begin{proof}
Assuming that $T$ is not in the $1-$REF:
\begin{description}
\item[Step 1] For each leader $t_{rc}$ such that $c=1$ we simply apply the operation $\calQ_{(r,1)}(\quad,\beta)$ with $\beta=\frac{t{r1}}{2}$ obtaining a new matrix $T,$ whose corresponding entry $t_{r1}$ is equal to zero. Note that this action do not introduce new leaders in the first column.
  We repeat progressively this process until we obtain a new matrix $T$ with no leaders in the first column. Note that this process only change the entries of the first column of the original matrix $T.$
\item [Step 2] (Starting with the new matrix $T$ obtained in Step 1) For each leader $t_{rc}$ such that $c>1$ and $\Delta_{j,c,r}^{(2)}=0$ for all $j\in[1,c[,$ we apply the operation $\calQ_{(r,c)}(\quad,\beta)$ with $\beta=\frac{t_{rc}}{2},$ obtaining a new matrix $T$ whose corresponding entry $t_{rc}$ is zero. Note that this action affects the entries $t_{kc}$ for $k\geq r$ whenever the entry $t_{kr}\neq 0.$ Any other entry of $T$ keeps its original value, therefore this action do not introduce new leaders for the new matrix $T.$ Applying repeatedly this action we obtain a new matrix $T$ with no leaders in the first column and for each leader $t_{rc}$ with $c>1$ there will be a $j_0\in[1,c[$ such that $ \Delta_{j_0,c,r}^{(2)}\neq 0.$
\item[Step 3] (Starting with the new matrix $T$ obtained in Step 2) Applying a finite sequence of admissible $\calF-$operations, we can put all the zero rows above the other. After that we apply a second sequence of admissible $\calF-$operations in order to obtain a new matrix $T$ where the columns corresponding to the leader positions are disposed increasingly.
\item [Step 4] We define $R$ as the last matrix $T$ obtained in Step 3.
\end{description}
By construction, $R$ is a matrix in $1-$REF that is in the same class (under relation $\sim$) as the original matrix $T.$
\end{proof}


\begin{example}\label{ex-algorithm1}
  Consider the following matrix:
  \begin{equation*}
    T=\left[\begin{matrix}
              0 & 0 & 0 & 0 & 0 & 0 & 0 & 0 & 0 \\
              -1 & 0 & 0 & 0 & 0 & 0 & 0 & 0 & 0 \\
              0 & 0 & 0 & 0 & 0 & 0 & 0 & 0 & 0 \\
              0 & 0 & -1 & 0 & 0 & 0 & 0 & 0 & 0 \\
              -1 & -1 & -1 & -1 & 0 & 0 & 0 & 0 & 0 \\
              -1 & -1 & 0 & 1 & -1 & 0 & 0 & 0 & 0 \\
              -1 & -1 & 1 & 0 & -1 & -1 & 0 & 0 & 0 \\
              0 & 1 & -1 & -1 & -1 & 0 & 0 & 0 & 0 \\
              1 & 0 & -1 & -1 & -1 & 0 & 0 & -1 & 0
            \end{matrix}\right]
  \end{equation*}
  \begin{description}
    \item[Step 1] The objective of Steps 1 and 2 is to clean all the \emph{noise} that we can find in the matrix $T,$ that is, for each row, we ask if it is possible to apply a $\calQ-$operation that make zero the leader of that row.

        In this case, in the second row of $T$ we find a leader in the first column, then we can apply $\calQ_{(2,1)}(\quad,\frac{-1}{2})$ to $T.$ Let $T_1=\calQ_{(2,1)}(\quad,\frac{-1}{2}),$ then we obtain:
        \begin{equation*}
    T_1=\left[\begin{matrix}
              0 & 0 & 0 & 0 & 0 & 0 & 0 & 0 & 0 \\
              \quad &\quad &\quad &\quad &\quad &\quad &\quad &\quad &\quad \\
              0 & 0 & 0 & 0 & 0 & 0 & 0 & 0 & 0 \\
              \quad &\quad &\quad &\quad &\quad &\quad &\quad &\quad &\quad \\
              0 & 0 & 0 & 0 & 0 & 0 & 0 & 0 & 0 \\
              \quad &\quad &\quad &\quad &\quad &\quad &\quad &\quad &\quad \\
              0 & 0 & -1 & 0 & 0 & 0 & 0 & 0 & 0 \\
              \quad &\quad &\quad &\quad &\quad &\quad &\quad &\quad &\quad \\
              -\frac{1}{2} & -1 & -1 & -1 & 0 & 0 & 0 & 0 & 0 \\
              \quad &\quad &\quad &\quad &\quad &\quad &\quad &\quad &\quad \\
               -\frac{1}{2} & -1 & 0 & 1 & -1 & 0 & 0 & 0 & 0 \\
               \quad &\quad &\quad &\quad &\quad &\quad &\quad &\quad &\quad \\
               -\frac{1}{2} & -1 & 1 & 0 & -1 & -1 & 0 & 0 & 0 \\
               \quad &\quad &\quad &\quad &\quad &\quad &\quad &\quad &\quad \\
              -\frac{1}{2} & 1 & -1 & -1 & -1 & 0 & 0 & 0 & 0 \\
              \quad &\quad &\quad &\quad &\quad &\quad &\quad &\quad &\quad \\
              1 & 0 & -1 & -1 & -1 & 0 & 0 & -1 & 0
            \end{matrix}\right]
  \end{equation*}
  \item[Step 2]
  Now in the fourth row, the leader position is $(4,3).$ Is not difficult to see that $\Delta_{(1,3,4)}^{(2)}=\Delta_{(2,3,4)}^{(2)}=0,$ then we can apply $\calQ_{(4,3)}(\quad,\frac{-1}{2})$ to $T_1.$ Make $T_2=\calQ_{(4,3)}(T_1,\frac{-1}{2}),$ then
   \begin{equation*}
    T_2=\left[\begin{matrix}
              0 & 0 & 0 & 0 & 0 & 0 & 0 & 0 & 0 \\
              \quad &\quad &\quad &\quad &\quad &\quad &\quad &\quad &\quad \\
              0 & 0 & 0 & 0 & 0 & 0 & 0 & 0 & 0 \\
              \quad &\quad &\quad &\quad &\quad &\quad &\quad &\quad &\quad \\
              0 & 0 & 0 & 0 & 0 & 0 & 0 & 0 & 0 \\
              \quad &\quad &\quad &\quad &\quad &\quad &\quad &\quad &\quad \\
              0 & 0 & 0 & 0 & 0 & 0 & 0 & 0 & 0 \\
              \quad &\quad &\quad &\quad &\quad &\quad &\quad &\quad &\quad \\
              -\frac{1}{2} & -1 & -\frac{1}{2} & -1 & 0 & 0 & 0 & 0 & 0 \\
              \quad &\quad &\quad &\quad &\quad &\quad &\quad &\quad &\quad \\
               -\frac{1}{2} & -1 & -\frac{1}{2} & 1 & -1 & 0 & 0 & 0 & 0 \\
               \quad &\quad &\quad &\quad &\quad &\quad &\quad &\quad &\quad \\
               -\frac{1}{2} & -1 & 1 & 0 & -1 & -1 & 0 & 0 & 0 \\
               \quad &\quad &\quad &\quad &\quad &\quad &\quad &\quad &\quad \\
              -\frac{1}{2} & 1 & -\frac{1}{2} & -1 & -1 & 0 & 0 & 0 & 0 \\
              \quad &\quad &\quad &\quad &\quad &\quad &\quad &\quad &\quad \\
              1 & 0 & -1 & -\frac{1}{2} & -1 & 0 & 0 & -1 & 0
            \end{matrix}\right]
  \end{equation*}
  One can check that there is no more leaders that can be eliminated by a $\calQ-$operation.
    \item[Step 3] The objective of Step 3 is reorder the positions of the rows via $\calF-$operations, in order to satisfy the third condition of Definition \ref{def-first-red-form}. In this case if we define $R=\calF_{(7,8)}(T_2)$ we obtain a matrix, in $1-$REF such that $T\sim R.$
     \begin{equation*}
    S=\left[\begin{matrix}
              0 & 0 & 0 & 0 & 0 & 0 & 0 & 0 & 0 \\
              \quad &\quad &\quad &\quad &\quad &\quad &\quad &\quad &\quad \\
              0 & 0 & 0 & 0 & 0 & 0 & 0 & 0 & 0 \\
              \quad &\quad &\quad &\quad &\quad &\quad &\quad &\quad &\quad \\
              0 & 0 & 0 & 0 & 0 & 0 & 0 & 0 & 0 \\
              \quad &\quad &\quad &\quad &\quad &\quad &\quad &\quad &\quad \\
              0 & 0 & 0 & 0 & 0 & 0 & 0 & 0 & 0 \\
              \quad &\quad &\quad &\quad &\quad &\quad &\quad &\quad &\quad \\
              -\frac{1}{2} & -1 & -\frac{1}{2} & -1 & 0 & 0 & 0 & 0 & 0 \\
              \quad &\quad &\quad &\quad &\quad &\quad &\quad &\quad &\quad \\
               -\frac{1}{2} & -1 & -\frac{1}{2} & 1 & -1 & 0 & 0 & 0 & 0 \\
               \quad &\quad &\quad &\quad &\quad &\quad &\quad &\quad &\quad \\
               -\frac{1}{2} & 1 & -\frac{1}{2} & -1 & -1 & 0 & 0 & 0 & 0 \\
              \quad &\quad &\quad &\quad &\quad &\quad &\quad &\quad &\quad \\
               -\frac{1}{2} & -1 & 1 & 0 & -1 & -1 & 0 & 0 & 0 \\
               \quad &\quad &\quad &\quad &\quad &\quad &\quad &\quad &\quad \\
              1 & 0 & -1 & -\frac{1}{2} & -1 & 0 & -1 & 0 & 0
            \end{matrix}\right]
  \end{equation*}
  Recall that operation $\calF_{(7,8)}$ also affects the $7$th and the $8$th columns of the matrix.
  \end{description}
\end{example}

In the following definition, we shall introduce our first \emph{Invariant} for the isomorphism classes in the category $\bcalNT:$

\begin{definition}\label{def-Blocks-1REF}
  Let $S$ be a \nSLTM, such that $S\nsim 0_n.$ Assuming that $S$ is in the $1-$REF:
  \begin{enumerate}
   \item We define recursively numbers $r_j$ and a vector $\mathbb{W}(S)$ associated to $S$ as follows:
      \begin{enumerate}
      \item Let $r_1:=\min\{r\in[1,n]:S^{[r]}\neq 0\}$.
      \item Assuming the number $r_j$ is already defined, define the set $A_j=\{r\in]r_j,n]:c_r>r_j\}.$
        \begin{itemize}
          \item If $A_j\neq \emptyset$ then we define $r_{j+1}=\min(A_j).$
          \item If $A_j=\emptyset,$ then we define the  vector $\mathbb{W}(S)=(r_1,\dots,r_j)$.
        \end{itemize}
        The vector $\mathbb{W}(S)$ will be called the \emph{wall} of $S,$
      \end{enumerate}
   \item If $\mathbb{W}(S)=(r_1),$ then we define the \emph{brick} of $S$ by:
     \begin{equation*}
       B_1(S):=S_{[1,r_1[}^{[r_1,n]}
     \end{equation*}
   \item If $\mathbb{W}(S)=(r_1,\dots,r_h),$ with $h>1,$ then for each $j\in[1,h]$ we define the $j-$th \emph{brick} of $S$ as the submatrix:
       \begin{equation*}
         B_{j}(S)=\left\{\begin{array}{cc}
                        S_{[1,r_1[}^{[r_1,r_{2}[} & \textrm{if}\quad j=1 \\
                        \quad & \quad \\
                        S_{[r_{j-1},r_j[}^{[r_j,r_{j+1}[} & \textrm{if}\quad 1<j<h\\
                        \quad & \quad \\
                        S_{[r_{h-1},r_h[}^{[r_h,n]}&\textrm{if}\quad j=h
                      \end{array}\right.
       \end{equation*}
  \end{enumerate}
  For the zero matrix we also define $\mathbb{W}(0_n)=(0).$ 
\end{definition}

Sometimes, if $n$ is fixed, when we write generically $\mathbb{W}(S)=(r_1,\dots,r_{h}),$ we will denote $r_0:=1$ and $r_{h+1}=n+1,$ then  $[r_0,r_1[=[1,r_1[$ and $[r_h,r_{h+1}[=[r_h,n+1[=[r_h,n].$

Note that by Definition \ref{def-Blocks-1REF}, $S$ looks like this:
\begin{equation*}
  S=\left[\begin{matrix}
                0 & 0 & \cdots & 0&0 \\
                B_1(S) & 0 &\cdots & 0&0 \\
                \bullet & B_{2}(S)&\cdots & 0&0\\
                \vdots & \vdots &\ddots&\vdots&\vdots\\
                \bullet & \bullet & \bullet& B_{h}(S)&0
              \end{matrix}\right]
\end{equation*}

 In particular, for any $r\in[r_j,r_{j+1}[$ we have that the leader $s_{rc_r}$ of the $r-$th row of $S$ belongs to the brick $B_{j}(S),$ that is $c_r<r_j.$

\begin{example}
  Considering the matrices $S_j$ of equation \ref{eq2-ex-1-red-form}. As we have mentioned in example \ref{ex-1-red-form}, the three matrices are in the $1-$REF, so we can calculate the wall $\mathbb{W}(S_j)$ in each case. Is not difficult to check for each case that $\mathbb{W}(S_j)=(4)$ then $S_j$ contains only one brick. In the following equation we have painted in blue for each case the entries of the corresponding brick of the matrix $S_j:$

   \begin{equation}\label{eq3-ex-1-red-form}
   S_1=\left[\begin{matrix}
            0 & 0 & 0 & 0& 0 \\
            0 & 0 & 0 & 0& 0 \\
            0& 0& 0& 0& 0\\
            {\color{blue}2} & {\color{blue}2} & {\color{blue}0} & 0& 0 \\
            {\color{blue}-1} & {\color{blue}-1} & {\color{blue}0} & 0& 0
          \end{matrix}\right],\quad
   S_2=\left[\begin{matrix}
            0 & 0 & 0 & 0& 0 \\
            0 & 0 & 0 & 0& 0 \\
            0 & 0 & 0 & 0& 0 \\
            {\color{blue}1} & {\color{blue}1} & {\color{blue}0} & 0& 0\\
            {\color{blue}1}& {\color{blue}1}& {\color{blue}1}& 0& 0
          \end{matrix}\right],\quad
   S_3=\left[\begin{matrix}
            0 & 0 & 0 & 0& 0 \\
            0 & 0 & 0 & 0& 0 \\
            0 & 0 & 0 & 0& 0 \\
            {\color{blue}3} & {\color{blue}2} & {\color{blue}0} & 0& 0\\
            {\color{blue}4}& {\color{blue}4}& {\color{blue}3}& 0& 0
          \end{matrix}\right],
  \end{equation}
\end{example}

\begin{example}\label{ex-walls}
  Considering the matrices
  \begin{equation*}
    S=\left[\begin{matrix}
              0 & 0 & 0& 0& 0& 0 \\
              0 & 0 & 0& 0& 0& 0 \\
              1 & 1 & 0& 0 & 0& 0 \\
              2 & 2 & 0& 0 & 0& 0 \\
              3 & 3 & 3& 3 & 0& 0 \\
              4 & 4 & 4& 4 & 0& 0
            \end{matrix}\right],\quad T=\left[\begin{matrix}
              0 & 0 & 0& 0& 0& 0 \\
              0 & 0 & 0& 0& 0& 0 \\
              1 & 1 & 0& 0& 0& 0 \\
              2 & 2 & 2& 0& 0& 0 \\
              3 & 3 & 3& 0& 0& 0 \\
              4 & 4 & 4& 4& 4& 0
            \end{matrix}\right],\quad U=\left[\begin{matrix}
              0 & 0 & 0& 0& 0& 0 \\
              0 & 0 & 0& 0& 0& 0 \\
              1 & 1 & 0& 0& 0& 0 \\
              2 & 2 & 2& 0& 0& 0 \\
              3 & 3 & 3& 3& 0& 0 \\
              4 & 4 & 4& 4& 4& 0
            \end{matrix}\right],
  \end{equation*}
  Note that they all are in the $1-$REF. Following Definition \ref{def-Blocks-1REF} we can easily calculate:
  \begin{equation*}
  \begin{array}{c}
    \mathbb{W}(S)=(3,5),
  \end{array}\quad \begin{array}{c}
    \mathbb{W}(T)=(3,4,6),
  \end{array}\quad \begin{array}{c}
    \mathbb{W}(U)=(3,4,5,6)
  \end{array}
  \end{equation*}
  With this information we can visualize the corresponding bricks:

  \begin{equation*}
    \begin{array}{ccccc}
      B_1(S)={S_{[1,2]}^{[3,4]}}, & {B_2(S)=S_{[3,4]}^{[5,6]}}, & \quad & \quad & \quad \\
      \quad & \quad & \quad & \quad & \quad \\
       B_1(T)={T_{[1,2]}^{[3]}}, & {B_2(T)=T_{[3]}^{[4,5]}}, & B_3(T)={T_{[4,5]}^{[6]}}, & \quad & \quad \\
      \quad & \quad & \quad & \quad & \quad \\
       B_1(U)={U_{[1,2]}^{[3]}}, & {B_2(U)=U_{[3]}^{[4]}}, & B_3(U)={U_{[4]}^{[5]}},  &  {B_3(U)=U_{[5]}^{[6]}}.  & \quad
    \end{array}
  \end{equation*}

  \begin{equation*}
    S=\left[\begin{matrix}
              0 & 0 & 0& 0& 0& 0 \\
              0 & 0 & 0& 0& 0& 0 \\
              {\color{blue}1} & {\color{blue}1} & 0& 0 & 0& 0 \\
              {\color{blue}2} & {\color{blue}2} & 0& 0 & 0& 0 \\
              3 & 3 & {\color{red}3}& {\color{red}3} & 0& 0 \\
              4 & 4 & {\color{red}4}& {\color{red}4} & 0& 0
            \end{matrix}\right],\quad T=\left[\begin{matrix}
              0 & 0 & 0& 0& 0& 0 \\
              0 & 0 & 0& 0& 0& 0 \\
              {\color{blue}1} & {\color{blue}1}& 0& 0& 0& 0 \\
              2 & 2 & {\color{red}2}& 0& 0& 0 \\
              3 & 3 & {\color{red}3}& 0& 0& 0 \\
              4 & 4 & 4& {\color{green}4}& {\color{green}4}& 0
            \end{matrix}\right],\quad U=\left[\begin{matrix}
              0 & 0 & 0& 0& 0& 0 \\
              0 & 0 & 0& 0& 0& 0 \\
              {\color{blue}1} & {\color{blue}1} & 0& 0& 0& 0 \\
              2 & 2 & {\color{red}2}& 0& 0& 0 \\
              3 & 3 & 3& {\color{orange}3}& 0& 0 \\
              4 & 4 & 4& 4& {\color{green}4}& 0
            \end{matrix}\right],
  \end{equation*}
\end{example}

Note that when $S$ has more than one brick, there is an area of $S$ below the bricks, we call it the \emph{dark side} of $S.$


  \begin{lemma}\label{lemma-Wall-eq-0}
    Let $T$ be a \nSLTM, that is in $1-REF$. Then $T\sim 0_n$ implies that $T=0_n.$
  \end{lemma}

  \begin{proof}
     Is a direct consequence of Theorem \ref{theo-zero-class-n-gen} and  Proposition \ref{Prop-Agorithm-1REF}.
  \end{proof}

In Lemma \ref{lemma-Wall-eq-0}, the case of matrices $T$ in the class $\overline{0_n}$ is completely described. For the rest of this section, we declare that all the matrices mentioned in Lemmas and Propositions are outside of the class $\overline{0_n}.$ We are not going to mention this assumption in the statements of the Lemmas or Propositions of this section.

\begin{definition}\label{def-measure-mat}
Given any matrix $M=[m_{ij}],$ we denote by $\mu(M)=|\{(i,j):m_{ij}\neq 0\}|.$ We call $\mu(M)$ the \emph{measure} of the matrix $M.$
\end{definition}

\begin{lemma}\label{lemma-r1-geq-3}
Let $T$ be a \nSLTM, that is in $1-$REF. Let $\mathbb{W}(T)=(r_1,\dots,r_h)$ its wall. The following assertions are true:
\begin{enumerate}
  \item $r_1\geq 3.$
  \item If $r\in [r_1,r_2[,$ then $\mu(T^{[r]})\geq2.$
\end{enumerate}

\end{lemma}
\begin{proof}
\begin{enumerate}
  \item Since $T$ is a \nSLTM, that is in $1-$REF, necessarily  $T^{[r]}=0$ for $r=1,2.$ In particular we have $r_1\geq 3.$
  \item If we assume that $\mu(T^{[r]})=0,$ that is, the $r-$th row of $T$ is a zero row, then by Definition \ref{def-Blocks-1REF}, we have $r<r_1.$

  On the other hand, if we assume that $\mu(T^{[r]})=1,$ for some $r\in[r_1,r_2[,$ then the leader of this row, $t_{r,c_r}$ is the only nonzero entry of that row. By Definition \ref{def-Blocks-1REF}, the leader of the $r-$th row of $T$ belongs to $B_1(T),$ therefore $c_r<r_1.$ If $c_r=1$ we can apply the operation $\calQ_{(r,1)}(T,\frac{t_{r1}}{2})$ and make this entry to be equal to zero, but this contradicts our assumption, that $T$ is already in the $1-$REF. Otherwise, if $1<c_r<r_1,$ then for each $i\in[1,c_r[$ we will have $\Delta_{i,c_r,r}^{(2)}=2t_{ri}+t_{r,c_r}t_{c_r,i}=0,$ but this again contradicts our assumption that $T$ is already in the $1-$REF (see Definition \ref{def-first-red-form}).
\end{enumerate}
\end{proof}

\begin{proposition}\label{prop-Gamma-triang-1}
  Let $T,S$ be two nonzero \nSLTM, such that:
   \begin{itemize}
     \item $T$ and $S$ are in $1-$REF
     \item $\mathbb{W}(T)=\mathbb{W}(S)=(r_1,\dots,r_h).$
   \end{itemize}
    If $\bgamma:\calA(T)\rightarrow\calA(S)$ is an isomorphism, then the associated matrix $\Gamma=[\gamma_{ij}]$ satisfies the following relations:
  \begin{enumerate}
    \item If $r\in [r_j, r_{j+1}[$ then $\gamma_{kr}=0$ for all $k\geq r_{j+1}.$
    \item There is a bijection $\rrn:[1,n]\rightarrow [1,n]$ such that:
    \begin{enumerate}
      \item For each $j\in[1,h]$ we have $\rrn([r_j,r_{j+1}[)= [r_j,r_{j+1}[$
      \item For each $i,k\in [r_j,r_{j+1}[$ we have that
    \begin{equation*}
       \gamma_{ik}\neq 0 \quad \textrm{if and only if}\quad k=\rrn(i).
    \end{equation*}
    \end{enumerate}
  \end{enumerate}
\end{proposition}

The first assertion described in Proposition \ref{prop-Gamma-triang-1}, basically say that the matrix $\Gamma$ is \emph{upper triangular by blocks}:

\begin{equation}\label{eq0-Gamma-triangular}
          \Gamma=\left[\begin{matrix}
                         \Gamma_{[1,r_1[}^{[1,r_1[} & \Gamma_{[r_1,r_2[}^{[1,r_1[} & \cdots& \Gamma_{[r_{p-1},r_p[}^{[1,r_1[}&\Gamma_{[r_{p},r_{p+1}[}^{[1,r_1[}
                         &\cdots & \Gamma_{[r_{h},n]}^{[1,r_1[}\\
                         \quad & \quad & \quad & \quad & \quad& \quad & \quad\\
                         0 & \Gamma_{[r_1,r_2[}^{[r_1,r_2[} & \cdots& \Gamma_{[r_{p-1},r_p[}^{[r_1,r_2[}  &\Gamma_{[r_{p},r_{p+1}[}^{[r_1,r_2[}&\cdots & \Gamma_{[r_{h},n]}^{[r_1,r_2[} \\
                         \quad & \quad & \quad & \quad& \quad& \quad & \quad \\
                         \vdots & \vdots & \ddots & \vdots & \vdots&\cdots & \vdots\\
                         \quad & \quad & \quad & \quad&  \quad & \quad & \quad\\
                         0 & 0 & \cdots& \Gamma_{[r_{p-1},r_p[}^{[r_{p-1},r_p[} &\Gamma_{[r_{p},r_{p+1}[}^{[r_{p-1},r_p[}&\cdots & \Gamma_{[r_{h},n]}^{[r_{p-1},r_p[}  \\
                          \quad & \quad & \quad & \quad& \quad& \quad & \quad \\
                          0 & 0 & \cdots& 0 &\Gamma_{[r_{p},r_{p+1}[}^{[r_{p},r_{p+1}[}&\cdots & \Gamma_{[r_{h},n]}^{[r_p,r_{p+1}[}\\
                          \quad & \quad & \quad & \quad& \quad& \quad & \quad \\
                          0 & 0 & \cdots& 0 &0&\ddots & \vdots\\
                          \quad & \quad & \quad & \quad& \quad& \quad & \quad \\
                          0 & 0 & \cdots& 0 &0&\cdots & \Gamma_{[r_{h},n]}^{[r_{h},n]}\\
                       \end{matrix}\right]
        \end{equation}
          The second assertion implies that in each column of each block in the diagonal, there is only one entry different from zero, that is $\mu\left(\Gamma_{[r]}^{[r_{p},r_{p+1}[}\right)=1$ for each $r\in [r_p,r_{p+1}[$ and each $p=1,\dots,h.$ Equivalently we have $\mu\left(\Gamma_{[r_{p},r_{p+1}[}^{[r]}\right)=1$ for each $r\in [r_p,r_{p+1}[$ and each $p=1,\dots,h.$

\begin{proof}
  First we recall that the entries of the matrices $T,S$ and $\Gamma$ are all related by equation \ref{Eq1-theo-Key-condition} in Theorem \ref{theo-Key-condition}
  \begin{itemize}
    \item If $h=1$ then
    \begin{equation}
      T=\left[\begin{matrix}
                0 & 0\\
                B_1(T) & 0
              \end{matrix}\right],\quad S=\left[\begin{matrix}
                0 & 0\\
                B_1(S) & 0
              \end{matrix}\right].
    \end{equation}
    \begin{itemize}
     \item If we take $r\in[1,r_1[,$ then $t_{rj}=0$ for each $j,$ then equation \ref{Eq1-theo-Key-condition} reduces to:
      \begin{equation}\label{eq2-key-recall}
          (2\gamma_{ir}+s_{ki}\gamma_{kr})\gamma_{kr}
          =0,\quad (1\leq i<k\leq n,r\in[1,r_1[)
      \end{equation}
      If we take $1\leq i_1<i_2<r_1$ we also have $s_{i_2i_1}=0,$ then we obtain:
      \begin{equation}\label{eq3-key-recall}
          2\gamma_{i_1r}\gamma_{i_2r}=0\quad (1\leq i_1<i_2<r_1,r\in[1,r_1[)
      \end{equation}
       Now if we take $k\geq r_1,$ by lemma \ref{lemma-r1-geq-3}, we have $\mu(S^{[k]})\geq 2,$ therefore we can pick $i_1,i_2\in [1,r_1[$ such that $s_{ki_1}\neq 0$ and $s_{ki_2}\neq 0.$
       Now if we replace respectively $i_1,i_2\in [1,r_1[$ and $k\geq r_1$ in equation \ref{eq2-key-recall}, we obtain
      \begin{equation}\label{eq4-key-recall}
       \left\{\begin{array}{c}
               (2\gamma_{i_1r}+s_{ki_1}\gamma_{kr})\gamma_{kr}=0 \\
               (2\gamma_{i_2r}+s_{ki_2}\gamma_{kr})\gamma_{kr}=0
             \end{array}\right.
      \end{equation}
      then equation \ref{eq3-key-recall} implies that $\gamma_{kr}=0$ for $k\geq r_1,$ that is
      \begin{equation}\label{eq1-Gamma-triangular}
       \Gamma=\left[\begin{matrix}
                      \Gamma_{[1,r_1[}^{[1,r_1[} & \Gamma_{[r_1,n]}^{[1,r_1]} \\
                      \quad&\quad\\
                      0 & \Gamma_{[r_1,n]}^{[r_1,n]}
                    \end{matrix}\right].
      \end{equation}
       By the above result, we have that
     \begin{equation}\label{eq2-Gamma-triangular}
      \det(\Gamma)=\det\left(\Gamma_{[1,r_1[}^{[1,r_1[}\right)
      \det\left(\Gamma_{[r_1,n]}^{[r_1,n]}\right)
     \end{equation}
     Since $\det(\Gamma)\neq 0$ and equation \ref{eq3-key-recall} holds, we conclude that there is a bijection $\rrn:[1,r_1[\mapsto[1,r_1[$ such that for each $i,k\in[1,r_1[$ we have that $\gamma_{ik}\neq 0$ if and only if $k=\rrn(i).$
     \item If we take $r\in[r_1,n],$ the equation \ref{Eq1-theo-Key-condition} takes the form:
         \begin{equation}\label{eq5-key-recall}
          2\gamma_{ir}\gamma_{kr}+s_{ki}\gamma_{kr}^2
          =\sum_{j<r_1}t_{rj}\left(\gamma_{kj}\gamma_{kr}s_{ki}+\gamma_{kj}\gamma_{ir}+
          \gamma_{ij}\gamma_{kr}\right)
     \end{equation}

      Now if we replace $i,k\in[r_1,n]$ with $i<k$ in equation \ref{eq5-key-recall}, we obtain:
         \begin{equation}\label{eq6-key-recall}
           2\gamma_{ir}\gamma_{kr}=0,\quad (r_1\leq i<k\leq n,r\in [r_1,n])
         \end{equation}
         since $s_{ki}=0$ and by the above case, $\gamma_{kj}=\gamma_{ij}=0$ for each $j<r_1.$ Therefore by equation \ref{eq2-Gamma-triangular} and since $\det(\Gamma)\neq 0,$ we obtain that the bijection $\rrn:[1,r_1[\rightarrow[1,r_1[$ can be extended to a bijection $\rrn:[1,n]\rightarrow[1,n]$ satisfying items (2.a) and (2.b).
     \end{itemize}
    \item If $h>1$ we have:
    \begin{equation}\label{eq1-two-SLTM-by-blocks}
      T=\left[\begin{matrix}
                0 & 0 & \cdots & 0&0 \\
                B_1(T) & 0 &\cdots & 0&0 \\
                \bullet & B_{2}(T)&\cdots & 0&0\\
                \vdots & \vdots &\ddots&\vdots&\vdots\\
                \bullet & \bullet & \bullet& B_{h}(T)&0
              \end{matrix}\right],\quad
      S=\left[\begin{matrix}
                0 & 0 & \cdots & 0&0 \\
                B_1(S) & 0 &\cdots & 0&0 \\
                \bullet & B_{2}(S)&\cdots & 0&0\\
                \vdots & \vdots &\ddots&\vdots&\vdots\\
                \bullet & \bullet & \bullet& B_{h}(S)&0
              \end{matrix}\right]
    \end{equation}
    \begin{itemize}
      \item If we take $r\in[1,r_1[$ the analysis is basically the same as the case $h=1.$ First note that equation \ref{Eq1-theo-Key-condition} reduces to equation \ref{eq2-key-recall}, then if we replace $1\leq i<k<r_1$ in this equation, we obtain \ref{eq3-key-recall}. Now if we take $k\in[r_1,r_2[$ then since $\mu(S^{[k]})\geq 2$ we obtain analogously to the case $h=1$ that $\gamma_{kr}=0.$

          Suppose that we already have prove for $1\leq j<h$ and $i\in [r_j,r_{j+1}[$ that $\gamma_{ir}=0,$ then take $k_0\in [r_{j+1},r_{j+2}[.$ By definition of the bricks, we have that the leader $s_{k_0,i_0}$ of the $k_0-$th row of $S$ is in the brick $B_{j+1}(S),$ that is $i_0\in [r_j,r_{j+1}[.$ Replacing $k=k_0,i=i_0$ and $r\in [1,r_1[$ in equation \ref{eq2-key-recall}, we obtain
          \begin{equation*}
            (2\gamma_{i_0r}+s_{k_0i_0}\gamma_{k_0r})\gamma_{k_0r}=0
          \end{equation*}
          but by our hypothesis, we have $\gamma_{i_0r}=0,$ and we conclude that $\gamma_{k_0r}=0$ (since $s_{k_0i_0}\neq 0$). By induction we conclude that $\gamma_{kr}=0$ whenever $k\geq r_1$ and $r\in [1,r_1[.$ Therefore we have that equations \ref{eq1-Gamma-triangular} and \ref{eq2-Gamma-triangular} hold in this case. The same analysis that we did in the case $h=1$ allows us to conclude that there is a bijection $\rrn_0:[1,r_1[\rightarrow[1,r_1[$ such that for every pair $i,k\in [1,r_1[$
          \begin{equation*}
            \gamma_{ik}\neq 0\quad\textrm{iff}\quad k=\rrn_0(i)
          \end{equation*}

      \item Let $p\in [1,h[$ and assume that we already have proved for each $q<p$ that:
       \begin{enumerate}
         \item $
            \gamma_{kr}=0\quad\textrm{whenever}\quad k\geq r_{q+1}\wedge r\in[r_q,r_{q+1}[$
         \item There is a bijection $\rrn_q:[r_q,r_{q+1}[\rightarrow [r_q,r_{q+1}[$ such that, for every pair $i,k\in [r_q,r_{q+1}[$
        \begin{equation*}
        \gamma_{ik}\neq 0\quad\textrm{iff}\quad k=\rrn_q(i)
        \end{equation*}
       \end{enumerate}
        Therefore, at this point, $\Gamma$ looks like this:
        \begin{equation}\label{eq3-Gamma-triangular}
          \Gamma=\left[\begin{matrix}
                         \Gamma_{[1,r_1[}^{[1,r_1[} & \Gamma_{[r_1,r_2[}^{[1,r_1[} & \cdots& \Gamma_{[r_{p-1},r_p[}^{[1,r_1[}&\Gamma_{[r_{p},r_{p+1}[}^{[1,r_1[}
                         &\cdots & \Gamma_{[r_{h},n]}^{[1,r_1[}\\
                         \quad & \quad & \quad & \quad & \quad& \quad & \quad\\
                         0 & \Gamma_{[r_1,r_2[}^{[r_1,r_2[} & \cdots& \Gamma_{[r_{p-1},r_p[}^{[r_1,r_2[}  &\Gamma_{[r_{p},r_{p+1}[}^{[r_1,r_2[}&\cdots & \Gamma_{[r_{h},n]}^{[r_1,r_2[} \\
                         \quad & \quad & \quad & \quad& \quad& \quad & \quad \\
                         \vdots & \vdots & \ddots & \vdots & \vdots&\cdots & \vdots\\
                         \quad & \quad & \quad & \quad&  \quad & \quad & \quad\\
                         0 & 0 & \cdots& \Gamma_{[r_{p-1},r_p[}^{[r_{p-1},r_p[} &\Gamma_{[r_{p},r_{p+1}[}^{[r_{p-1},r_p[}&\cdots & \Gamma_{[r_{h},n]}^{[r_{p-1},r_p[}  \\
                          \quad & \quad & \quad & \quad& \quad& \quad & \quad \\
                          0 & 0 & \cdots& 0 &\Gamma_{[r_{p},r_{p+1}[}^{[r_{p},r_{p+1}[}&\cdots & \Gamma_{[r_{h},n]}^{[r_p,r_{p+1}[}\\
                         \quad & \quad & \quad & \quad& \quad& \quad & \quad \\
\vdots & \vdots & \vdots& \vdots&\vdots &\vdots &\vdots\\
\quad & \quad & \quad & \quad& \quad& \quad & \quad \\
                          0 & 0 & \cdots& 0 &\Gamma_{[r_{p},r_{p+1}[}^{[r_{p+1},n]}&\cdots & \Gamma_{[r_{h},n]}^{[r_{p+1},n]}
                       \end{matrix}\right]
        \end{equation}
        Define for $i\in[1,r_p[$

        \begin{equation}\label{eq1-extending-rrn}
         \rrn(i)=\rrn_q(i),\quad \textrm{whenever}\quad i\in [r_q,r_{q+1}[
        \end{equation}
         Now take a $k_0\in [r_{p+1},r_{p+2}[$ and suppose that $s_{k_0,i_0}\neq 0$ is the leader of the $k_0-$th row of $S.$ By definition of the bricks, $i_0\in [r_{p},r_{p+1}[,$ then if we replace $k=k_0,i=i_0$ and $r\in[r_{p},r_{p+1}[$ in equation \ref{Eq1-theo-Key-condition}, we obtain
        \begin{equation}\label{eq7-key-recall}
          (2\gamma_{i_0,r}+s_{k_0,i_0}\gamma_{k_0,r})\gamma_{k_0,r}
          =\sum_{j<r_p}t_{rj}(\gamma_{k_0j}\gamma_{k_0r}s_{k_0i_0}+\gamma_{k_0j}
          \gamma_{i_0r}+\gamma_{i_0j}\gamma_{k_0r})
          =0,
        \end{equation}
        since by our hypothesis, $\gamma_{i_0,j}=\gamma_{k_0,j}=0$ for each $j<r_p.$ If we suppose that $\gamma_{k_0,r}\neq 0$ then equation \ref{eq7-key-recall} implies that $\gamma_{i_0,r}\neq 0.$ Moreover we have
        \begin{equation}\label{eq8-key-recall}
          \gamma_{i_0,r}=-\frac{s_{k_0,i_0}}{2} \gamma_{k_0,r}
        \end{equation}
       If we replace $i\in [1,i_0[,r\in[r_p,r_{p+1}[$ in equation \ref{Eq1-theo-Key-condition}, for $k=k_0$ and $k=i_0$ respectively, we obtain
       \begin{equation}\label{eq9-key-recall}
         \left\{\begin{array}{c}
                  2\gamma_{ir}+s_{k_0i}\gamma_{k_0r}-b_{ri}=0 \\
                  \quad\\
                  2\gamma_{ir}+s_{i_0i}\gamma_{i_0r}-b_{ri}=0
                \end{array}\right.
       \end{equation}
       where $b_{ri}:=\sum_{j<r}t_{rj}\gamma_{ij}.$ Therefore we have

       \begin{equation}\label{eq10-key-recall}
         s_{k_0i}\gamma_{k_0r}=s_{i_0i}\gamma_{i_0r}=-\frac{s_{k_0,i_0}s_{i_0,i}}{2} \gamma_{k_0,r}
       \end{equation}
       (the last equation comes from equation \ref{eq8-key-recall}). Finally, since $\gamma_{k_0r}\neq 0$ we conclude that $\Delta_{i,i_0,k_0}^{(2)}(S)=0$ for each $i\in [1,i_0[$ and we obtain a contradiction on our assumption that $S$ is in the $1-$REF. Therefore we have that $\gamma_{k_0r}=0,$ for any $k_0\in [r_{p+1},r_{p+2}[.$


       Let $l\in[p+2,h]$ and assume that we have already proved that $\gamma_{ir}=0$ whenever $i\in [r_{p+1},r_{l}[$ and $r\in [r_{p},r_{p+1}[.$ Take any $k_1\in [r_l,r_{l+1}[$ and assume that $s_{k_1i_1}$ is the leader of the $k_1-$th row of $S,$ then $s_{k_1i_1}\neq 0$ and $i_1\in [r_{l-1},r_{l}[,$ then when we replace $i=i_1,k=k_1$ and $r\in [r_{p},r_{p+1}[$ in equation \ref{Eq1-theo-Key-condition}, we obtain
       \begin{equation}\label{eq12-key-recall}
          (2\gamma_{i_1,r}+s_{k_1,i_1}\gamma_{k_1,r})\gamma_{k_1,r}
          =0,
        \end{equation}
       but this equation implies that $\gamma_{k_1,r}=0,$ since $s_{k_1i_1}\neq 0$ and by our hypothesis $\gamma_{i_1,r}=0.$ By induction we conclude that $\gamma_{kr}=0,$ whenever $k\geq r_{p+1}$ and $r\in [r_p,r_{p+1}[.$

       In particular we have

        \begin{equation}\label{eq5-Gamma-triangular}
          \det(\Gamma)=\left(
          \prod_{q=0}^{p}\det\left(\Gamma_{[r_q,r_{q+1}[}^{[r_q,r_{q+1}[}\right)\right)
          \det\left(\Gamma_{[r_{p+1},n]}^{[r_{p+1},n]}\right)
        \end{equation}
        and since $\det(\Gamma)\neq 0,$ in particular we have
        \begin{equation}\label{eq6-Gamma-triangular}
        \det\left(\Gamma_{[r_p,r_{p+1}[}^{[r_p,r_{p+1}[}\right)\neq 0.
        \end{equation}

        Finally if we replace $i,k,r\in[r_p,r_{p+1}[$ (with $i<k$) in equation \ref{Eq1-theo-Key-condition}, we obtain:
        \begin{equation}\label{eq13-key-recall}
          2\gamma_{ir}\gamma_{kr}=0
        \end{equation}
        since $s_{ki}=0$ and $\gamma_{ij}=\gamma_{kj}=0$ for every $j\in [1,r_p[.$ Then equation \ref{eq6-Gamma-triangular} and \ref{eq13-key-recall} imply that there is  bijection $\rrn_{p}:[r_p,r_{p+1}[\rightarrow[r_p,r_{p+1}[$ such that for every $i,k\in[r_p,r_{p+1}[$
        \begin{equation*}
          \gamma_{ik}\neq 0\quad \textrm{iff}\quad k=\rrn_{p}(i),
        \end{equation*}
        that is, the bijection $\rrn:[1,r_p[\rightarrow [1,r_p[$ given in equation \ref{eq1-extending-rrn}, can be extended to a bijection $\rrn:[1,r_{p+1}[\rightarrow [1,r_{p+1}[,$ where $\rrn(i)=\rrn_p(i),$ whenever $i\in [r_p,r_{p+1}[.$

     The above recursive argument implies that all the assertions of Proposition \ref{prop-Gamma-triang-1} are true.
    \end{itemize}
  \end{itemize}
\end{proof}

\begin{example}\label{ex1-prop-Gamma-triang1}
Consider the following matrices
  \begin{equation}\label{eq1-ex1-prop-Gamma-triang1}
    T=\left[\begin{matrix}
              0 & 0 & 0 & 0 \\
              0 & 0 & 0 & 0 \\
              t_{31} & t_{32} & 0 & 0 \\
              t_{41} & t_{42} & 0 & 0
            \end{matrix}\right],\quad
    S=\left[\begin{matrix}
              0 & 0 & 0 & 0 \\
              0 & 0 & 0 & 0 \\
               s_{31} & s_{32} & 0 & 0 \\
              s_{41} & s_{42} & 0 & 0
            \end{matrix}\right].
  \end{equation}
  and assume that $\mathbb{W}(T)=\mathbb{W}(S)=(3).$ If $\bgamma:\calA(T)\rightarrow \calA(S)$ is an isomorphism, then Proposition \ref{prop-Gamma-triang-1} says that the associated matrix $\Gamma=[\gamma_{ij}]$ and the corresponding bijection $\rrn$ have one of the following forms:
   \begin{equation*}
     \Gamma_1=\left[\begin{matrix}
                    \gamma_{11} & 0 &  \gamma_{13} &  \gamma_{14} \\
                    0 &  \gamma_{22} &  \gamma_{23} &  \gamma_{24} \\
                    0 & 0 &  \gamma_{33} & 0 \\
                    0 & 0 & 0 &  \gamma_{44}
                  \end{matrix}\right],\quad
     \rrn:\left\{\begin{array}{c}
                   1\mapsto 1 \\
                   2\mapsto 2 \\
                    3\mapsto 3\\
                   4\mapsto 4
                 \end{array}\right.
   \end{equation*}
   \begin{equation*}
   \Gamma_2=\left[\begin{matrix}
                    0 &  \gamma_{12} &  \gamma_{13} &  \gamma_{14} \\
                     \gamma_{21} & 0 &  \gamma_{23} &  \gamma_{24} \\
                    0 & 0 &  \gamma_{33} & 0 \\
                    0 & 0 & 0 &  \gamma_{44}
                  \end{matrix}\right],\quad
     \rrn:\left\{\begin{array}{c}
                   1\mapsto 2 \\
                   2\mapsto 1 \\
                    3\mapsto 3\\
                   4\mapsto 4
                 \end{array}\right.
   \end{equation*}
   \begin{equation*}
   \Gamma_3=\left[\begin{matrix}
                     \gamma_{11} & 0 &  \gamma_{13} &  \gamma_{14} \\
                    0 &  \gamma_{22} &  \gamma_{23} &  \gamma_{24} \\
                    0 & 0 & 0 &  \gamma_{34} \\
                    0 & 0 &  \gamma_{43} & 0
                  \end{matrix}\right],\quad
     \rrn:\left\{\begin{array}{c}
                   1\mapsto 1 \\
                   2\mapsto 2 \\
                    3\mapsto 4\\
                   4\mapsto 3
                 \end{array}\right.
   \end{equation*}

   \begin{equation*}
   \Gamma_4=\left[\begin{matrix}
                    0 &  \gamma_{12} &  \gamma_{13} &  \gamma_{14} \\
                     \gamma_{21} & 0 &  \gamma_{23} &  \gamma_{24} \\
                    0 & 0 & 0 &  \gamma_{34}\\
                    0 & 0 &  \gamma_{43} & 0
                  \end{matrix}\right],\quad
     \rrn:\left\{\begin{array}{c}
                   1\mapsto 2 \\
                   2\mapsto 1 \\
                    3\mapsto 4\\
                   4\mapsto 3
                 \end{array}\right.
   \end{equation*}
\end{example}

\begin{example}\label{ex2-prop-Gamma-triang1}
  Consider the following matrices:
  \begin{equation}\label{eq1-ex2-prop-Gamma-triang1}
    T=\left[\begin{matrix}
              0 & 0 & 0 & 0 \\
              0 & 0 & 0 & 0 \\
              t_{31} & t_{32} & 0 & 0 \\
              t_{41} & t_{42} & t_{43}  & 0
            \end{matrix}\right],\quad
   S=\left[\begin{matrix}
             0 & 0 & 0 & 0 \\
              0 & 0 & 0 & 0 \\
              s_{31} & s_{32} & 0 & 0 \\
              s_{41} & s_{42} & s_{43}  & 0
           \end{matrix}\right],
  \end{equation}
  and assume that $\mathbb{W}(T)=\mathbb{W}(S)=(3,4).$ If $\bgamma:\calA(T)\rightarrow\calA(S)$ is an isomorphism, then Proposition \ref{prop-Gamma-triang-1} says that the associated matrix $\Gamma=[\gamma_{ij}]$ and the corresponding bijection $\rrn$ have one of the following forms:
\begin{equation*}
  \Gamma_1=\left[\begin{matrix}
                 \gamma_{11} & 0 & \gamma_{13} & \gamma_{14} \\
                 0 & \gamma_{22} & \gamma_{23} & \gamma_{24} \\
                 0 & 0 & \gamma_{33} & \gamma_{34} \\
                 0 & 0 & 0 & \gamma_{44}
               \end{matrix}\right],\quad
  \rrn:\left\{\begin{array}{c}
                1\mapsto 1 \\
                2\mapsto 2 \\
                3\mapsto 3 \\
                4\mapsto 4
              \end{array}\right.
\end{equation*}

\begin{equation*}
\Gamma_2=\left[\begin{matrix}
                 0 & \gamma_{12} & \gamma_{13} & \gamma_{14} \\
                 \gamma_{21} & 0 & \gamma_{23} & \gamma_{24} \\
                 0 & 0 & \gamma_{33} & \gamma_{34} \\
                 0 & 0 & 0 & \gamma_{44}
               \end{matrix}\right],\quad
  \rrn:\left\{\begin{array}{c}
                1\mapsto 2 \\
                2\mapsto 1 \\
                3\mapsto 3 \\
                4\mapsto 4
              \end{array}\right.
\end{equation*}
\end{example}

Note that in both, example \ref{ex1-prop-Gamma-triang1} and example \ref{ex2-prop-Gamma-triang1}, in principle, to find out if there exist an isomorphisms $\bgamma:\calA(T)\rightarrow\calA(S),$ we have to study the system of equations given in equation \ref{Eq1-theo-Key-condition}, which is composed by $16$ unknowns variables, but with the information provided in Proposition \ref{prop-Gamma-triang-1}, we obtain a reduction of those variables to just $8$ of them.

\begin{lemma}\label{lemma-r1-fails}
  Let $T,S$ be two \nSLTM, such that both are in $1-$REF, $\mathbb{W}(T)=(r_1,\dots,r_h)$ and $\mathbb{W}(S)=(r_1',\dots,r_g').$ If $r_1\neq r_1'$ then $T\nsim S.$
\end{lemma}
\begin{proof}
  We can assume, without loss of generality that  $r_1>r_1'.$

  If we suppose that there is an isomorphism $\bgamma:\calA(T)\rightarrow \calA(S),$ then by Theorem \ref{theo-Key-condition}, the associated matrix $\Gamma=[\gamma_{ij}]$ satisfies equation \ref{Eq1-theo-Key-condition}.

  If we replace $i,k\in[1,r_1'[$ and $r\in[1,r_1[$ in equation \ref{Eq1-theo-Key-condition}, we obtain:
  \begin{equation}\label{eq1-proof-lemma-r1-fails}
    2\gamma_{ir}\gamma_{kr}=0.
  \end{equation}
  If we replace $i\in[1,r_1'[, k\in [r_1',r_2'[$ and  $r\in[1,r_1[$ in equation \ref{Eq1-theo-Key-condition}, we obtain:

  \begin{equation}\label{eq2-proof-lemma-r1-fails}
    (2\gamma_{ir}+s_{ki}\gamma_{kr})\gamma_{kr}=0.
  \end{equation}
  By Lemma \ref{lemma-r1-geq-3} $\mu(S^{[k]})\geq 2,$ therefore there are at least two (different) indices $i_1,i_2\in[1,r_1'[$ such that $s_{ki_1}\neq 0$ and $s_{ki_2}\neq 0.$ By equation \ref{eq1-proof-lemma-r1-fails} $\gamma_{i_1r}=0$ or $\gamma_{i_2r}=0,$ then equation \ref{eq2-proof-lemma-r1-fails} implies that $\gamma_{kr}=0.$

  Assuming that we have prove that $\gamma_{ir}=0$ for $i\in [r_p',r_{p+1}'[$ and $r\in [1,r_1[$ then if we replace $k\in [r_{p+1}',r_{p+2}'[,i=c_k$ and $r\in [1,r_1[$ in equation \ref{Eq1-theo-Key-condition}, we obtain

  \begin{equation}\label{eq3-proof-lemma-r1-fails}
    (2\gamma_{ir}+s_{ki}\gamma_{kr})\gamma_{kr}=s_{ki}\gamma_{kr}^{2}=0.
  \end{equation}
   Since  $s_{k,c_k}$ is the leader of the $k-$th row of $S$ (see Subsection \ref{ssec-Notations} for notations), therefore $s_{k,c_k}\neq 0$ and $c_k\in [r_p',r_{p+1}'[$ then $\gamma_{c_kr}=0$ and we conclude that $\gamma_{kr}=0$ for all $k\in [r_{p+1}',r_{p+2}'[.$ By induction on $p$ we conclude that $\gamma_{kr}=0$ for all $k\geq r_1'$ and all $r\in [1,r_1[.$  Now by our assumption $r_1>r_1',$ then we conclude that $\det(\Gamma)=0$ and that is a contradiction with Theorem \ref{theo-Key-condition}
\end{proof}

Note that in Lemma \ref{lemma-r1-fails}, the number of bricks (the length of the Wall) is not relevant, we only have to see the first component of the corresponding walls.  Lemma \ref{lemma-r1-fails} implies that the number of zero rows (resp. the number of nonzero rows) of matrices in the $1-$REF is an invariant for the isomorphism classes of objects on the category $\bcalNT,$ therefore they deserve to have a name:

\begin{definition}\label{def-Tnulity}
  Let $T$ be a \SLTM, that is in $1-$REF. We define:
  \begin{enumerate}
    \item The \emph{Triangular nullity} of $T$ as the number of zero rows of $T.$ We denote that number by $\textrm{TNul}(T).$
    \item The \emph{Triangular rank} of $T$ as the number of non zero rows of the matrix $T.$ We denote that number by $\textrm{TRank}(T).$
  \end{enumerate}
\end{definition}

\begin{corollary}\label{coro-TNul-fails}
  Let $T,S$ be two \nSLTM, such that, both are in $1-$REF:
  \begin{enumerate}
    \item If $\textrm{TNul}(T)\neq \textrm{TNul}(S),$ then $T\nsim S.$
    \item If $\textrm{TRank}(T)\neq \textrm{TRank}(S),$ then $T\nsim S.$
  \end{enumerate}
\end{corollary}
\begin{proof}
  It follows directly from Lemma \ref{lemma-r1-fails}.
\end{proof}

\begin{example}\label{ex1-Tnul-neq}
  Consider the matrices
  \begin{equation*}
    S=\left[\begin{matrix}
              0 & 0 & 0 & 0 & 0 & 0 \\
              0 & 0 & 0 & 0 & 0 & 0 \\
              0 & 0 & 0 & 0 & 0 & 0 \\
              1 & 1 & 1 & 0 & 0 & 0 \\
              1 & 1 & 1 & 0 & 0 & 0 \\
              1 & 1 & 1 & 0 & 0 & 0
            \end{matrix}\right],\quad
    T=\left[\begin{matrix}
              0 & 0 & 0 & 0 & 0 & 0 \\
              0 & 0 & 0 & 0 & 0 & 0 \\
              0 & 0 & 0 & 0 & 0 & 0 \\
              0 & 0 & 0 & 0 & 0 & 0 \\
              1 & 1 & 1 & 1 & 0 & 0 \\
              1 & 1 & 1 & 1 & 0 & 0
            \end{matrix}\right].
  \end{equation*}
  It is clear that $S$ and $T$ are in the $1-$REF. They both have only one brick, but we can easily see that $\textrm{TNul}(S)=3$ and $\textrm{TNul}(T)=4,$ therefore by Corollary \ref{coro-TNul-fails}, we conclude that $S\nsim T.$
\end{example}

\begin{lemma}\label{lemma-rd-fails}
  Let $T,S$ be two \nSLTM, such that, both are in $1-$REF, $\mathbb{W}(T)=(r_1,\dots,r_h)$ and  $\mathbb{W}(S)=(r_1',\dots,r_{g}'),$ let $m=\min\{h,g\}.$ If there is a $p\in[1,m]$ such that $r_p\neq r_p',$ then $T\nsim S.$
\end{lemma}

\begin{proof}
  By lemma \ref{lemma-r1-fails} we have that if $r_1\neq r_1',$ then $T\nsim S.$ Lets assume that $r_1=r_1'$ and define $p=\min\{j\in[2,h]:r_j\neq r_j'\}.$ We can assume without loss of generality, that $r_p>r_p'.$

  If we assume that $T\sim S,$ then there is an invertible matrix $\Gamma=[\gamma_{ij}]$ satisfying equation \ref{Eq1-theo-Key-condition} of Theorem \ref{theo-Key-condition}. By our hypothesis $r_j=r_j'$ for every $j\in [1,p[,$ then if we restrict the parameter $r$ in equation \ref{Eq1-theo-Key-condition} to the interval $[1,r_p'[$ we can apply the same argument as we used in the proof of Proposition \ref{prop-Gamma-triang-1} to conclude that the submatrix $\Gamma_{[1,r_p'[}^{[1,n]}$  satisfies the following relations:
  \begin{enumerate}
    \item If $r\in [r_j, r_{j+1}[$ and $j\in[1,p[,$ then $\gamma_{kr}=0$ for all $k\geq r_{j+1}.$
    \item There is a bijection $\rrn:[1,r_p'[\rightarrow [1,r_p'[$ such that:
    \begin{enumerate}
      \item For each $j\in[1,p[$ we have $\rrn([r_j,r_{j+1}[)= [r_j,r_{j+1}[$
      \item For each $j\in[1,p[$ and each $i,k\in [r_j,r_{j+1}[$ we have that
    \begin{equation*}
       \gamma_{ik}\neq 0 \quad \textrm{if and only if}\quad k=\rrn(i).
    \end{equation*}
    \end{enumerate}
  \end{enumerate}
  In particular the matrix $\Gamma$ looks like this:

\begin{equation}\label{eq03-Gamma-triangular}
\Gamma=
\left[\begin{matrix}
\Gamma_{[1,r_1[}^{[1,r_1[} & \Gamma_{[r_1,r_2[}^{[1,r_1[} & \cdots& \Gamma_{[r_{p-1},r_p'[}^{[1,r_1[}&\Gamma_{[r_{p}',r_{p+1}'[}^{[1,r_1[}
&\cdots & \Gamma_{[r_{h}',n]}^{[1,r_1[}\\
\quad & \quad & \quad & \quad & \quad& \quad & \quad\\
0 & \Gamma_{[r_1,r_2[}^{[r_1,r_2[} & \cdots& \Gamma_{[r_{p-1},r_p'[}^{[r_1,r_2[}  &\Gamma_{[r_{p}',r_{p+1}'[}^{[r_1,r_2[}&\cdots & \Gamma_{[r_{h}',n]}^{[r_1,r_2[} \\
\quad & \quad & \quad & \quad& \quad& \quad & \quad \\
\vdots & \vdots & \ddots & \vdots & \vdots&\cdots & \vdots\\
\quad & \quad & \quad & \quad&  \quad & \quad & \quad\\
0 & 0 & \cdots& \Gamma_{[r_{p-1},r_p'[}^{[r_{p-1},r_p'[} &\Gamma_{[r_{p}',r_{p+1}'[}^{[r_{p-1},r_p'[}&\cdots & \Gamma_{[r_{h}',n]}^{[r_{p-1},r_p'[}  \\
\quad & \quad & \quad & \quad& \quad& \quad & \quad \\
0 & 0 & \cdots& 0 &\Gamma_{[r_{p}',r_{p+1}'[}^{[r_{p}',r_{p+1}'[}&\cdots & \Gamma_{[r_{h}',n]}^{[r_p',r_{p+1}'[}\\
\quad & \quad & \quad & \quad& \quad& \quad & \quad \\
\vdots & \vdots & \vdots& \vdots&\vdots &\vdots &\vdots\\
\quad & \quad & \quad & \quad& \quad& \quad & \quad \\
0 & 0 & \cdots& 0 &\Gamma_{[r_{p}',r_{p+1}'[}^{[r_{h}',n]}&\cdots & \Gamma_{[r_{h}',n]}^{[r_{h}',n]}
\end{matrix}\right]
\end{equation}

  Now if we replace $r=r_p'$ in equation \ref{Eq1-theo-Key-condition}, we obtain
  \begin{equation}\label{eq1-recall-d-KeyEq}
    2\gamma_{ir}\gamma_{kr}+s_{ki}\gamma_{kr}^2=
    \sum_{j<r_{p}'}{t_{rj}\left(\gamma_{kj}\gamma_{kr}s_{ki}+\gamma_{kj}\gamma_{ir}+
    \gamma_{ij}\gamma_{kr}\right)}.
  \end{equation}
  But since $r_p>r_{p}'>r_{p-1},$ we have that $t_{rj}=0$ for each $j\in[r_{p-1},r_p'[.$ Then for $r=r_{p}'$ we have:

  \begin{equation}\label{eq2-recall-d-KeyEq}
    2\gamma_{ir}\gamma_{kr}+s_{ki}\gamma_{kr}^2=
    \sum_{j<r_{p-1}}{t_{rj}\left(\gamma_{kj}\gamma_{kr}s_{ki}+\gamma_{kj}\gamma_{ir}+
    \gamma_{ij}\gamma_{kr}\right)}.
  \end{equation}

  Now we fix a $k_0\in[r_p',r_{p+1}'[$ and take $i_0$ such that $s_{k_0i_0}$ is the leader of the $k_0-$th row of $S.$ If we replace $k=k_0,i=i_0$ in equation \ref{eq2-recall-d-KeyEq}, we obtain (for $r=r_{p}'$):

  \begin{equation}\label{eq3-recall-d-KeyEq}
    (2\gamma_{{i_0}r}+s_{k_0i_0}\gamma_{k_0r})\gamma_{k_0r}=0,
  \end{equation}
  since by definition of the bricks (Definition \ref{def-Blocks-1REF}), we have that  ${i_0}\in [r_{p-1},r_p'[,$ therefore $\gamma_{kj}=\gamma_{ij}=0$ for each $j\in [1,r_{p-1}[.$

  Now if we assume that $\gamma_{k_0,r_p'}\neq0,$ then equation \ref{eq3-recall-d-KeyEq} implies that $\gamma_{i_0,r_p'}\neq 0.$ Moreover we have

  \begin{equation}\label{eq4-recall-d-KeyEq}
    \gamma_{i_0,r_p'}=-\left(\frac{s_{k_0,i_0}}{2}\right)\gamma_{k_0,r_p'}.
  \end{equation}

  Taking any $i_1\in[1,i_0[,$ we replace $i=i_1, k=k_0, r=r_{p}'$  and respectively $i=i_1, k=i_0, r=r_{p}'$ in equation \ref{Eq1-theo-Key-condition}, we obtain

  \begin{equation}\label{eq5-recall-d-KeyEq}
  \left\{\begin{array}{c}
           \left(2\gamma_{i_1,r_{p}'}+s_{k_0,i_1}\gamma_{k_0,r_{p}'}-b_{r_p',i_1}\right)
    \gamma_{k,r_{p}'}=0 \\
           \quad\\
           \left(2\gamma_{i_1,r_{p}'}+s_{i_0,i_1}\gamma_{i_0,r_{p}'}-b_{r_{p}',i_1}\right)
    \gamma_{i_0,r_{p}'}=0
         \end{array}\right.
  \end{equation}
  where $b_{ri}=\sum_{j<r}t_{rj}\gamma_{ij}.$

  But $\gamma_{k_0,r_{p}'}\neq 0$ and $\gamma_{i_0,r_{p}'}\neq 0,$ then

  \begin{equation}\label{eq6-recall-d-KeyEq}
  s_{r_p',i_1}\gamma_{k_0,r_{p}'}=s_{i_0,i_1}\gamma_{i_0,r_{p}'}
  \end{equation}
  Now by equation \ref{eq4-recall-d-KeyEq}, we conclude that

  \begin{equation}\label{eq7-recall-d-KeyEq}
    \Delta_{i_1,i_0,k_0}^{(2)}(S)=2s_{r_p',i_1}+s_{r_p',i_0}s_{i_0,i_1}=0
  \end{equation}
  for each $i_1\in[1,i_0[.$ But this contradicts our assumption that $S$ is already in the $1-$REF. Then we have that $\gamma_{k,r_p'}=0,$ for all $k\in [r_p',r_{p+1}'[.$ By a recursive argument (same as in the proof of \ref{prop-Gamma-triang-1}), we can prove that $\gamma_{k,r_p'}=0$ for all $k\geq r_p'.$ But this implies that $\det(\Gamma)=0,$ another contradiction.

  The final conclusion is that there is no an invertible matrix $\Gamma$ satisfying equation \ref{Eq1-theo-Key-condition}, and then $T\nsim S.$
\end{proof}

\begin{example}\label{ex-rp-fails}
Consider the matrices
\begin{equation*}
  S=\left[\begin{matrix}
              0 & 0 & 0& 0& 0& 0 \\
              0 & 0 & 0& 0& 0& 0 \\
              {\color{blue}1} & {\color{blue}1}& 0& 0& 0& 0 \\
              2 & 2 & {\color{red}2}& 0& 0& 0 \\
              3 & 3 & 3 & {\color{green}\boldsymbol{3}}& 0& 0 \\
              4 & 4 & 4& {\color{green}\boldsymbol{4}}& 0& 0
            \end{matrix}\right],\quad
T=\left[\begin{matrix}
              0 & 0 & 0& 0& 0& 0 \\
              0 & 0 & 0& 0& 0& 0 \\
              {\color{blue}1} & {\color{blue}1}& 0& 0& 0& 0 \\
              2 & 2 & {\color{red}2}& 0& 0& 0 \\
              3 & 3 & {\color{red}3}& 0& 0& 0 \\
              4 & 4 & 4& {\color{green}\boldsymbol{4}}& {\color{green}\boldsymbol{4}}& 0
            \end{matrix}\right].
\end{equation*}
It is clear that they are in the $1-$REF, they also have the same number of bricks, but $\mathbb{W}(S)=(3,4,5),$ while $\mathbb{W}(T)=(3,4,6).$ The walls $\mathbb{W}(S)$ and $\mathbb{W}(T)$ are different in the third component, therefore Lemma \ref{lemma-rd-fails}, implies that $S\nsim T.$
\end{example}

\begin{lemma}\label{lemma-h-fails}
Let $T,S$ be two \nSLTM, such that both are in $1-$REF, $\mathbb{W}(T)=(r_1,\dots,r_h)$ and  $\mathbb{W}(S)=(r_1',\dots,r_{g}').$ If $h\neq g,$ then $T\nsim S.$
\end{lemma}

\begin{proof}

 If $h\neq g,$ we can assume without loss of generality, that $h<g,$ that means that $T$ has fewer bricks than $S.$ We also assume that $r_j=r_j'$ for all $j\in [1,h],$ since Lemma  \ref{lemma-rd-fails} covers other cases.

  Let $p=h+1$  and suppose that  $\Gamma=[\gamma_{ij}]$ defines an isomorphism $\bgamma:\calA(T)\rightarrow \calA(S),$ then we have that $\Gamma$ has the same form, given in equation  \ref{eq03-Gamma-triangular}.
  Fix $k_0\in [r_p',r_{p+1}[$,and take $i_0$ such that $s_{k_0i_0}$ is the leader of the $k_0-$th row of $S.$ If we replace $k=k_0$ $i=i_0$ and $r=r_p'$ in equation \ref{Eq1-theo-Key-condition}, we obtain

  \begin{equation}\label{eq3-coro-recall-d-KeyEq}
    (2\gamma_{{i}r}+s_{ki}\gamma_{kr})\gamma_{kr}=0,
  \end{equation}
  (compare with equation \ref{eq3-recall-d-KeyEq}) the rest of the proof follows using exactly the same arguments that we used in the proof of Lemma \ref{lemma-rd-fails} from equation \ref{eq3-recall-d-KeyEq} onwards.
\end{proof}

\begin{example}\label{ex-h-fails}
  Consider the matrices
  \begin{equation*}
    S=\left[\begin{matrix}
              0 & 0 & 0& 0& 0& 0 \\
              0 & 0 & 0& 0& 0& 0 \\
              {\color{blue}1} & {\color{blue}1} & 0& 0 & 0& 0 \\
              {\color{blue}2} & {\color{blue}2} & 0& 0 & 0& 0 \\
              3 & 3 & {\color{red}3}& {\color{red}3} & 0& 0 \\
              4 & 4 & {\color{red}4}& {\color{red}4} & 0& 0
            \end{matrix}\right],\quad T=\left[\begin{matrix}
              0 & 0 & 0& 0& 0& 0 \\
              0 & 0 & 0& 0& 0& 0 \\
              {\color{blue}1} & {\color{blue}1}& 0& 0& 0& 0 \\
              2 & 2 & {\color{red}2}& 0& 0& 0 \\
              3 & 3 & {\color{red}3}& 0& 0& 0 \\
              4 & 4 & 4& {\color{green}\boldsymbol{4}}& {\color{green}\boldsymbol{4}}& 0
            \end{matrix}\right].
  \end{equation*}
  It is clear that $S$ and $T$ are in the $1-$REF. It is also clear that $S$ contains 2 bricks while $T$ contains $3$ bricks. Therefore Lemma \ref{lemma-h-fails} implies that $S\nsim T.$
\end{example}


\begin{theorem}\label{theo-invariants-1REF}$(\mathbb{W}-\textrm{\textbf{Test}})$
  Let $S,T$ be two \nSLTM, both in $1-$REF. If $S\sim T,$ then $\mathbb{W}(S)=\mathbb{W}(T).$
\end{theorem}

\begin{proof}(of Theorem \ref{theo-invariants-1REF})
  Lemmas \ref{lemma-r1-fails}, \ref{lemma-rd-fails} and \ref{lemma-h-fails} implies that if $\mathbb{W}(T)\neq \mathbb{W}(S),$ then $T\nsim S.$
\end{proof}

\begin{corollary}\label{coro-theo-invariants-1-REF}
Let $T,S_1,S_2$ be three \nSLTM. Suppose that $S_1,S_2$ are in $1-$REF, $T\sim S_1$ and $T\sim S_2,$ then $\mathbb{W}(S_1)=\mathbb{W}(S_2).$
\end{corollary}
\begin{proof}
Since $\sim$ is an equivalence relation, if $T\sim S_1$ and $T\sim S_2$ implies that $S_1\sim S_2.$ Now Theorem \ref{theo-invariants-1REF} implies $\mathbb{W}(S_1)=\mathbb{W}(S_2).$
\end{proof}

Corollary \ref{coro-theo-invariants-1-REF} and Proposition \ref{Prop-Agorithm-1REF} allows us to give a more general definition of \emph{Walls}:

\begin{definition}\label{def-general-wall}
Let $T$ be any \SLTM. Suppose that $S$ is a \SLTM, that is in $1-$REF and $T\sim S.$ We define the \emph{wall} of $T$ by $\mathbb{W}(T)=\mathbb{W}(S).$
\end{definition}

Corollary \ref{coro-theo-invariants-1-REF} implies that the wall of any \SLTM, is well defined. We also have:

\begin{corollary}\label{coro2-theo-invariants-1-REF}
Let $T,S$ be two \SLTM. If $T\sim S,$ then $\mathbb{W}(T)=\mathbb{W}(S).$
\end{corollary}

\begin{proof}
Is a direct consequence of Corollary \ref{coro-theo-invariants-1-REF} and definition \ref{def-general-wall}.
\end{proof}

 Note that Theorem \ref{theo-invariants-1REF} asserts that if $S\sim T,$ then their walls are equal, that means in particular that the sizes of the corresponding bricks are equal, but it does not asserts that the bricks are equal. For example if we consider the matrix $S$ of example \ref{ex-walls}, and $S'$ in equation \ref{eq-ex-walls-equal1} below, one can check that they have the same wall, that is $\mathbb{W}(S)=\mathbb{W}(S')=(3,5)$ but their corresponding bricks are not equals, despite that, we have $S'\sim S.$ In fact $S'=\calP_6(\calF_{(3,4)}(S),4)$.

 On the other hand, to have the same wall is not a guarantee that the matrices will be in the same class (with respect to $\sim$). For example if we consider the matrices $S'$ and $S''$ of equation \ref{eq-ex-walls-equal1}, then we have $\mathbb{W}(S')=\mathbb{W}(S'')=(3,5),$ but in this case we have $S'\nsim S'',$ as we shall show in Subsection \ref{ssec-invariants-2REF}.

 \begin{equation}\label{eq-ex-walls-equal1}
   S'=\left[\begin{matrix}
              0 & 0 & 0& 0& 0& 0 \\
              0 & 0 & 0& 0& 0& 0 \\
              {\color{blue}2} & {\color{blue}2} & 0& 0 & 0& 0 \\
              {\color{blue}1} & {\color{blue}1} & 0& 0 & 0& 0 \\
              3 & 3 & {\color{red}3}& {\color{red}3} & 0& 0 \\
              1 & 1 & {\color{red}1}& {\color{red}1} & 0& 0
            \end{matrix}\right],\quad S''=\left[\begin{matrix}
              0 & 0 & 0& 0& 0& 0 \\
              0 & 0 & 0& 0& 0& 0 \\
              {\color{blue}2} & {\color{blue}2} & 0& 0 & 0& 0 \\
              {\color{blue}1} & {\color{blue}1} & 0& 0 & 0& 0 \\
              3 & 3 & {\color{red}3}& {\color{red}3} & 0& 0 \\
              1 & 1 & {\color{red}0}& {\color{red}1} & 0& 0
            \end{matrix}\right]
 \end{equation}

 Now we naturally can ask, what if the walls and the corresponding bricks are equal, are the matrices necessarily related? Well, the answer is not, for example consider the following three matrices:

 \begin{equation}\label{eq-ex-walls-equal2}
   T_1=\left[\begin{matrix}
               0 & 0 & 0 & 0 \\
               0 & 0 & 0 & 0 \\
               {\color{blue}1} & {\color{blue}1} & 0 & 0 \\
               0 & 0 & {\color{red}1} & 0
             \end{matrix}\right],\quad
 T_2=\left[\begin{matrix}
               0 & 0 & 0 & 0 \\
               0 & 0 & 0 & 0 \\
               {\color{blue}1} & {\color{blue}1} & 0 & 0 \\
               -1 & -1 & {\color{red}1} & 0
             \end{matrix}\right],\quad
 T_3=\left[\begin{matrix}
               0 & 0 & 0 & 0 \\
               0 & 0 & 0 & 0 \\
               {\color{blue}1} & {\color{blue}1} & 0 & 0 \\
               0 & -\frac{1}{2} & {\color{red}1} & 0
             \end{matrix}\right],
 \end{equation}

 In this case we have $\mathbb{W}(T_1)=\mathbb{W}(T_2)=\mathbb{W}(T_3)=(3,4),$ moreover we have $B_1(T_1)=B_1(T_2)=B_1(T_3)=[{\color{blue}1},{\color{blue}1}]$ and $B_2(T_1)=B_2(T_2)=B_2(T_3)=[{\color{red}1}]$ but they are different in their dark sides. The reader can check, by using Theorem \ref{theo-Key-condition}, that $T_1\sim T_2,$ but $T_2\nsim T_3.$ At this point we do not have any criteria to control, in general, the information coming from the dark side of a matrix, we expect to address this analysis in future works.

 The correct way to use Theorem \ref{theo-invariants-1REF} is this: If we have two \nSLTM, with different walls, then we can conclude immediately that they are not related under relation $\sim.$ For example if $S,T,U$ are the matrices of example \ref{ex-walls}, then since their walls are pairwise different, then we can conclude that their classes $\overline{S},\overline{T},\overline{U}$ are pairwise different. As we can see in that example, Theorem \ref{theo-invariants-1REF} is very visual, it concerns more on the shapes and distribution of the bricks than the specific data contained in each matrix:

 \begin{equation*}
    S=\left[\begin{matrix}
              \circ & \circ & \circ& \circ& \circ& \circ \\
              \circ & \circ & \circ& \circ& \circ& \circ \\
              {\color{blue}\bullet} & {\color{blue}\bullet} & \circ& \circ& \circ& \circ \\
              {\color{blue}\bullet} & {\color{blue}\bullet} & \circ& \circ& \circ& \circ \\
              \bullet & \bullet & {\color{red}\bullet}& {\color{red}\bullet} & \circ& \circ \\
              \bullet & \bullet & {\color{red}\bullet}& {\color{red}\bullet} & \circ& \circ
            \end{matrix}\right],\quad T=\left[\begin{matrix}
              \circ & \circ & \circ& \circ& \circ& \circ \\
              \circ & \circ & \circ& \circ& \circ& \circ \\
              {\color{blue}\bullet} & {\color{blue}\bullet} & \circ& \circ& \circ& \circ \\
              \bullet & \bullet & {\color{red}\bullet}& \circ& \circ& \circ \\
              \bullet & \bullet & {\color{red}\bullet}& \circ& \circ& \circ \\
              \bullet & \bullet & \bullet& {\color{green}\bullet}& {\color{green}\bullet}& \circ
            \end{matrix}\right],\quad U=\left[\begin{matrix}
              {\circ} & \circ & \circ& \circ& \circ& \circ \\
              \circ & \circ & \circ& \circ& \circ& \circ \\
              {\color{blue}\bullet} & {\color{blue}\bullet} & \circ& \circ& \circ& \circ\\
              \bullet & \bullet & {\color{red}\bullet}& \circ& \circ& \circ \\
              \bullet & \bullet & \bullet& {\color{orange}\bullet}&  \circ& \circ \\
              \bullet & \bullet & \bullet& \bullet& {\color{green}\bullet}&  \circ
            \end{matrix}\right],
  \end{equation*}

  We can conclude the same if we consider any three matrices $S,T,U$ with the same corresponding walls as the matrices $S,T,U$ of example \ref{ex-walls}.


\subsection{Measure Sequences}\label{ssec-invariants-2REF}

Recall that if $M$ is a matrix $\mu(M)$ denotes the number of nonzero entries of $M.$ in particular $\mu(M^{[r]})$ denotes the number of nonzero entries of the $r-$th row of $M$ (see Definition \ref{def-measure-mat}).

 Let $T,S$ be two \SLTM, such that $T\sim S$ and $S$ is in $1-$REF.  If $\mathbb{W}(T)=(r_1,\dots,r_h),$ then the bricks of $T$ are given by:
\begin{equation*}
  B_1(T)=S_{[1,r_1[}^{[r_1,r_2[},\quad B_{2}(T)=S_{[r_1,r_2[}^{[r_2,r_3[},\dots,
  B_{h}(T)=S_{[r_{h-1},r_h[}^{[r_h,n]}.
\end{equation*}

\begin{definition}\label{def-second-red-form}
 Let $T\nsim 0_n$ be a \nSLTM, and let $\mathbb{W}(T)=(r_1,\dots,r_h)$ its wall. We say that $T$ is in the \emph{Second Reduced Echelon Form} ($2-$REF) if it satisfies the following conditions:
 \begin{enumerate}
   \item $T$ contains exactly $r_1-1$ zero rows and they are above of the nonzero rows.
   \item If $(r,c_r)$ is the leader position of the $r-$th row of $T,$ then $c_r>1$ and there is a $i_0\in [1,c_r[$ such that $\Delta^{(2)}_{i_0,c,r}\neq 0.$ Moreover for each $j\in [1,h]$ and $r\in[r_j,r_{j+1}[$ we have that, $c_r\in [r_{j-1},r_j[.$
   \item For each $j\in [1,h]$ and $r\in[r_j,r_{j+1}-1[$ we have that
   \begin{equation*}
     \mu\left(T_{[r_j,r_{j+1}[}^{[r]}\right)\leq \mu\left(T_{[r_j,r_{j+1}[}^{[r+1]}\right).
   \end{equation*}
   We also say that the zero matrix is in the $2-$REF.
 \end{enumerate}
\end{definition}

Definition \ref{def-second-red-form}, basically say that $T$ is in the $2-$REF, then $T$ is almost in the $1-$REF. The first and the second conditions in Definition \ref{def-second-red-form} are more specific versions of the first and the second condition of Definition \ref{def-first-red-form}. In particular, the second condition in Definition \ref{def-second-red-form} say that leader position $(r,c_r)$ for $r\in [r_j,r_{j+1}[$ is located in the corresponding brick $B_j$ of $T.$ The main difference between $1-$REF and $2-$REF is the third condition of Definition \ref{def-first-red-form} and Definition \ref{def-second-red-form} respectively. In $1-$REF we ask that leaders have to be located in a \emph{increasingly} way, with respect to the column $c_r,$ but in $2-$REF we ask that in each brick, the rows are located in a \emph{increasingly} way with respect to the measure $\mu.$ Note that the third condition in Definition \ref{def-second-red-form}, only consider the measure of the rows of the brick, it does not consider the measure of the complete rows of $T.$

\begin{example}\label{ex1-ref2}
  Consider the matrices
  \begin{equation*}
  S=\left[\begin{matrix}
              0 & 0 & 0& 0& 0& 0 \\
              0 & 0 & 0& 0& 0& 0 \\
              {\color{blue}2} & {\color{blue}2} & 0& 0 & 0& 0 \\
              {\color{blue}1} & {\color{blue}1} & 0& 0 & 0& 0 \\
              0 & 0 & {\color{red}3}& {\color{red}3} & 0& 0 \\
              1 & 1 & {\color{red}0}& {\color{red}1} & 0& 0
            \end{matrix}\right],\quad
            T=\left[\begin{matrix}
              0 & 0 & 0& 0& 0& 0 \\
              0 & 0 & 0& 0& 0& 0 \\
              {\color{blue}2} & {\color{blue}2} & 0& 0 & 0& 0 \\
              {\color{blue}1} & {\color{blue}1} & 0& 0 & 0& 0 \\
              1 & 1 & {\color{red}0}& {\color{red}1} & 0& 0 \\
              0 & 0 & {\color{red}3}& {\color{red}3} & 0& 0
            \end{matrix}\right].
 \end{equation*}
 Note that both $S$ and $T$ are in the $1-$REF, but only $T$ is in the $2-$REF. Even when $\mu(T^{[5]})>\mu(T^{[6]}),$ we only have to check that $\mu(T_{[3,4]}^{[5]})\leq \mu(T_{[3,4]}^{[6]}).$ Note that  $\mathbb{W}(S)=\mathbb{W}(T)=(3,5)$ and $S\sim T,$ in fact $T=\calF_{(5,6)}(S).$
\end{example}

\begin{example}\label{ex2-ref2}
  Consider the matrices
  \begin{equation*}
    S=\left[\begin{matrix}
               0 & 0 & 0 & 0 & 0 & 0 \\
               0 & 0 & 0 & 0 & 0 & 0 \\
               0 & 0 & 0 & 0 & 0 & 0 \\
              0 & 0 & 0 & 0 & 0 & 0 \\
              {\color{blue}1} & {\color{blue}1} & {\color{blue}1} & {\color{blue}0} & 0 & 0 \\
              {\color{blue}0} & {\color{blue}1} & {\color{blue}0} & {\color{blue}1} & 0 & 0
            \end{matrix}\right],\quad
            T=\left[\begin{matrix}
               0 & 0 & 0 & 0 & 0 & 0 \\
               0 & 0 & 0 & 0 & 0 & 0 \\
               0 & 0 & 0 & 0 & 0 & 0 \\
              0 & 0 & 0 & 0 & 0 & 0 \\
              {\color{blue}0} & {\color{blue}1} & {\color{blue}0} & {\color{blue}1} & 0 & 0\\
              {\color{blue}1} & {\color{blue}1} & {\color{blue}1} & {\color{blue}0} & 0 & 0
            \end{matrix}\right]
  \end{equation*}
  In this case $S$ is in $1-$REF but not in $2-$REF, while $T$ is in $2-$REF but not in $1-$REF. Note that $\mathbb{W}(S)=\mathbb{W}(T)=(5)$ and $S\sim T,$ in fact $T=\calF_{(5,6)}(S).$
\end{example}

\begin{proposition}\label{Prop-Agorithm-2REF}
  Let $T$ be a \nSLTM, such that $T\nsim 0_n.$ Then there is a sequence of ETO's that transform $T$ into a matrix $V$ that is in the second reduced echelon form.
\end{proposition}

\begin{proof}
By proposition \ref{Prop-Agorithm-1REF}, we can transform $T$ into a matrix $R$ that is in $1-$REF. By Definition \ref{def-second-red-form}, we only have to reorganize the rows of $R$ in order to satisfy condition (3) and we can make this by applying a finite sequence of $\calF-$operations, obtaining a matrix $V$ in $2-$REF such that $T\sim V.$
\end{proof}

\begin{remark}\label{remark1-Algorithm2}
  Note that each time we apply one of the $\calF-$operations described in the proof of Proposition \ref{Prop-Agorithm-2REF}, we only affect the rows of the matrix that correspond to a brick $B_j(R),$ and the columns of the matrix that correspond to the next brick $B_{j+1}(R),$ In particular we do not affect the measure of any row of the brick $B_{j+1}(R).$ Moreover we do not change any other entry of $R.$
\end{remark}

If $\mathbb{W}(T)=(r_1,\dots,r_{h})$ and $V$ is the matrix obtained in the proof of Proposition \ref{Prop-Agorithm-2REF}, then for each $j\in[1,h],$ we define the brick $\overline{B}_{j}(V)$ as follows:

\begin{equation}\label{eq-def-new-bricks}
  \overline{B}_{j}(V)=V_{[r_{j-1},r_j[}^{[r_j,r_{j+1}[}.
\end{equation}

Then the matrix $V,$ looks like this:

\begin{equation}\label{eq-form-of-2REF}
  V=\left[\begin{matrix}
                0 & 0 & \cdots & 0&0 \\
                \overline{B}_1(V) & 0 &\cdots & 0&0 \\
                \bullet & \overline{B}_{2}(V)&\cdots & 0&0\\
                \vdots & \vdots &\ddots&\vdots&\vdots\\
                \bullet & \bullet & \bullet& \overline{B}_{h}(V)&0
              \end{matrix}\right]
\end{equation}

Note that, if $R\sim T,$ is a matrix in $1-$REF obtained in Proposition \ref{Prop-Agorithm-1REF}, and   $V\sim T$ is the matrix obtained in the proof of Proposition \ref{Prop-Agorithm-2REF}, then each brick  $\overline{B}_{j}(V)$ of $V,$ covers the same area of the corresponding brick $B_j(R)$ of $R,$ but where the position of the entries was changed.

\begin{example}\label{ex3-ref2}
  Consider the matrix
  \begin{equation*}
    S=\left[\begin{matrix}
              0 & 0 & 0 & 0 & 0 & 0 & 0 &0\\
              0 & 0 & 0 & 0 & 0 & 0 & 0 &0\\
              0 & 0 & 0 & 0 & 0 & 0 & 0 &0\\
              1 & 1 & 1 & 0 & 0 & 0 & 0 &0\\
              1 & 0 & 1 & 0 & 0 & 0 & 0 &0\\
              0 & 0 & 0 & 2 & 2 & 0 & 0 &0\\
              1 & 1 & 1 & 0 & 3 & 0 & 0 &0\\
              1 & 2 & 3 & 4 & 5 & 7 & 0 &0
            \end{matrix}\right]
  \end{equation*}
  Note that $S$ is already in the $1-$REF. Also note that $\mathbb{W}(S)=(4,6,8),$ then

  \begin{equation*}
  {B}_1(S)=S_{[1,4[}^{[4,6[}=\left[\begin{matrix}
                            1 & 1 & 1 \\
                            1 & 0 & 1
                          \end{matrix}\right],\quad
  {B}_2(S)=S_{[4,6[}^{[6,8[}=\left[\begin{matrix}
                            2 & 2 \\
                             0 & 3
                          \end{matrix}\right],\quad
  {B}_3(S)S_{[6,8[}^{[8]}=\left[\begin{matrix}
                             6 & 7
                          \end{matrix}\right].
\end{equation*}

  Since $\mu\left(S_{[1,4[}^{[4]}\right)>\mu\left(S_{[1,4[}^{[5]}\right),$ we have to apply $\calF_{(4,5)}$ to $S:$

  \begin{equation*}
    S=\left[\begin{matrix}
              0 & 0 & 0 & 0 & 0 & 0 & 0 &0\\
              0 & 0 & 0 & 0 & 0 & 0 & 0 &0\\
              0 & 0 & 0 & 0 & 0 & 0 & 0 &0\\
              {\color{blue}1} & {\color{blue}1} & {\color{blue}1} & 0 & 0 & 0 & 0 &0\\
              {\color{blue}1} & {\color{blue}0} & {\color{blue}1} & 0 & 0 & 0 & 0 &0\\
              0 & 0 & 0 & {\color{red}2} & {\color{red}2} & 0 & 0 &0\\
              1 & 1 & 1 & {\color{red}0} & {\color{red}3} & 0 & 0 &0\\
              1 & 2 & 3 & 4 & 5 & {\color{green}\boldsymbol{6}} & {\color{green}\boldsymbol{7}} &0
            \end{matrix}\right]\rightarrow
  S_1=\calF_{(4,5)}(S)=
              \left[\begin{matrix}
              0 & 0 & 0 & 0 & 0 & 0 & 0 &0\\
              0 & 0 & 0 & 0 & 0 & 0 & 0 &0\\
              0 & 0 & 0 & 0 & 0 & 0 & 0 &0\\
              {\color{blue}1} & {\color{blue}0} & {\color{blue}1} & 0 & 0 & 0 & 0 &0\\
              {\color{blue}1} & {\color{blue}1} & {\color{blue}1} & 0 & 0 & 0 & 0 &0\\
              0 & 0 & 0 & {\color{red}2} & {\color{red}2} & 0 & 0 &0\\
              1 & 1 & 1 & {\color{red}3} & {\color{red}0} & 0 & 0 &0\\
              1 & 2 & 3 & 5 & 4 & {\color{green}\boldsymbol{6}} & {\color{green}\boldsymbol{7}} &0
            \end{matrix}\right]
  \end{equation*}
  Now note that $\mu\left(S_{[4,6[}^{[6]}\right)>\mu\left(S_{[4,6[}^{[7]}\right)$ we have to apply $\calF_{(6,7)}$ to $S_1:$

  \begin{equation*}
    S_1=
              \left[\begin{matrix}
              0 & 0 & 0 & 0 & 0 & 0 & 0 &0\\
              0 & 0 & 0 & 0 & 0 & 0 & 0 &0\\
              0 & 0 & 0 & 0 & 0 & 0 & 0 &0\\
              {\color{blue}1} & {\color{blue}0} & {\color{blue}1} & 0 & 0 & 0 & 0 &0\\
              {\color{blue}1} & {\color{blue}1} & {\color{blue}1} & 0 & 0 & 0 & 0 &0\\
              0 & 0 & 0 & {\color{red}2} & {\color{red}2} & 0 & 0 &0\\
              1 & 1 & 1 & {\color{red}3} & {\color{red}0} & 0 & 0 &0\\
              1 & 2 & 3 & 5 & 4 & {\color{green}\boldsymbol{6}} & {\color{green}\boldsymbol{7}} &0
            \end{matrix}\right]\rightarrow
  S_2=\calF_{(6,7)}(S_1)=
              \left[\begin{matrix}
              0 & 0 & 0 & 0 & 0 & 0 & 0 &0\\
              0 & 0 & 0 & 0 & 0 & 0 & 0 &0\\
              0 & 0 & 0 & 0 & 0 & 0 & 0 &0\\
              {\color{blue}1} & {\color{blue}0} & {\color{blue}1} & 0 & 0 & 0 & 0 &0\\
              {\color{blue}1} & {\color{blue}1} & {\color{blue}1} & 0 & 0 & 0 & 0 &0\\
              1 & 1 & 1 & {\color{red}3} & {\color{red}0} & 0 & 0 &0\\
              0 & 0 & 0 & {\color{red}2} & {\color{red}2} & 0 & 0 &0\\
              1 & 2 & 3 & 5 & 4 & {\color{green}\boldsymbol{7}} & {\color{green}\boldsymbol{6}} &0
            \end{matrix}\right]
  \end{equation*}
One can check that the matrix $S_2$ is in the $2-$REF, therefore we have $V=S_2.$ The new bricks are given as follows:
\begin{equation*}
  \overline{B}_1(V)=V_{[1,2]}^{[4,5]}=\left[\begin{matrix}
                            1 & 0 & 1 \\
                            1 & 1 & 1
                          \end{matrix}\right],\quad
  \overline{B}_2(V)=V_{[4,5]}^{[6,7]}=\left[\begin{matrix}
                            3 & 0 \\
                             2 & 2
                          \end{matrix}\right],\quad
  \overline{B}_3(V)=V_{[6,7]}^{[8]}=\left[\begin{matrix}
                             7 & 6
                          \end{matrix}\right].
\end{equation*}
\end{example}

\begin{definition}\label{def-Blocks-2REF}
  Let $T$ be a \SLTM, that is in $2-$REF, let $\mathbb{W}(T)=(r_1,\dots,r_h).$ For each $j\in[1,h]$ we define recursively numbers $r_{jk}$ and $\mu_{jk}$ and a matrix $M_j(T),$ as follows:
  \begin{enumerate}
    \item Let $r_{j1}:=r_j,$ and $\mu_{j1}=\mu\left(T_{[r_{j-1},r_j[}^{[r_{j1}]}\right).$
    \item Assuming that the numbers $r_{jk}$ and $\mu_{jk}$ are already defined. Define the set
         \begin{equation*}
           A_{jk}=\left\{r\in]r_{jk},r_{j+1}[\quad : \quad \mu\left(T_{[r_{j-1},r_j[}^{[r]}\right)>\mu_{jk}\right\}.
         \end{equation*}
         \begin{itemize}
           \item If $A_{jk}\neq \emptyset,$ then define $r_{j,k+1}:=\min(A_{jk})$ and $\mu_{j,k+1}:=\mu\left(T_{[r_{j-1},r_j[}^{[r_{j,k+1}]}\right)$
           \item If $A_{jk}=\emptyset,$ then define the matrix $M_j(T)$ as follows:
               \begin{equation*}
                 M_j(T):=\left[\begin{matrix}
                                 r_{j1} & r_{j2} & \cdots & r_{jk} \\
                                 \mu_{j1} & \mu_{j2} & \cdots & \mu_{jk}
                               \end{matrix}\right].
               \end{equation*}
         \end{itemize}
          We call the sequence $\mathbb{M}(T):=(M_1(T),\dots,M_h(T))$ the \emph{measure sequence} of $T.$
  \end{enumerate}
\end{definition}

\begin{example}
  Consider the matrix $T$ of Example \ref{ex1-ref2}. As we mentioned before $T$ is already in both the $1-$REF and the $2-$REF. We also know that $\mathbb{W}(T)=(3,5),$ then the bricks of $T$ are:
  \begin{equation*}
    B_1(T)=T_{[1,3[}^{[3,5[}=\left[\begin{matrix}
                                     2 & 2 \\
                                     1 & 1
                                   \end{matrix}\right],\quad
   B_2(T)=T_{[3,5[}^{[5,6]}=\left[\begin{matrix}
                                     0 & 1 \\
                                     3 & 3
                                   \end{matrix}\right]
  \end{equation*}
  For $j=1$ we can see that both rows of $B_1(T)$ have the same measure, $\mu_{1,1}=2.$ On the other hand for $j=2$ the measure of the first row of $B_2(T)$ is $\mu_{2,1}=1$ while the measure of the second row is $\mu_{2,2}=2.$ Therefore we have:
  \begin{equation*}
    M_1(T)=\left[\begin{matrix}
                   3 \\
                   2
                 \end{matrix}\right],\quad
    M_2(T)=\left[\begin{matrix}
                   5 & 6 \\
                   1 & 2
                 \end{matrix}\right],
  \end{equation*}
  then the measure sequence of $T$ is:
  \begin{equation*}
    \mathbb{M}(T)=\left(\left[\begin{matrix}
                   3 \\
                   2
                 \end{matrix}\right],
    \left[\begin{matrix}
                   5 & 6 \\
                   1 & 2
                 \end{matrix}\right]\right)
  \end{equation*}
\end{example}

\begin{example}
  Consider the matrix $T$ of example \ref{ex2-ref2}. As we mentioned before, $T$ is already in the $2-$REF, but in this case $T$ is not in the $1-$REF. We also know that $\mathbb{W}(T)=(5).$ To obtain the table of measures of $T,$ we have to watch the block:
  \begin{equation*}
    T_{[1,5[}^{[5,6]}=\left[\begin{matrix}
                              0 & 1 & 0 & 1 \\
                              1 & 1 & 1 & 0
                            \end{matrix}\right]
  \end{equation*}
  In this case we have $r_{1,1}=r_1=5,$ and $\mu_{1,1}=2.$ Since $\mu\left(T_{[1,5[}^{[6]}\right)=3>2$ we have $r_{1,2}=6$ and $\mu_{1,2}=3.$ Therefore we have:
  \begin{equation*}
    \mathbb{M}(T)=\left(\left[\begin{matrix}
                   5 & 6 \\
                   2 & 3
                 \end{matrix}\right]\right).
  \end{equation*}
\end{example}

\begin{example}
  Consider the matrix $V=S_2$ of Example \ref{ex3-ref2}. Then $V$ is in $2-$REF and $\mathbb{W}(V)=(3,6,8).$ The reader can check that
  \begin{equation*}
    \mathbb{M}(V)=\left(\left[\begin{matrix}
                   4 & 5 \\
                   2 & 3
                 \end{matrix}\right],\left[\begin{matrix}
                   6 & 7 \\
                   1 & 2
                 \end{matrix}\right],\left[\begin{matrix}
                   8  \\
                   2
                 \end{matrix}\right]\right)
  \end{equation*}
\end{example}

\begin{lemma}\label{lemma-from-1REF-to-2REF}
  Let $S,T$ be Let $T,S$ be two non-zero \nSLTM, such that:
   \begin{itemize}
     \item $T$ and $S$ are in the $2-$REF,
     \item $\mathbb{W}(T)=\mathbb{W}(S)=(r_1,\dots,r_h).$
     \end{itemize}
   If we assume that there is an isomorphism $\bgamma:\calA(T)\rightarrow \calA(S),$ then the associated matrix $\Gamma=[\gamma_{ij}]$ has the form given in equation \ref{eq0-Gamma-triangular} and there is a bijection $\rrn$ satisfying conditions (2.a) and (2.b) of Proposition \ref{prop-Gamma-triang-1}.
\end{lemma}

\begin{proof}
  Since we are assuming that $T$ and $S$ are already in the $2-$REF, we have that
  \begin{equation}\label{eq2-two-SLTM-by-blocks}
      T=\left[\begin{matrix}
                0 & 0 & \cdots & 0&0 \\
                \overline{B}_1(T) & 0 &\cdots & 0&0 \\
                \bullet & \overline{B}_{2}(T)&\cdots & 0&0\\
                \vdots & \vdots &\ddots&\vdots&\vdots\\
                \bullet & \bullet & \bullet& \overline{B}_{h}(T)&0
              \end{matrix}\right],\quad
      S=\left[\begin{matrix}
                0 & 0 & \cdots & 0&0 \\
                \overline{B}_1(S) & 0 &\cdots & 0&0 \\
                \bullet & \overline{B}_{2}(S)&\cdots & 0&0\\
                \vdots & \vdots &\ddots&\vdots&\vdots\\
                \bullet & \bullet & \bullet& \overline{B}_{h}(S)&0
              \end{matrix}\right]
    \end{equation}
    (compare with equation \ref{eq1-two-SLTM-by-blocks}, in the proof of Proposition \ref{prop-Gamma-triang-1}). If we suppose that there is an isomorphism $\bgamma:\calA(T)\rightarrow \calA(S),$ then, even if $S$ or $T$ are not in their $1-$REF, one can pass from $1-$REF to $2-$REF and vice-versa, only by applying $\calF-$operations that do not change the form of $\Gamma$ given in equation \ref{eq0-Gamma-triangular} (see the proof of Proposition \ref{Prop-Agorithm-2REF}).  This also implies that, after we have applied the sequence of $\calF-$operations needed to pass from $1-$REF to $2-$REF, we obtain a bijection $\rrn,$ that is not exactly the same bijection obtained in Proposition \ref{prop-Gamma-triang-1}, but it satisfies the desired conditions as well.
\end{proof}

\begin{proposition}\label{prop-Gamma-triang-2}
  Let $T,S$ be two nonzero \nSLTM, such that:
   \begin{itemize}
     \item $T$ and $S$ are in the $2-$REF,
     \item $\mathbb{W}(T)=\mathbb{W}(S)=(r_1,\dots,r_h).$
     \item $\mathbb{M}(T)=\mathbb{M}(S)=(M_1,\dots,M_h).$
   \end{itemize}

     If $\bgamma:\calA(T)\rightarrow\calA(S)$ is an isomorphism, then the associated matrix $\Gamma=[\gamma_{ij}]$ satisfies the following relations:
  \begin{enumerate}
    \item If $r\in [r_j, r_{j+1}[$ then $\gamma_{kr}=0$ for all $k\geq r_{j+1}.$
    \item There is a bijection $\rrn:[1,n]\rightarrow [1,n]$ such that:
    \begin{enumerate}
      \item For each $j\in[1,h]$ we have $\rrn([r_j,r_{j+1}[)= [r_j,r_{j+1}[$
      \item For each $i,k\in [r_j,r_{j+1}[$ we have that
    \begin{equation*}
       \gamma_{ik}\neq 0 \quad \textrm{if and only if}\quad k=\rrn(i).
    \end{equation*}
    \item For each $j\in[1,h],$ if
    \begin{equation*}
      M_j=\left[\begin{matrix}
                     r_{j1} & \cdots & r_{j,l(j)} \\
                     \mu_{j1} & \cdots  & \mu_{j,l(j)}
                   \end{matrix}\right]
    \end{equation*}
    for some $l=l(j)>1.$  Then for each $i\in[1,l(j)[$ we have that $\rrn([r_{j,i},r_{j,i+1}[)= [r_{j,i},r_{j,i+1}[.$
    \end{enumerate}
  \end{enumerate}
\end{proposition}

Note that Proposition \ref{prop-Gamma-triang-2} is a more specific version of Proposition \ref{prop-Gamma-triang-1}. Again we have that the associated matrix $\Gamma$ have a \emph{triangular by blocks} form, just like in equation \ref{eq0-Gamma-triangular}, but this time we are providing more specific information on each diagonal block $\Gamma_{[r_{j},r_{j+1}[}^{[r_{j},r_{j+1}[}$ for $j\in[1,h].$ Condition (2.c) means that the block $\Gamma_{[r_{j},r_{j+1}[}^{[r_{j},r_{j+1}[}$ has the following \emph{diagonal by blocks} form:

\begin{equation}\label{eq1-Gamma-trangular-blocks}
  \Gamma_{[r_{j},r_{j+1}[}^{[r_{j},r_{j+1}[}
  =\left[\begin{matrix}
           \Gamma_{[r_{j,1},r_{j,2}[}^{[r_{j,1},r_{j,2}[} & 0 & 0 & \cdots & 0 \\
           \quad & \quad & \quad & \quad & \quad \\
           0 & \Gamma_{[r_{j,2},r_{j,3}[}^{[r_{j,2},r_{j,3}[} & 0 & \cdots & 0 \\
           \quad & \quad & \quad & \quad & \quad \\
           \vdots & \vdots & \vdots & \ddots & \vdots\\
           \quad & \quad & \quad & \quad & \quad \\
           0 & 0 & 0 & \cdots & \Gamma_{[r_{j,l},r_{j+1}[}^{[r_{j,l},r_{j+1}[}
         \end{matrix}\right]
\end{equation}
where $\mu\left(\Gamma_{[r]}^{[r_{j,i},r_{j,i+1}[}\right)=1$ for each $r\in [r_{j,i},r_{j,i+1}[$ and each $i\in[1,l(1)].$


\begin{proof}
Assertions (1), (2.a) and (2.b) follow directly from Lemma \ref{lemma-from-1REF-to-2REF}, therefore we only have to prove assertion (2.c). In particular we are assuming that there is a bijection $\rrn:[1,n]\rightarrow [1,n]$ satisfying conditions (2.a) and (2.b).
    \begin{itemize}
      \item If
      \begin{equation*}
        M_1=\left[\begin{matrix}
                    r_{1,1} & \cdots & r_{1,l(1)} \\
                    \mu_{1,1}&  \cdots & \mu_{1,l(1)}
                  \end{matrix}\right]
      \end{equation*}
      with $l(1)>1.$ For each $r\in[r_{1},r_{2}[$ define the sets:
      \begin{equation*}
        J_r^{(1)}:=\{j\in[1,r_1[:t_{rj}\neq 0\},\quad I_r^{(1)}:=\{j\in[1,r_1[:s_{rj}\neq 0\}
      \end{equation*}
      By assumption  $|J_r^{(1)}|=|I_r^{(1)}|\geq 2$ (the inequality comes from Lemma \ref{lemma-r1-geq-3}). Moreover we have  $|J_r^{(1)}|=|I_r^{(1)}|=\mu_{1,p},$ whenever $r\in[r_{1,p},r_{1,p+1}[$ (for some $p\in[1,l(1)]$).
      \begin{itemize}
      \item If we replace $i,k\in [1,r_1[$ and $r\in[r_1,r_2[$ in equation \ref{Eq1-theo-Key-condition}, we obtain:
      \begin{equation}\label{eq2-recalling-Key}
        2\gamma_{ir}\gamma_{kr}=
        \sum_{j\in J_r^{(1)}}{t_{rj}\left(\gamma_{kj}\gamma_{ir}+
        \gamma_{ij}\gamma_{kr}\right)}=a_{r,i}\gamma_{kr}+a_{r,k}\gamma_{i,r}.
      \end{equation}
      where (as in the proof of Proposition \ref{prop-Gamma-triang-1}) $a_{r,i}:=t_{r,\rrn(i)}\gamma_{i,\rrn(i)},$ for each $i\in [1,r_1[.$
      \item If we replace $i\in [1,r_1[$ and $k,r\in[r_1,r_2[$ in equation \ref{Eq1-theo-Key-condition}, we obtain:
      \begin{equation*}
        2\gamma_{ir}\gamma_{kr}+s_{ki}\gamma_{kr}^2=
        \sum_{j\in J_r^{(1)}}{t_{rj}\left(\gamma_{kj}\gamma_{kr}s_{ki}+\gamma_{kj}\gamma_{ir}+
        \gamma_{ij}\gamma_{kr}\right)}=a_{r,i}\gamma_{kr}
      \end{equation*}
      (since, $\gamma_{kj}=0$ for each $j\in [1,r_1[$). Therefore
      \begin{equation}\label{eq3-recalling-Key}
        \left(2\gamma_{ir}+s_{ki}\gamma_{kr}-a_{r,i}\right)\gamma_{kr}=0.
      \end{equation}
      \end{itemize}
      If we assume that $\gamma_{kr}\neq 0$ for $k\in[r_1,r_2[$ and $r\in[r_{1,p},r_{1,p+1}[,$ then we can pick two different elements $i,j\in I_k^{(1)}\subset[1,r_1[$ (again by lemma \ref{lemma-r1-geq-3}) and replace them in equation \ref{eq3-recalling-Key} obtaining:
      \begin{equation}\label{eq4-recalling-Key}
        \left\{\begin{array}{cc}
                 2\gamma_{ir}= & a_{ri}-s_{ki}\gamma_{kr}  \\
                 \quad & \quad\\
                 2\gamma_{jr}= & a_{rj}-s_{kj}\gamma_{kr}
               \end{array}\right.
      \end{equation}
      By equation \ref{eq2-recalling-Key} we have:
      \begin{equation}\label{eq5-recalling-Key}
        (2\gamma_{ir})(2\gamma_{jr})=a_{ri}(2\gamma_{jr})+a_{rj}(2\gamma_{ir}).
      \end{equation}
      Combining equations \ref{eq4-recalling-Key} and \ref{eq5-recalling-Key}, we obtain:
      \begin{equation*}
        (a_{ri}-s_{ki}\gamma_{kr})(a_{rj}-s_{kj}\gamma_{kr})
        =a_{ri}(a_{rj}-s_{kj}\gamma_{kr})+a_{rj}(a_{ri}-s_{ki}\gamma_{kr}),
      \end{equation*}
      from where:
      \begin{equation}\label{eq6-recalling-Key}
        s_{ki}s_{kj}\gamma_{kr}^{2}=a_{ri}a_{rj}.
      \end{equation}
      By our assumptions, the left side of equation \ref{eq6-recalling-Key} is different from zero. Therefore $a_{ri}\neq 0$ whenever $i\in I_k^{(1)},$ but this means that $\rrn(i)\in J_r^{(1)}.$

      Conversely, if we assume that there is a $i\in [1,r_1[$ such that $\rrn(i)\in J_r^{(1)}$ but $i\notin I_k^{(1)},$ then by equation \ref{eq3-recalling-Key}, we have $2\gamma_{ir}=a_{r,i}\neq 0.$ If we replace this in equation \ref{eq2-recalling-Key}, we obtain:
      \begin{equation*}
        0=a_{r,j}\gamma_{ir}
      \end{equation*}
      for any $j\in[1,r_1[$ with $j\neq i.$ In particular for $j\in I_{k}^{(1)},$ what contradicts the fact that $\rrn(j)\in J_r^{(1)}.$

      The last analysis implies that $i\in I_{k}^{(1)}$ if and only if $\rrn(i)\in J_{r}^{(1)},$ and since $\rrn$ is bijective, we conclude that $\rrn(I_k^{(1)})=J_{r}^{(1)},$ in particular $|I_{k}^{(1)}|=|J_{r}^{(1)}|$ and therefore, having $\gamma_{kr}\neq 0$ implies that $k,r$ both belongs to the same subinterval $[r_{1,p},r_{1,p+1}[,$ for certain $p\in[1,k],$ and that means $\rrn([r_{1,p},r_{1,p+1}[)= [r_{1,p},r_{1,p+1}[.$
      We have proved that the block $\Gamma_{[r_1,r_2[}^{[r_{1},r_{2}[}$ has the form given in equation \ref{eq1-Gamma-trangular-blocks}, and using item (2.b) we conclude that $\mu\left(\Gamma_{[r]}^{[r_{1,p},r_{1,p+1}[}\right)=1$ for each $r\in [r_{1,p-1},r_{1,p}[$ and $p\in[1,l(1)].$
      \item For $q\in]1,h],$ the proof is almost analogous to the case $q=1:$

      Let
      \begin{equation*}
        M_q=\left[\begin{matrix}
                    r_{q1} & \cdots & r_{ql(q)} \\
                    \mu_{q1} & \cdots & \mu_{ql(q)}
                  \end{matrix}\right]
      \end{equation*}
      with $l(q)>1.$ For each $r\in[r_q,r_{q+1}[$ we define the sets:
          \begin{equation*}
            J_r^{(q)}=\{j\in[r_{q-1},r_q[:t_{rj}\neq 0\},\quad I_r^{(q)}=\{j\in[r_{q-1},r_q[:s_{rj}\neq 0\}
          \end{equation*}

    Then, by the hypothesis, we have $|J_r^{(q)}|=|I_r^{(q)}|.$ Moreover $|J_r^{(q)}|=|I_r^{(q)}|=\mu_{q,p}$ whenever $r\in[r_{q,p},r_{q,p+1}[$ for certain $p\in[1,l(q)].$
    \begin{itemize}
      \item If we replace $i,k\in [r_{q-1},r_{q}[$ (with $i<k$) in equation \ref{Eq1-theo-Key-condition}, we obtain
          \begin{equation}\label{eq7-recalling-Key}
            2\gamma_{ir}\gamma_{kr}=
            \sum_{j\in J_r^{(q)}}{t_{rj}\left(\gamma_{kj}\gamma_{kr}s_{ki}+
            \gamma_{kj}\gamma_{ir}+\gamma_{ij}\gamma_{kr}\right)},
          \end{equation}
          since by condition (1), we have $\gamma_{kj}=\gamma_{ij}=0$ for all $j<r_{q-1}.$
          Now note that $s_{ki}=0$ therefore equation \ref{eq7-recalling-Key} takes the form
          \begin{equation}\label{eq8-recalling-Key}
            2\gamma_{ir}\gamma_{kr}=a_{ri}\gamma_{kr}+a_{rk}\gamma_{ir}.
          \end{equation}
          where (again) $a_{ri}=t_{r,\rrn(i)}\gamma_{i,\rrn(i)}.$
      \item If we replace $i\in [r_{q-1},r_{q}[$ and $k\in [r_q,r_{q+1}[,$ we obtain
          \begin{equation}\label{eq9-recalling-Key}
            \left(2\gamma_{ir}+s_{ki}\gamma_{kr}-a_{ri}\right)\gamma_{kr}=0.
          \end{equation}
    \end{itemize}
    If we assume that $\gamma_{kr}\neq 0$ for $k\in[r_q,r_{q+1}[$ and $r\in[r_{q,p},r_{q,p+1}[$ for certain $p\in[1,l(q)].$

    If we have that $|I_{k}^{(q)}|\geq 2,$ we can follow the same arguments that we use in the case $q=1$ to obtain $|I_{k}^{(q)}|=|J_{r}^{(q)}|.$ Otherwise, if $|I_{k}^{(p)}|=1$ and we assume that $|J_{r}^{(q)}|>1,$ then we can pick $j_1,j_2\in J_{r}^{(q)}$ with $j_1\neq j_2.$ Let $i_1:=\rrn^{-1}(j_1)$ and $i_2:=\rrn^{-1}(j_2),$ then $i_1,i_2\in [r_{q-1},r_q[$ and $i_1\neq i_2.$ By our hypothesis we have $s_{ki_1}s_{ki_2}=0.$

    If we assume that $s_{ki_2}=0,$ then by equation \ref{eq9-recalling-Key}, we have $2\gamma_{i_2r}=a_{r,i_2}.$ Since $\rrn(i_2)\in J_{r}^{(q)},$ we have that $a_{r,i_2}\neq 0,$ then $\gamma_{i_2r}\neq0.$ If we replace $i=i_1$ and $k=i_2$ in equation \ref{eq8-recalling-Key} we obtain $a_{r,i_1}\gamma_{i_2r}=0,$ therefore $a_{r,i_1}=0,$ but that means that $t_{r,\rrn(i_1)}=0,$ and this is a contradiction with our assumption that $j_1=\rrn(i_1)\in J_{r}^{(q)}.$ So we conclude that $|J_{r}^{(q)}|=1.$

    The above analysis says that having $\gamma_{kr}\neq 0$ implies that $k,r$ both belongs to the same subinterval $[r_{q,p},r_{q,p+1}[,$ and therefore $\rrn([r_{q,p},r_{q,p+1}[)=[r_{q,p},r_{q,p+1}[$ as desired.
    \end{itemize}
\end{proof}

\begin{example}\label{ex1-prop-Gamma-triang2}
  Consider the matrices
  \begin{equation}\label{eq1-ex1-prop-Gamma-triang2}
  T=\left[\begin{matrix}
              0 & 0 & 0& 0& 0& 0 \\
              0 & 0 & 0& 0& 0& 0 \\
              t_{31} & t_{32} & 0& 0 & 0& 0 \\
              t_{41} & t_{42} & 0& 0 & 0& 0 \\
              0 & 0 & t_{53}& {0} & 0& 0 \\
              0 & 0 & t_{63}& t_{64} & 0& 0
            \end{matrix}\right],\quad
  S=\left[\begin{matrix}
              0 & 0 & 0& 0& 0& 0 \\
              0 & 0 & 0& 0& 0& 0 \\
              s_{31} & s_{32} & 0& 0 & 0& 0 \\
              s_{41} & s_{42} & 0& 0 & 0& 0 \\
              0 & 0 & s_{53}& {0} & 0& 0 \\
              0 & 0 & s_{63}& s_{64} & 0& 0
            \end{matrix}\right]
 \end{equation}
 And assume that
 \begin{equation*}
   \mathbb{W}(T)=\mathbb{W}(S)=(3,5),\quad \mathbb{M}(T)=\mathbb{M}(S)=
   \left(\left[\begin{matrix}
           3 \\
           2
         \end{matrix}\right],\left[\begin{matrix}
                               5 & 6 \\
                               1 & 2
                             \end{matrix}\right]\right).
 \end{equation*}
 If $\bgamma:\calA(T)\rightarrow\calA(S)$ is an isomorphism, then Proposition \ref{prop-Gamma-triang-2} says that the associated matrix $\Gamma=[\gamma_{ij}]$ and the corresponding bijection $\rrn$ have one of the following forms:
 \begin{equation*}
   \Gamma_1=\left[\begin{matrix}
           \gamma_{11} & 0 & \gamma_{13} & \gamma_{14} & \gamma_{15} & \gamma_{16} \\
           0 & \gamma_{22} & \gamma_{23} & \gamma_{24} & \gamma_{25} & \gamma_{26} \\
           0 & 0 & \gamma_{33} & 0 & \gamma_{35} & \gamma_{36} \\
           0 & 0 & 0 & \gamma_{44} & \gamma_{45} & \gamma_{46} \\
           0 & 0 & 0 & 0 & \gamma_{55} & 0 \\
           0 & 0 & 0 & 0 & 0 & \gamma_{66} \\
         \end{matrix}\right],\quad
   \rrn:\left\{\begin{matrix}
                 1\mapsto 1 \\
                 2\mapsto 2 \\
                 3\mapsto 3 \\
                 4\mapsto 4 \\
                 5\mapsto 5 \\
                 6\mapsto 6
               \end{matrix}\right.
 \end{equation*}
 \begin{equation*}
   \Gamma_2=\left[\begin{matrix}
           0 & \gamma_{22} & \gamma_{13} & \gamma_{14} & \gamma_{15} & \gamma_{16} \\
           \gamma_{21} & 0 & \gamma_{23} & \gamma_{24} & \gamma_{25} & \gamma_{26} \\
           0 & 0 & \gamma_{33} & 0 & \gamma_{35} & \gamma_{36} \\
           0 & 0 & 0 & \gamma_{44} & \gamma_{45} & \gamma_{46} \\
           0 & 0 & 0 & 0 & \gamma_{55} & 0 \\
           0 & 0 & 0 & 0 & 0 & \gamma_{66} \\
         \end{matrix}\right],\quad
   \rrn:\left\{\begin{matrix}
                 1\mapsto 2 \\
                 2\mapsto 1 \\
                 3\mapsto 3 \\
                 4\mapsto 4 \\
                 5\mapsto 5 \\
                 6\mapsto 6
               \end{matrix}\right.
 \end{equation*}
 \begin{equation*}
   \Gamma_3=\left[\begin{matrix}
           \gamma_{11} & 0 & \gamma_{13} & \gamma_{14} & \gamma_{15} & \gamma_{16} \\
           0 & \gamma_{22} & \gamma_{23} & \gamma_{24} & \gamma_{25} & \gamma_{26} \\
           0 & 0 & 0 & \gamma_{34} & \gamma_{35} & \gamma_{36} \\
           0 & 0 & \gamma_{43} & 0 & \gamma_{45} & \gamma_{46} \\
           0 & 0 & 0 & 0 & \gamma_{55} & 0 \\
           0 & 0 & 0 & 0 & 0 & \gamma_{66} \\
         \end{matrix}\right],\quad
   \rrn:\left\{\begin{matrix}
                 1\mapsto 1 \\
                 2\mapsto 2 \\
                 3\mapsto 4 \\
                 4\mapsto 3 \\
                 5\mapsto 5 \\
                 6\mapsto 6
               \end{matrix}\right.
 \end{equation*}
 \begin{equation*}
   \Gamma_4=\left[\begin{matrix}
           0 & \gamma_{12} & \gamma_{13} & \gamma_{14} & \gamma_{15} & \gamma_{16} \\
           \gamma_{21} & 0 & \gamma_{23} & \gamma_{24} & \gamma_{25} & \gamma_{26} \\
           0 & 0 & 0 & \gamma_{34} & \gamma_{35} & \gamma_{36} \\
           0 & 0 & \gamma_{43} & 0 & \gamma_{45} & \gamma_{46} \\
           0 & 0 & 0 & 0 & \gamma_{55} & 0 \\
           0 & 0 & 0 & 0 & 0 & \gamma_{66} \\
         \end{matrix}\right],\quad
   \rrn:\left\{\begin{matrix}
                 1\mapsto 2 \\
                 2\mapsto 1 \\
                 3\mapsto 4 \\
                 4\mapsto 3 \\
                 5\mapsto 5 \\
                 6\mapsto 6
               \end{matrix}\right.
 \end{equation*}
 Note that while Proposition \ref{prop-Gamma-triang-1} give us $8$ possible forms for the matrix $\Gamma,$ when we incorporate the information of measures using Proposition \ref{prop-Gamma-triang-2}, we reduce those possibilities to just $4$ of them.
\end{example}

\begin{definition}\label{def-sur-les-sets-I-J}
  Let $T,S$ be two nonzero \nSLTM, such that:
   \begin{enumerate}
     \item $T$ and $S$ are in the $2-$REF,
     \item $\mathbb{W}(T)=\mathbb{W}(S)=(r_1,\dots,r_h).$
     \item $\mathbb{M}(T)=\mathbb{M}(S)=(M_1,\dots,M_h).$
   \end{enumerate}
   For each $q\in[1,h]$ and each $r\in [r_q,r_{q+1}[$  define the sets:
   \begin{equation*}
     J_{r}^{(q)}=\{j\in[r_{q-1},r_{q}[\quad :t_{rj}\neq 0\},\quad
     I_{r}^{(q)}=\{j\in[r_{q-1},r_{q}[\quad:s_{rj}\neq 0\}
   \end{equation*}
\end{definition}

\begin{theorem}\label{coro-sur-les-sets-I-J}$((I,J)-\textrm{\textbf{Test}})$
   Let $T,S$ be two nonzero \nSLTM, satisfying conditions (1),(2),(3) of Definition \ref{def-sur-les-sets-I-J}. If $\bgamma:\calA(T)\rightarrow\calA(S)$ is an isomorphism and $\rrn:[1,n]\rightarrow [1,n]$ is a bijection satisfying all the conditions described in Proposition \ref{prop-Gamma-triang-2} (relative to the associated matrix $\Gamma=[\gamma_{ij}]),$ then we have that
    \begin{equation}\label{eq1-coro-sur-les-sets}
      \rrn(I_{k}^{(q)})=J_{\rrn(k)}^{(q)},\quad\textrm{for each}\quad k\in[r_q,r_{q+1}[.
    \end{equation}
\end{theorem}
\begin{proof}
  It follows from the proof of Proposition \ref{prop-Gamma-triang-2}.
\end{proof}

Theorem \ref{coro-sur-les-sets-I-J} provide a considerably reduction of options for possible isomorphisms between the algebras $\calA(T)$ and $\calA(S):$
\begin{example}
  Consider the following matrices
  \begin{equation*}
    T=\left[\begin{matrix}
              0 & 0 & 0 & 0 & 0 & 0 \\
              0 & 0 & 0 & 0 & 0 & 0 \\
              0 & 0 & 0 & 0 & 0 & 0 \\
              t_{41} & t_{42} & 0 & 0 & 0 & 0 \\
              t_{51} & t_{52} & 0 & 0 & 0 & 0 \\
              t_{61} & t_{62} & t_{63} & 0 & 0 & 0
            \end{matrix}\right],\quad
  S=\left[\begin{matrix}
              0 & 0 & 0 & 0 & 0 & 0 \\
              0 & 0 & 0 & 0 & 0 & 0 \\
              0 & 0 & 0 & 0 & 0 & 0 \\
              s_{41} & 0 & s_{43} & 0 & 0 & 0 \\
              s_{51} & 0 & s_{43} & 0 & 0 & 0 \\
              s_{61} & s_{62} & s_{63} & 0 & 0 & 0
            \end{matrix}\right],
  \end{equation*}
  and assume that
  \begin{equation*}
   \mathbb{W}(T)=\mathbb{W}(S)=(4),\quad \textrm{and}\quad  M_1(T)=M_1(S)=\left[\begin{matrix}
                   4 & 6 \\
                   2 & 3
                 \end{matrix}\right]
  \end{equation*}
  Following Definition \ref{def-sur-les-sets-I-J}, we have:
  \begin{equation*}
   \begin{array}{ccc}
     J_{4}^{(1)}=\{1,2\} & J_{5}^{(1)}=\{1,2\} & J_{6}^{(1)}=\{1,2,3\} \\
     \quad & \quad & \quad\\
     I_{4}^{(1)}=\{1,3\} & I_{5}^{(1)}=\{1,3\} & I_{6}^{(1)}=\{1,2,3\}.
   \end{array}
  \end{equation*}
  If  $\bgamma:\calA(T)\rightarrow \calA(S),$ is an isomorphism, then the bijection $\rrn:[1,n]\rightarrow [1,n]$ of Proposition \ref{prop-Gamma-triang-2}, satisfies the following rules:
  \begin{equation*}
    \rrn([1,3])=[1,3],\quad \rrn([4,5])=[4,5],\quad \rrn(6)=6.
  \end{equation*}
  Moreover by Theorem \ref{coro-sur-les-sets-I-J}, we have a more specific requirement:
  \begin{equation*}
    \rrn(\{1,3\})=\{1,2\},\quad \rrn(2)=3
  \end{equation*}
  Therefore the associated matrix $\Gamma$ and the bijection $\rrn$ have one of the following forms:
  \begin{equation*}
    \Gamma_1=\left[\begin{matrix}
                   \gamma_{11} & 0 & 0 & \gamma_{14} & \gamma_{15} & \gamma_{16} \\
                   0 & 0 & \gamma_{23} & \gamma_{24} & \gamma_{25} & \gamma_{26} \\
                   0 & \gamma_{32} & 0 & \gamma_{34} & \gamma_{35} & \gamma_{36} \\
                   0 & 0 & 0 & \gamma_{44} & 0 & 0 \\
                   0 & 0 & 0 & 0 & \gamma_{55} & 0 \\
                   0 & 0  & 0 & 0 & 0 & \gamma_{66} \\
                 \end{matrix}\right],\quad
    \rrn:\left\{\begin{array}{c}
                  1\mapsto 1 \\
                  2\mapsto 3 \\
                  3\mapsto 2 \\
                   4\mapsto 4 \\
                   5\mapsto 5 \\
                  6\mapsto 6
                \end{array}\right.
  \end{equation*}
  \begin{equation*}
    \Gamma_2=\left[\begin{matrix}
                   \gamma_{11} & 0 & 0 & \gamma_{14} & \gamma_{15} & \gamma_{16} \\
                   0 & 0 & \gamma_{23} & \gamma_{24} & \gamma_{25} & \gamma_{26} \\
                   0 & \gamma_{32} & 0 & \gamma_{34} & \gamma_{35} & \gamma_{36} \\
                   0 & 0 & 0 & 0 & \gamma_{45} & 0 \\
                   0 & 0 & 0 & \gamma_{54} & 0 & 0 \\
                   0 & 0  & 0 & 0 & 0 & \gamma_{66} \\
                 \end{matrix}\right],\quad
    \rrn:\left\{\begin{array}{c}
                  1\mapsto 1 \\
                  2\mapsto 3 \\
                  3\mapsto 2 \\
                   4\mapsto 5 \\
                   5\mapsto 4 \\
                  6\mapsto 6
                \end{array}\right.
  \end{equation*}
  \begin{equation*}
    \Gamma_3=\left[\begin{matrix}
                   0 & \gamma_{12} & 0 & \gamma_{14} & \gamma_{15} & \gamma_{16} \\
                   0 & 0 & \gamma_{23} & \gamma_{24} & \gamma_{25} & \gamma_{26} \\
                   \gamma_{31} & 0 & 0 & \gamma_{34} & \gamma_{35} & \gamma_{36} \\
                   0 & 0 & 0 & 0 & \gamma_{45} & 0 \\
                   0 & 0 & 0 & \gamma_{54} & 0 & 0 \\
                   0 & 0  & 0 & 0 & 0 & \gamma_{66} \\
                 \end{matrix}\right],\quad
    \rrn:\left\{\begin{array}{c}
                  1\mapsto 2 \\
                  2\mapsto 3 \\
                  3\mapsto 1 \\
                   4\mapsto 5 \\
                   5\mapsto 4 \\
                  6\mapsto 6
                \end{array}\right.
  \end{equation*}
  \begin{equation*}
    \Gamma_4=\left[\begin{matrix}
                   0 & \gamma_{12} & 0 & \gamma_{14} & \gamma_{15} & \gamma_{16} \\
                   0 & 0 & \gamma_{23} & \gamma_{24} & \gamma_{25} & \gamma_{26} \\
                   \gamma_{31} & 0 & 0 & \gamma_{34} & \gamma_{35} & \gamma_{36} \\
                   0 & 0 & 0 & \gamma_{44} & 0 & 0 \\
                   0 & 0 & 0 & 0 & \gamma_{55} & 0 \\
                   0 & 0  & 0 & 0 & 0 & \gamma_{66} \\
                 \end{matrix}\right],\quad
    \rrn:\left\{\begin{array}{c}
                  1\mapsto 2 \\
                  2\mapsto 3 \\
                  3\mapsto 1 \\
                   4\mapsto 4 \\
                   5\mapsto 5 \\
                  6\mapsto 6
                \end{array}\right.
  \end{equation*}
\end{example}

Theorem \ref{coro-sur-les-sets-I-J} can be also used to prove that two matrices $S,T$ are not related under relation $\sim,$ even when $\mathbb{W}(T)=\mathbb{W}(S)$ and $\mathbb{M}(T)=\mathbb{M}(S).$ If we are able to prove that there is no  bijection $\rrn:[1,n]\rightarrow [1,n]$ satisfying item (2) of Proposition \ref{prop-Gamma-triang-2} and equation \ref{eq1-coro-sur-les-sets}, then we can conclude that $S\nsim T.$
\begin{example}
  Consider the following matrices
  \begin{equation*}
    T=\left[\begin{matrix}
              0 & 0 & 0 & 0 & 0 & 0 \\
              0 & 0 & 0 & 0 & 0 & 0 \\
              0 & 0 & 0 & 0 & 0 & 0 \\
              1 & 1 & 0 & 0 & 0 & 0 \\
              1 & 0 & 1 & 0 & 0 & 0 \\
              1 & 1 & 1 & 0 & 0 & 0
            \end{matrix}\right],\quad
  S=\left[\begin{matrix}
              0 & 0 & 0 & 0 & 0 & 0 \\
              0 & 0 & 0 & 0 & 0 & 0 \\
              0 & 0 & 0 & 0 & 0 & 0 \\
              1 & 0 & 1 & 0 & 0 & 0 \\
              1 & 0 & 1 & 0 & 0 & 0 \\
              1 & 1 & 1 & 0 & 0 & 0
            \end{matrix}\right],
  \end{equation*}
  Note that
  \begin{equation*}
   \mathbb{W}(T)=\mathbb{W}(S)=(4),\quad \textrm{and}\quad  M_1(T)=M_1(S)=\left[\begin{matrix}
                   4 & 6 \\
                   2 & 3
                 \end{matrix}\right]
  \end{equation*}
  Following Definition \ref{def-sur-les-sets-I-J}, we have:
  \begin{equation*}
   \begin{array}{ccc}
     J_{4}^{(1)}=\{1,2\} & J_{5}^{(1)}=\{1,3\} & J_{6}^{(1)}=\{1,2,3\} \\
     \quad & \quad & \quad\\
     I_{4}^{(1)}=\{1,3\} & I_{5}^{(1)}=\{1,3\} & I_{6}^{(1)}=\{1,2,3\}.
   \end{array}
  \end{equation*}
  If we suppose that $T\sim S,$ that is, there is an isomorphism $\bgamma:\calA(T)\rightarrow \calA(S),$ and a bijection $\rrn:[1,n]\rightarrow [1,n]$ satisfying item (2) of Proposition \ref{prop-Gamma-triang-2}, that is:
  \begin{equation*}
    \rrn([1,3])=[1,3],\quad \rrn([4,5])=[4,5],\quad \rrn(6)=6.
  \end{equation*}
  Since $\{1,3\}= I_4\cap I_5\cap I_6,$ Theorem \ref{coro-sur-les-sets-I-J} implies that $\{\rrn(1),\rrn(3)\}= J_4\cap J_5\cap J_6=\{1\},$ but this contradicts the fact that $\rrn$ is a bijection. Therefore there is no such  bijection $\rrn$ and we conclude that $S\nsim T.$
\end{example}

\begin{lemma}\label{lemma-M1T-neq-M1S}
  Let $S,T$ be two \nSLTM, such that both are in $2-$REF and $\mathbb{W}(S)=\mathbb{W}(T)=(r_1,\dots,r_h).$ If $M_1(S)\neq M_1(T),$ then $S\nsim T.$
\end{lemma}
\begin{proof}
If we suppose that there is an isomorphism $\bgamma:\calA(T)\rightarrow \calA(S),$ then by lemma \ref{lemma-from-1REF-to-2REF} the associated matrix $\Gamma=[\gamma_{ij}]$ has the form given in equation \ref{eq0-Gamma-triangular} and there is a bijection $\rrn$ satisfying conditions (2.a) and (2.b) of Proposition \ref{prop-Gamma-triang-1}.
  \begin{description}
    \item[Case 1]
      If $\mu\left(S_{[1,r_1[}^{[r_1]}\right)\neq \mu\left(T_{[1,r_1[}^{[r_1]}\right).$ We can assume without loss of generality that $\mu\left(S_{[1,r_1[}^{[r_1]}\right)> \mu\left(T_{[1,r_1[}^{[r_1]}\right).$
        Let $k_1:=\rrn^{-1}(r_1),$ then by Lemma \ref{lemma-from-1REF-to-2REF}, we have that $k_1\in [r_1,r_2[.$ Following the same analysis and same notations as  we did in the proof of Proposition \ref{prop-Gamma-triang-2}, we obtain the following relations:
        \begin{equation}\label{eq1-proof-lemma-M1T-neq-M1S}
          \left\{\begin{array}{c}
                   2\gamma_{ir_1}\gamma_{jr_1}=a_{r_1,i}\gamma_{jr_1}
                   +a_{r_1,j}\gamma_{ir_1}\\
                   \quad\\
                   2\gamma_{ir_1}=a_{r_1,i}-s_{k_1i} \\
                   \quad\\
                   2\gamma_{jr_1}=a_{r_1,j}-s_{k_1j}
                 \end{array}  \right.
        \end{equation}
        for each $i,j\in[1,r_1[$ such that $i,j\in I_{k_1}^{(1)}$ (By Lemma  \ref{lemma-r1-geq-3}, we always can find at least two of them). And this implies that
        \begin{equation}\label{eq2-proof-lemma-M1T-neq-M1S}
          s_{k_1i}s_{k_1j}\gamma_{k_1r_1}^2=a_{r_1,i}a_{r_1,j}.
        \end{equation}
        Since the left side of equation \ref{eq2-proof-lemma-M1T-neq-M1S} is different from $0,$ it implies that $a_{r_1i}=t_{r_1\rrn(i)}\gamma_{i\rrn(i)}\neq0$ for each $i\in I_{k_1}^{(1)}.$ In other words, $\rrn(I_{k_1}^{(1)})\subset J_{r_1}^{(1)}$ and since $\rrn$ is injective we have $|I_{k_1}^{(1)}|\leq |J_{r_1}^{(1)}|,$ but this is a contradiction since by our hypothesis, we have that  $|I_{r_1}^{(1)}|> |J_{r_1}^{(1)}|,$ and since $S$ is in the $2-$REF, we have that $|I_{k}^{(1)}|\geq |I_{r_1}^{(1)}|,$ for each $k\in [r_1,r_2[.$
    \item[Case 2] If $\mu\left(S_{[1,r_1[}^{[r_1]}\right)= \mu\left(T_{[1,r_1[}^{[r_1]}\right),$ define $$l_1:=\min\left\{l\in]r_1,r_2[\quad:\quad \mu\left(S_{[1,r_1[}^{[l]}\right)\neq \mu\left(T_{[1,r_1[}^{[l]}\right)\right\}.$$
     Since $M_1(S)\neq M_1(T),$ $l_1$ is well defined.  We can assume without loss of generality that $\mu\left(S_{[1,r_1[}^{[l_1]}\right)> \mu\left(T_{[1,r_1[}^{[l_1]}\right).$
     For each $r\in [r_1,l_1[,$ let $k:=\rrn^{-1}(r).$ Following the same arguments used in  the proof of Proposition \ref{prop-Gamma-triang-2}, we can conclude that $|I_k^{(1)}|=|J_{r}^{(1)}|$ and this implies that $k\in [r_1,l_1[.$ Therefore if we define $k_1:=\rrn^{-1}(l_1),$ then necessarily $k_1\in [l_1,r_2[$ and then $|I_{k_1}^{(1)}|\geq |I_{l_1}^{(1)}|>|J_{l_1}^{(1)}|.$
     If we use equation \ref{Eq1-theo-Key-condition}  we obtain the following relations
     \begin{equation}\label{eq3-proof-lemma-M1T-neq-M1S}
          \left\{\begin{array}{c}
                   2\gamma_{il_1}\gamma_{jl_1}=a_{l_1,i}\gamma_{jl_1}
                   +a_{l_1,j}\gamma_{il_1}\\
                   \quad\\
                   2\gamma_{il_1}=a_{l_1,i}-s_{k_1i} \\
                   \quad\\
                   2\gamma_{jl_1}=a_{l_1,j}-s_{k_1j}
                 \end{array}  \right.
        \end{equation}
        for each $i,j\in[1,r_1[$ such that $i,j\in I_{k_1}^{(1)}.$ Following the same arguments we used in case  1, we conclude that $|I_{k_1}^{(1)}|\leq |J_{l_1}^{(1)}|,$ obtaining a contradiction just like in case 1.
  \end{description}
  Cases 1 and 2, imply that there is no isomorphism $\bgamma:\calA(T)\rightarrow \calA(S),$ and therefore $T\nsim S$ as desired.
\end{proof}

\begin{lemma}\label{lemma-MqT-neq-MqS}
  Let $S,T$ be two \nSLTM, such that both are in $2-$REF and $\mathbb{W}(S)=\mathbb{W}(T)=(r_1,\dots,r_h).$ If there is a $q\in]1,h]$ such that $M_q(S)\neq M_q(T),$ then $S\nsim T.$
\end{lemma}

 The proof of Lemma \ref{lemma-MqT-neq-MqS} is analogous to the proof of Lemma \ref{lemma-M1T-neq-M1S}, we only have to pay attention in some few details:

\begin{proof}
  If $\bgamma:\calA(T)\rightarrow \calA(S),$ is an isomorphism, Lemma \ref{lemma-from-1REF-to-2REF} implies that the associated matrix $\Gamma$ has the form given in equation \ref{eq0-Gamma-triangular} and there is a bijection $\rrn:[1,n]\rightarrow [1,n]$ satisfying conditions (2.a) and (2.b) of Proposition \ref{prop-Gamma-triang-1}.

  Following the proof of Proposition \ref{prop-Gamma-triang-2}, we can see that for all $p\in [1,q-1[$ and all $k\in[r_{p},r_{p+1}[$ we have that $\rrn(I_{k}^{(p)})=J_{\rrn(k)}^{(p)},$ in particular $|I_{k}^{(p)}|=|J_{\rrn(k)}^{(p)}|.$

  \begin{description}
    \item[Case 1] If $\mu\left(S_{[r_{q-1},r_q[}^{[r_q]}\right)\neq \mu\left(T_{[r_{q-1},r_q[}^{[r_q]}\right),$ we can assume without loss of generality that $\mu\left(S_{[r_{q-1},r_q[}^{[r_q]}\right)> \mu\left(T_{[r_{q-1},r_q[}^{[r_q]}\right).$ Let $k_1:=\rrn^{-1}(r_q),$ then $k_1\in[r_q,r_{q+1}[.$ Using equation \ref{Eq1-theo-Key-condition}, we obtain
        \begin{equation}\label{eq1-proof-lemma-MqT-neq-MqS}
          \left\{\begin{array}{c}
                   2\gamma_{ir_q}\gamma_{jr_q}=a_{r_q,i}\gamma_{jr_q}
                   +a_{r_q,j}\gamma_{ir_q}\\
                   \quad\\
                   2\gamma_{ir_q}=a_{r_q,i}-s_{k_1i} \\
                   \quad\\
                   2\gamma_{jr_q}=a_{l_q,j}-s_{k_1j}
                 \end{array}  \right.
        \end{equation}
        for all $i,j\in [r_{q-1},r_{q}[.$ In particular, if we take $i,j\in I_{k_1}^{(q)}.$ Note that, since we are assuming that $\mu\left(S_{[r_{q-1},r_q[}^{[r_q]}\right)> \mu\left(T_{[r_{q-1},r_q[}^{[r_q]}\right)$ we have that $|I_{k_1}^{(q)}|\geq 2.$ Applying the same analysis that we used in the proof of Lemma \ref{lemma-M1T-neq-M1S}, we conclude that $\rrn(I_{k_1}^{(q)})\subset J_{r_q}^{(q)}$ and we obtain a contradiction.
    \item[Case 2] If $\mu\left(S_{[r_{q-1},r_q[}^{[r_q]}\right)= \mu\left(T_{[r_{q-1},r_q[}^{[r_q]}\right),$ define
        $$l_1=\min\left\{l\in [r_q,r_{q+1}[\quad:\quad \mu\left(S_{[r_{q-1},r_q[}^{[l]}\right)\neq \mu\left(T_{[r_{q-1},r_q[}^{[l]}\right)\right\}$$
        By our hypothesis $M_q(S)\neq M_{q}(T),$ therefore $l_1$ is well defined. We can assume without loss of generality that $\mu\left(S_{[r_{q-1},r_q[}^{[l_1]}\right)> \mu\left(T_{[r_{q-1},r_q[}^{[l_1]}\right).$ Again if we follow the analysis given in the proof of Proposition \ref{prop-Gamma-triang-2}, we conclude that $|I_{k}^{(q)}|=|J_{\rrn(k)}^{(q)}|$ for all $k\in [r_q,l_1[,$ therefore if we define $k_1:=\rrn^{-1}(l_1),$ we have that $k_1\in[l_1,r_{q+1}[.$ Since $S$ is in the $2-$REF, we have that $|I_{k_1}^{(q)}|\geq |I_{l_1}^{(q)}|>|J_{l_1}^{(q)}|.$ Now replacing $k=k_1,r=r_q$ in equation \ref{Eq1-theo-Key-condition}, we obtain
        \begin{equation}\label{eq2-proof-lemma-MqT-neq-MqS}
          \left\{\begin{array}{c}
                   2\gamma_{il_1}\gamma_{jl_1}=a_{l_1,i}\gamma_{jl_1}
                   +a_{l_1,j}\gamma_{il_1}\\
                   \quad\\
                   2\gamma_{il_1}=a_{l_1,i}-s_{k_1i} \\
                   \quad\\
                   2\gamma_{jl_1}=a_{l_1,j}-s_{k_1j}
                 \end{array}  \right.
        \end{equation}
        for each $i,j\in I_{k_1}^{(q)}$ (Note that our hypothesis imply that $|I_{k_1}^{(q)}|\geq 2).$ If we repeat the same analysis we did in case 2 of Lemma \ref{lemma-M1T-neq-M1S}, we obtain that $|I_{k_1}^{(q)}|\leq |J_{l_1}^{(q)}|$ and this is a contradiction.
  \end{description}
\end{proof}

\begin{theorem}\label{theo-invariants-2REF} $(\mathbb{M}-\textrm{\textbf{Test}})$
  Let $S,T$ two \nSLTM, such that both are in $2-$REF and $\mathbb{W}(S)=\mathbb{W}(T)=(r_1,\dots,r_h).$ If $S\sim T,$ then $\mathbb{M}(S)=\mathbb{M}(T).$
\end{theorem}
\begin{proof}
  Lemmas \ref{lemma-M1T-neq-M1S} and \ref{lemma-MqT-neq-MqS} implies that if $\mathbb{M}(T)\neq \mathbb{M}(S)$ then $T\nsim S.$
\end{proof}

In the following example, we show how we can use Theorem \ref{theo-invariants-2REF} to identify when two given matrices $T,S$ are not related under relation $\sim.$

\begin{example}\label{ex1-theo-ref2}
  Consider the matrices
  \begin{equation*}
  S=\left[\begin{matrix}
              0 & 0 & 0& 0& 0& 0 \\
              0 & 0 & 0& 0& 0& 0 \\
              {\color{blue}1} & {\color{blue}1} & 0& 0 & 0& 0 \\
              {\color{blue}1} & {\color{blue}1} & 0& 0 & 0& 0 \\
              s_{51} & s_{52} & {\color{red}1}& {\color{red}0} & 0& 0 \\
              s_{61} & s_{62} & {\color{red}0}& {\color{red}1} & 0& 0
            \end{matrix}\right],\quad
            T=\left[\begin{matrix}
              0 & 0 & 0& 0& 0& 0 \\
              0 & 0 & 0& 0& 0& 0 \\
              {\color{blue}1} & {\color{blue}1} & 0& 0 & 0& 0 \\
              {\color{blue}1} & {\color{blue}1} & 0& 0 & 0& 0 \\
              t_{51} & t_{52} & {\color{red}1}& {\color{red}0} & 0& 0 \\
              t_{61} & t_{62} & {\color{red}1}& {\color{red}1} & 0& 0
            \end{matrix}\right].
 \end{equation*}
 Note that both $S$ and $T$ are in the $2-$REF and $\mathbb{W}(S)=\mathbb{W}(T)=(3,5).$ Following Definition \ref{def-Blocks-2REF}, one can check that

 \begin{equation*}
   \mathbb{M}(S)=\left(\left[\begin{matrix}
                         3 \\
                         2
                       \end{matrix}\right],\left[\begin{matrix}
                         5 \\
                         1
                       \end{matrix}\right]\right),\quad
 \mathbb{M}(T)=\left(\left[\begin{matrix}
                         3 \\
                         2
                       \end{matrix}\right],
                       \left[\begin{matrix}
                         5 &6 \\
                         1 & 2
                       \end{matrix}\right]\right).
 \end{equation*}
 Since $\mathbb{M}(S)\neq \mathbb{M}(T),$  Theorem \ref{theo-invariants-2REF} asserts that $S\nsim T.$
\end{example}

\begin{corollary}\label{coro1-theo-invariants-2REF}
Let $T$ be a \SLTM. If $V_1,V_2$ are two matrices in $2-$REF such that $T\sim V_1$ and $T\sim V_2,$ then $\mathbb{M}(V_1)=\mathbb{M}(V_2).$
\end{corollary}
\begin{proof}
It follows directly from Theorem \ref{theo-invariants-2REF}.
\end{proof}

In view of Corollary \ref{coro1-theo-invariants-2REF}, we can extend the definition of Measure sequences to any  \SLTM:

\begin{definition}
Let $T$ be a \SLTM. We define the \emph{measure sequence} of $T$ by $\mathbb{M}(T)=\mathbb{M}(V),$ where $V$ is any matrix in $2-$REF, such that $T\sim V.$
\end{definition}

By Corollary \ref{coro1-theo-invariants-2REF}, the measure sequence of any \SLTM, is well defined.

\begin{corollary}\label{coro2-theo-invariants-2REF}
Let $T,S$ be two \nSLTM. If $T\sim S,$ then $\mathbb{M}(T)=\mathbb{M}(S).$
\end{corollary}
\begin{proof}
It follows directly from Corollary \ref{coro1-theo-invariants-2REF}.
\end{proof}


\section{On the number of isomorphism classes}


\subsection{A new lower bound}\label{ssec-lower-bound}

 Let $n,$ be a positive integer grater than $2.$ Lets denote by $N_n$ the number of different classes in \nSLTMn, under relation $\sim.$ In \cite{Lobos3} we show that (independently on $\F$),  $N_n\geq n-1.$ With the advances that we have made in this article, we are able to improve our lower bound, again independently on the ground field.

Recall that the number $C_n:\frac{1}{n+1}\binom{2n}{n}$ is known as the $n-$th \emph{Catalan number}. This numbers appear in several combinatorial situations (see \cite{Stanley2015}). For example, it is well known that the Catalan number $C_n$ is equal to the number of \emph{monotonic lattice paths} along the edges of a grid composed by $n\times n$ square cells, which do not pass above the diagonal. A monotonic path is one which starts in the upper left corner, finishes in the lower right corner, and consists entirely of edges pointing rightwards or downwards. The following are two examples of monotonic path for $n=5:$
  \begin{equation}\label{eq-paths-Catalan}
  \begin{tikzpicture}[xscale=0.5,yscale=0.5]
\node[] at (-1,2) {$P_0=$};
\draw[thin,gray] (0,0) --(0,4);
\draw[thin,gray] (1,0) --(1,4);
\draw[thin,gray] (2,0) --(2,4);
\draw[thin,gray] (3,0) --(3,4);
\draw[thin,gray] (4,0) --(4,4);
\draw[thin,gray] (0,0) --(4,0);
\draw[thin,gray] (0,1) --(4,1);
\draw[thin,gray] (0,2) --(4,2);
\draw[thin,gray] (0,3) --(4,3);
\draw[thin,gray] (0,4) --(4,4);
\draw[thin,red,dashed] (0,4) --(4,0);
\draw[thick,blue] (0,4) --(0,0);
\draw[thick,blue] (0,0) --(4,0);
\node[gray] at (4.2,4) {$2$};
\node[gray] at (4.2,3) {$3$};
\node[gray] at (4.2,2) {$4$};
\node[gray] at (4.2,1) {$5$};
\node[gray,below] at (0,0) {$1$};
\node[gray,below] at (1,0) {$2$};
\node[gray,below] at (2,0) {$3$};
\node[gray,below] at (3,0) {$4$};
\node[gray,below] at (4,0) {$5$};
\end{tikzpicture}\quad
    \begin{tikzpicture}[xscale=0.5,yscale=0.5]
\node[] at (-1,2) {$P_1=$};
\draw[thin,gray] (0,0) --(0,4);
\draw[thin,gray] (1,0) --(1,4);
\draw[thin,gray] (2,0) --(2,4);
\draw[thin,gray] (3,0) --(3,4);
\draw[thin,gray] (4,0) --(4,4);
\draw[thin,gray] (0,0) --(4,0);
\draw[thin,gray] (0,1) --(4,1);
\draw[thin,gray] (0,2) --(4,2);
\draw[thin,gray] (0,3) --(4,3);
\draw[thin,gray] (0,4) --(4,4);
\draw[thin,red,dashed] (0,4) --(4,0);
\draw[thick,blue] (0,4) --(0,3);
\draw[thick,blue] (0,3) --(1,3);
\draw[thick,blue] (1,3) --(1,0);
\draw[thick,blue] (1,0) --(4,0);
\node[gray] at (4.2,4) {$2$};
\node[gray] at (4.2,3) {$3$};
\node[gray] at (4.2,2) {$4$};
\node[gray] at (4.2,1) {$5$};
\node[gray,below] at (0,0) {$1$};
\node[gray,below] at (1,0) {$2$};
\node[gray,below] at (2,0) {$3$};
\node[gray,below] at (3,0) {$4$};
\node[gray,below] at (4,0) {$5$};
\end{tikzpicture}\quad
\begin{tikzpicture}[xscale=0.5,yscale=0.5]
\node[] at (-1,2) {$P_2=$};
\draw[thin,gray] (0,0) --(0,4);
\draw[thin,gray] (1,0) --(1,4);
\draw[thin,gray] (2,0) --(2,4);
\draw[thin,gray] (3,0) --(3,4);
\draw[thin,gray] (4,0) --(4,4);
\draw[thin,gray] (0,0) --(4,0);
\draw[thin,gray] (0,1) --(4,1);
\draw[thin,gray] (0,2) --(4,2);
\draw[thin,gray] (0,3) --(4,3);
\draw[thin,gray] (0,4) --(4,4);
\draw[thin,red,dashed] (0,4) --(4,0);
\draw[thick,blue] (0,4) --(0,2);
\draw[thick,blue] (0,2) --(2,2);
\draw[thick,blue] (2,2) --(2,1);
\draw[thick,blue] (2,1) --(3,1);
\draw[thick,blue] (3,1) --(3,0);
\draw[thick,blue] (3,0) --(4,0);
\node[gray] at (4.2,4) {$2$};
\node[gray] at (4.2,3) {$3$};
\node[gray] at (4.2,2) {$4$};
\node[gray] at (4.2,1) {$5$};
\node[gray,below] at (0,0) {$1$};
\node[gray,below] at (1,0) {$2$};
\node[gray,below] at (2,0) {$3$};
\node[gray,below] at (3,0) {$4$};
\node[gray,below] at (4,0) {$5$};
\end{tikzpicture}
  \end{equation}

  Note that in general a monotonic path $P$ determines a \emph{lower region} of the grid, that is, the path itself and the region that is below of the path .
 We shall denote by $P_0,$ the only path with only one change of direction (see equation \ref{eq-paths-Catalan}), we call $P_0$ the \emph{trivial path}. If we denote by $\calP_n$ the set of all monotonic paths in a grid compose by $n\times n$ square cells, then $|\calP_n|=C_n.$

\begin{theorem}\label{theo-new-lower-bound}
  There are at least $C_{n-1}$ different classes in \nSLTMn. That is $N_n\geq C_{n-1}.$
\end{theorem}

\begin{proof}
  We shall prove our claim, by defining an explicit function $\Psi:\calP_{n-1}\rightarrow \mathcal{TM}_n$ such that $\Psi(P_1)\nsim \Psi(P_2),$ whenever, $P_1,P_2$ are different paths.

  We first define for the trivial path $P_0$ the matrix $\Psi(P_0)=0_n.$ In order to define $\Psi(P)$ for nontrivial paths, we enumerate the vertical lines of the grid from $1$ to $n,$ from left to right, and the horizontal lines, from $2$ to $n$ from top to bottom, with the exception of the lowest line (see equation \ref{eq-paths-Catalan}). With this enumeration, each  node $(i,j)$ of the grid can be associated to an entry $s_{ij}$ of a $n\times n-$matrix $S=\Psi(P)$ (we insist that we are not considering nodes of the lowest horizontal line of the grid).
  Now we define the values of the entries $s_{ij}$ by the following rule:
  \begin{itemize}
    \item $s_{ij}=1,$ if $(i,j)$ is a node of the grid and it is part of the lower zone of $P.$
    \item $s_{ij}=0$ in any other case.
  \end{itemize}

  Note that, by construction, we always have that $s_{ij}=0$ for $i\in [1,2].$ Moreover it is clear that $\Psi(P)\in\mathcal{TM}_n$ and it is in both the $1-$REF and $2-$REF.

  Using the notation (with certain modifications) introduced in \cite{Juyumaya-Lobos-2025}, we can associate to any nontrivial path $P\in \calP_n$ an \emph{indexing matrix} $A(P)$ defined as follows:
  \begin{equation}\label{eq-indexing-matrix1}
    A(P)=\left[\begin{matrix}
                 r_1 & r_2 & \cdots & r_m \\
                 c_1 & c_2 & \cdots & c_m
               \end{matrix}\right],\quad
  \end{equation}
  where:
  \begin{itemize}
    \item  $3\leq r_1<r_2<\cdots <r_m\leq n;$ $2\leq c_1<c_2<\cdots <c_m<n;$ $c_j<r_j$ for $j\in [1,m].$
    \item Each column of $A(P)$ represents a node $(r_j,c_j)$ of the grid, where the path $P$ change its direction from rightwards to downward
  \end{itemize}

  On the other hand it is clear that a given indexing matrix $A,$ satisfying the above conditions, uniquely determines a path $P\in \calP_{n-1}.$ Therefore $A(P_1)\neq A(P_2)$ whenever $P_1\neq P_2.$

  We also define a vector $\mathbb{W}(P)=(w_1,\dots,w_h)$ associated to the matrix $A(P),$ given by the following rules:
  \begin{enumerate}
  \item Let $w_1=r_1.$
    \item If $m=1,$ then define $h=1$ and $\mathbb{W}(P):=(r_1).$
    \item If $m>1$ but $c_m<r_1$ then define $h=1$ and $\mathbb{W}(P):=(r_1).$ Otherwise, if $m>1$ and there is at least a $j\in[2,m]$ such that $c_j\geq r_1$ Take $i_2:=\min\{j\in]1,m]:c_j\geq r_1\}$ then define $w_2:=r_{i_2}.$
    \item Assuming that the index $i_p$ was already defined, and $w_p:=r_{i_p},$ then:
          \begin{itemize}
            \item If $c_j<w_p$ for all $j\in ]p,m]$ then define $h=p$ and $\mathbb{W}(P)=(w_1,\dots,w_p).$
            \item Otherwise, let $i_{p+1}:=\min\{j\in]i_p,m]:c_j\geq r_{i_p}\}$ and define $w_{p+1}:=r_{i_{p+1}}.$
          \end{itemize}
  \end{enumerate}
   Following recursively this process, we obtain $\mathbb{W}(P)=(w_1,\dots,w_h)$ as desired. (Note that the vector $\mathbb{W}(P)$ is totally determined by the indexing matrix $A(P).$)

  Having defined the indices $i_2,\dots,i_h$ as above, then for each $j\in[1,h]$ we define a matrix $M_j(P)$ associated to the path $P$ as follows:

  \begin{equation}\label{eq-measure-tables-for-paths0}
    M_j(P)=\left\{\begin{array}{cc}
                    \left[\begin{matrix}
                   r_1 & r_2 & \cdots & r_{i_2-1} \\
                   c_1 & c_2 & \cdots & c_{i_2-1}
                 \end{matrix}\right], & \textrm{if}\quad j=1\\
                    \quad & \quad \\
                    \left[\begin{matrix}
                   r_{i_j} & r_{i_j+1} & \cdots & r_{i_{j+1}-1} \\
                   c_{i_j}' & c_{i_j+1}' & \cdots & c_{i_{j+1}-1}'
                 \end{matrix}\right], & \textrm{if}\quad j\in]1,h[ \\
                    \quad & \quad \\
                    \left[\begin{matrix}
                   r_{i_h} & r_{i_h+1} & \cdots & r_m \\
                   c_{i_h}' & c_{i_h+1}' & \cdots & c_{m}'
                 \end{matrix}\right], & \textrm{if}\quad j=h
                  \end{array}\right.
  \end{equation}
   where $c_{p}':=c_p-r_{i_{j-1}}+1$ whenever $p\in[i_j,i_{j+1}[$ (assuming $i_{m+1}:=+\infty$)

  Comparing with definitions \ref{def-Blocks-1REF} and \ref{def-Blocks-2REF}, one can check that $\mathbb{W}(P)$ corresponds to the wall of the matrix $\Psi(P)$ and the sequence of matrices $\mathbb{M}(P):=(M_1(P),\dots,M_{h}(P))$ correspond to the measure sequence associated to the matrix $\Psi(P).$ By construction one can check that for two different paths $P_1,P_2\in\calP_{n-1}$ we have that one of the following assertions is true:
  \begin{itemize}
    \item $\mathbb{W}(P_1)\neq \mathbb{W}(P_2)$
    \item $\mathbb{W}(P_1)=\mathbb{W}(P_2)=(w_1,\dots,w_h)$ but there is a $j\in[1,h]$ such that $M_j(P_1)\neq M_{j}(P_2),$
  \end{itemize}
therefore by Theorems \ref{theo-invariants-1REF} and \ref{theo-invariants-2REF}, one can conclude that $\Psi(P_1)\nsim \Psi(P_2),$ whenever $P_1\neq P_2.$ This implies that $N_n\geq |\calP_{n-1}|=C_{n-1}$ as desired.
\end{proof}

\begin{example}\label{ex-proof-theo-new-lower-bound}
 Following definitions given in the proof of Theorem \ref{theo-new-lower-bound} for the paths in equation \ref{eq-paths-Catalan}, we have that the corresponding indexing matrices are:
  \begin{equation*}
    A(P_1)=\left[\begin{matrix}
                   3 \\
                   2
                 \end{matrix}\right],\quad
    A(P_2)=\left[\begin{matrix}
                   4&5 \\
                   3&4
                 \end{matrix}\right]
  \end{equation*}
  Therefore we have:
  \begin{equation*}
    \mathbb{W}(P_1)=(3),\quad  \mathbb{W}(P_2)=(4,5).
  \end{equation*}
  and
  \begin{equation*}
    \mathbb{M}(P_1)=\left(\left[\begin{matrix}
       3 \\
          2
        \end{matrix}\right]\right),\quad
        \mathbb{M}(P_2)=\left(\left[\begin{matrix}
                                                   4 \\
                                                   3
                                                 \end{matrix}\right],
                                                 \left[\begin{matrix}
                                                   5 \\
                                                   1
                                                 \end{matrix}\right]\right).
  \end{equation*}
  We can also construct the corresponding matrices:
  \begin{equation}\label{eq-paths-Catalan-matrix}
    \Psi({P_1})=\left[\begin{matrix}
                    0 & 0 & 0 & 0 & 0 \\
                    0 & 0 & 0 & 0 & 0 \\
                    1 & 1 & 0 & 0 & 0 \\
                    1 & 1 & 0 & 0 & 0 \\
                    1 & 1 & 0 & 0 & 0 \\
                  \end{matrix}\right],\quad
    \Psi({P_2})=\left[\begin{matrix}
                    0 & 0 & 0 & 0 & 0 \\
                    0 & 0 & 0 & 0 & 0 \\
                    0 & 0 & 0 & 0 & 0 \\
                    1 & 1 & 1 & 0 & 0 \\
                    1 & 1 & 1 & 1 & 0 \\
                  \end{matrix}\right]
  \end{equation}
  and we can verify that $\Psi(P_1)\nsim \Psi(P_2).$
\end{example}

Theorem \ref{theo-new-lower-bound} provide a beautiful number as lower bound for our number $N_n.$ Unfortunately, if we use all the invariants that we found in Section \ref{sec-Invariants}, we can verify that this lower bound is in general strict. This fact opens the problem of improving this lower bound, while we are not able to find the exact value of $N_n.$ In the following Proposition we prove that $N_n>C_{n-1}$ for each $n> 4:$



\begin{proposition}\label{prop-strict-lower-bound}
If $n>4$ then $N_n>C_{n-1}.$
\end{proposition}

\begin{proof}
  Let $T=[t_{rc}]$ be the \nSLTM, defined as follows:
  \begin{equation*}
    t_{rc}=\left\{\begin{array}{cc}
                    1 & \textrm{if}\quad r\in[4,n[\quad\textrm{and}\quad c\in\{1,2\} \\
                    1 & \textrm{if}\quad r=n\quad\textrm{and}\quad c\in\{1,3\} \\
                    0 & \textrm{otherwise}
                  \end{array}\right.
  \end{equation*}
  Is not difficult to check that $T$ is in the $2-\textrm{REF},$ its wall is composed for only one brick, whose rows have all the same measure $\mu=2.$ That is
  \begin{equation*}
    \mathbb{W}(T)=(4),\quad \mathbb{M}(T)=\left(\left[\begin{matrix}
                                                        4 \\
                                                        2
                                                      \end{matrix}\right]\right)
  \end{equation*}
  We claim that $T\nsim \Psi(P),$ for any path $P\in\calP_{n-1}.$ In fact, it is clear, by construction (see the proof of Theorem \ref{theo-new-lower-bound}) that $T\neq \Psi(P)$ for any $P\in\calP_{n-1}.$ On the other hand the only possible path $P_1\in \calP_{n-1}$ such that $\mathbb{W}(P_1)=\mathbb{W}(T)$ and $\mathbb{M}(P_1)=\mathbb{M}(T)$ is determined by the following indexing matrix:
  \begin{equation*}
    A(P_1)=\left[\begin{matrix}
                   4 \\
                   2
                 \end{matrix}\right]
  \end{equation*}
  and therefore the entries of the matrix $S=\Psi(P_1)$ are given as follows:
  \begin{equation*}
    s_{rc}=\left\{\begin{array}{cc}
                    1 & \textrm{if}\quad r\in[4,n]\quad\textrm{and}\quad c\in\{1,2\} \\
                    0 & \textrm{otherwise}
                  \end{array}\right.
  \end{equation*}
  If we apply Theorem \ref{coro-sur-les-sets-I-J} to the matrices $T$ and $S,$ we can see that for $r\in[4,n]$ we have $I_r^{(1)}=\{1,2\},$ while
  \begin{equation*}
    J_{r}^{(1)}=\left\{\begin{array}{cc}
                         \{1,2\} & \textrm{if}\quad r\in[4,n[ \\
                         \{1,3\} & \textrm{if}\quad r=n.
                       \end{array}\right.
  \end{equation*}
  If there were an isomorphism $\bgamma:\calA(T)\rightarrow \calA(S),$ Proposition \ref{prop-Gamma-triang-2} and Theorem \ref{coro-sur-les-sets-I-J}, would implied that there is a bijection $\rrn:[1,n]\rightarrow [1,n]$ such that $\rrn([1,3])=[1,3],$  $\rrn([4,n])=[4,n]$ and $\rrn(I_k^{(1)})=J_{\rrn(k)}^{(1)},$ for all  $k\in[1,n],$ but clearly this is impossible.

  We have proved that for $n>4,$ we can always find a class $[T]$ different from the $C_{n-1}$ classes $[\Psi(P)]$ described in the proof of Theorem \ref{theo-new-lower-bound}, therefore $N_n>C_{n-1}.$
\end{proof}

\begin{example}
  For $n=5$ the matrices $T$ and $S$ mentioned in the proof of Proposition \ref{prop-strict-lower-bound} are the following:
  \begin{equation*}
    T=\left[\begin{matrix}
              0 & 0 & 0 & 0 & 0 \\
              0 & 0 & 0 & 0 & 0 \\
              0 & 0 & 0 & 0 & 0 \\
              1 & 1 & 0 & 0 & 0 \\
              1 & 0 & 1 & 0 & 0
            \end{matrix}\right],\quad
  S=\left[\begin{matrix}
              0 & 0 & 0 & 0 & 0 \\
              0 & 0 & 0 & 0 & 0 \\
              0 & 0 & 0 & 0 & 0 \\
              1 & 1 & 0 & 0 & 0 \\
              1 & 1 & 0 & 0 & 0
            \end{matrix}\right]=\Psi(P_1)
  \end{equation*}
  where
  \begin{equation*}
    \begin{tikzpicture}[xscale=0.5,yscale=0.5]
\node[] at (-1,2) {$P_1=$};
\draw[thin,gray] (0,0) --(0,4);
\draw[thin,gray] (1,0) --(1,4);
\draw[thin,gray] (2,0) --(2,4);
\draw[thin,gray] (3,0) --(3,4);
\draw[thin,gray] (4,0) --(4,4);
\draw[thin,gray] (0,0) --(4,0);
\draw[thin,gray] (0,1) --(4,1);
\draw[thin,gray] (0,2) --(4,2);
\draw[thin,gray] (0,3) --(4,3);
\draw[thin,gray] (0,4) --(4,4);
\draw[thin,red,dashed] (0,4) --(4,0);
\draw[thick,blue] (0,4) --(0,2);
\draw[thick,blue] (0,2) --(1,2);
\draw[thick,blue] (1,2) --(1,0);
\draw[thick,blue] (1,0) --(4,0);
\node[gray] at (4.2,4) {$2$};
\node[gray] at (4.2,3) {$3$};
\node[gray] at (4.2,2) {$4$};
\node[gray] at (4.2,1) {$5$};
\node[gray,below] at (0,0) {$1$};
\node[gray,below] at (1,0) {$2$};
\node[gray,below] at (2,0) {$3$};
\node[gray,below] at (3,0) {$4$};
\node[gray,below] at (4,0) {$5$};
\end{tikzpicture}
  \end{equation*}
\end{example}

\textbf{Acknowledgments:}

\medskip
It is a pleasure to thank to Camilo Martinez Estay, for his collaboration on programming some of the algorithms developed in this article. I would especially want to thank Valentina Delgado and Rahma Salama, for their useful comments on the first draft of this paper.

\medskip
This work was supported by:
\begin{itemize}
  \item {FONDECYT de Postdoctorado 2024}, N°3240046, ANID, Chile.
  \item {Subvención a la Instalación en la Academia 2024}, N°85240053, ANID, Chile.
\end{itemize}




\sc
diego.lobosm@uv.cl, Universidad de Valpara\'iso, Chile.

\end{document}